\documentclass[twoside]{article}
\usepackage{latexsym}
\usepackage{pstricks}
\usepackage{amssymb}
\usepackage{amsmath}
\usepackage{amsfonts}
\usepackage{amsthm,array}
\usepackage{array}
\usepackage{epsfig}
\usepackage{graphicx}
\usepackage{float}
\usepackage{cite}

\def\noi{\noindent}
\def\vs{\vskip.4cm}

\newtheorem{remark-definition}{Remark and Definition}[section]
\newtheorem{rem-not}{Remark and Notation}[section]

\begin{document}
\title{\bf Higher-dimensional Golden Section, Generalized Fibonacci
Sequence, and Multiple Online Leasing}

 \date{}
\author{ Hu Maolin$^{1,*}$,~~~~ Luo Chu$^{2}$,~~~~ Xu Weidong$^{3,*}$\\
 {\small1~School of
 Mathematics and Statistics, Huaiyin Normal
 University,  Huaian, Jiangsu, 223300, China}\\
 {\small2~Great Bay University, Dongguan, Guangdong, 523000, China; luochu@gbu.edu}\\
 {\small3~School of Management, Zhejiang University, Hangzhou 310058, China}\\
{\small*~Correspondence: humaolin@hytc.edu.cn; weidxu@zju.edu.cn}\\}
 \maketitle
\vs
\begin{abstract}\noi
We discover a profound connection between the classical golden
section and the online leasing problem:
 When the number of skis to be rented or purchased increases
 from one to two, the optimal competitive algorithm is to rent two pairs of skis until
the golden section point $(\sqrt{5}-1)s/2$ within the time interval
[0,s], then buy one pair of skis at that point, and finally buy the
second pair of skis at time $s$, thereby achieving the optimal
competitive ratio of $(5+\sqrt{5})/4$.
  Building on this connection, we extend our research to
 explore higher-dimensional golden sections and multiple online leasing problem.

 For the generalization of the golden section, we provide definitions for the geometric-type and harmonic-type
 higher-dimensional
 golden sections, $\tau^{(m)}$  and $T^{(m)}$
 , respectively, based on the two equivalent propositions of the classical golden section, which divide the unit line
 segment into $m$ segments. We establish equations for the geometric-type higher-dimensional golden section ratio $\lambda$
 and the harmonic-type higher-dimensional golden section basis $\omega$. From the perspective of algebraic equations, we
 investigate the characteristics of the geometric-type higher-dimensional golden section ratio $\lambda$ and the
 harmonic-type
 higher-dimensional golden section basis $\omega$. Given the relationship between the classical golden section and the
 Fibonacci sequence, We extend the $k$-generalized Fibonacci sequence from existing literature to a bi-infinite
 $m$-generalized
 Fibonacci sequence and establish the relationship between two types of higher-dimensional golden sections and the
 bi-infinite $m$-generalized Fibonacci sequence. Furthermore, we introduce novel concepts including the generating
 sequence for harmonic-type higher-dimensional golden section points and the Fibonacci matrix, while
 investigating properties of submatrices of the Fibonacci matrix under row-column indexing systems.
 Results such as $det\mathbf{H}^{(m)}_{n}(x,y)=(-1)^{(n+h(xy))(m-1)}F^{(m)}_{x}F^{(m)}_{y}$ demonstrate the potential
 effectiveness of employing the Fibonacci matrix to study (generalized) Fibonacci sequences.

 Regarding the generalization of the online leasing problem, we consider the multiple online leasing of $m$ pairs of
 skis by a ski group. For the multiple online leasing problem, we design control and balancing strategy. Using
 the competitive analysis method for online problems, we prove that the control and balancing strategy is the
 optimal strategy for multiple online leasing. The results show that the optimal strategy for $m$-fold online
 leasing is precisely the $m$-dimensional harmonic-type golden section method, which means that the ski group rents skis
 and purchases one pair at each golden section point within the interval $[0,s]$. Considering that there
 may be no explicit solutions for the harmonic-type higher-dimensional golden section equations in dimensions higher than four, we
 propose and study the continuous relaxation problem for multiple online leasing. We explore the feasibility of
 using the relaxed solution corresponding to the optimal solution of the continuous relaxation problem
 to approximate the optimal solution of multiple
 online leasing or the largest real root of the harmonic-type higher-dimensional golden section equation. We correct the
 relaxation solution using the classical golden section method and propose the mantissa weighting method for
 solving the harmonic-type higher-dimensional golden section equation. Through error analysis, we show that the mantissa
 weighting method achieves high accuracy in solving the harmonic-type higher-dimensional golden section equation.
\end{abstract}
\vs \noi{\bf Keywords:} Golden section,
Higher-dimensional golden sections, Fibonacci sequence, Generalized
Fibonacci sequence, Generating sequence, Fibonacci matrix,
Full-dimensional submatrix, Online leasing, Multiple online leasing,
Control and balancing strategy, Competitive analysis, Competitive
ratio, Golden section-mantissa weighting method
 \medskip
\section{Introduction}

The golden section stands as one of the most classic problem in
geometry with a profound historical legacy, influencing diverse
fields spanning mathematics, art, and natural sciences. Whereas the
online leasing problem has emerged as a canonical model in
competitive algorithm research, driving three decades of
advancements in online decision theory. Despite their origins being
separated by over two millennia, we have discovered a fascinating
connection between these two classics. See Section 1.3 for details.
This discovery motivates our extension to higher-dimensional golden
sections, the bi-infinite $m$-generalized Fibonacci sequence, and
the multiple online leasing problem.

\subsection*{1.1\ \ A Survey of the Golden Section and the Fibonacci
Sequence}

\subsection*{1.1.1\ \ Principles and Applications of the Golden
Section}

The golden section, a quintessential concept imbued with profound
aesthetic significance in mathematical history, traces its
theoretical origins to the 5th-century BCE Pythagorean school in
ancient Greece. In Book VI of {\it Elements}(Heath trans. [13]),
Euclid systematically articulated the definition of {\bf extreme and
mean ratio}, marking the formal mathematical codification of this
theory:
 As shown in Figure 1, inserting point $C$ on a unit segment $AB$ such that
 $\frac{AB}{AC}=\frac{AC}{CB}$. Given $AB=AC+CB$,
 this yields the quadratic equation $AC^{2}+AC-1=0$. Its solution within $(0,1)$,
 $AC=\frac{\sqrt{5}-1}{2}\approx 0.618$,
 defines point $C$ as the golden section point of $AB$, with the ratio $\lambda=\omega=\frac{AB}{AC}
 =\frac{1+\sqrt{5}}{2}\approx1.618$
 representing the golden section ratio. This ratio $\lambda$ or $\omega$
 fundamentally reveals its essence as the
 largest real root of the functional equations
$\frac{1}{x}+\frac{1}{x^{2}}=1$ or $\frac{1}{x}+\frac{1}{x+1}=1$.
%%%%%%%%%%%%%%%%%%%%%%%%%%%%%%%%%%%%%%%%%%%%%%%%%%%%%%%%图象插入%%%%%%%%%%%%%%%%%%
\begin{figure}[H]
  \centering
  \includegraphics[width=\textwidth, trim=50 670 170 130, clip]{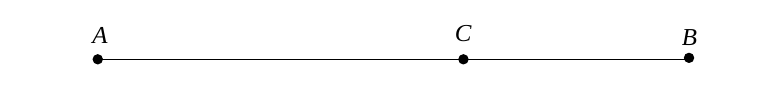}
  \caption{\small Classic Golden Section}
  \label{fig:fig1}
\end{figure}
%%%%%%%%%%%%%%%%%%%%%%%%%%%%%%%%%%%%%%%%%%%%%%%%%%%%%%%%%%%%%%%%%%%%%%%%%%%%%%%%%%%
The unique mathematical properties of the golden section ratio
establish it as a natural nexus between geometric intuition and
algebraic abstraction. As one of the earliest discovered irrational
proportions, it not only deepened humanity's understanding of
numerical essence (Huntley [26]) but also bridged discrete
mathematics and continuous geometry through the limit convergence
properties of Fibonacci sequences (Livio [38]). Since the Middle
Ages, interdisciplinary investigations have flourished: Fibonacci's
rabbit population model in {\it Liber Abaci} unveiled sequence
regularities, da Vinci's {\it Vitruvian Man} codified anthropometric
principles (Livio [38]), and Kepler's {\it Harmonices Mundi} (Kepler
[30]) constructed celestial motion models, collectively
 demonstrating the ratio's extraordinary potential in deciphering natural laws.

Modern mathematical advancements have extended the golden section's practical utility.
The 20th-century Golden Section Search (GSS) algorithm, leveraging
 its unique interval contraction mechanism, emerged as an efficient method for unimodal
 function optimization, reducing iterations by approximately $40\%$
  compared to traditional bisection (Kiefer [32], Tsai et al. [53], Cai et al. [5]).
  This optimization paradigm now permeates diverse fields: machine
  learning hyperparameter tuning (Koupaei et al. [33]), structural engineering design
  (Tang et al. [52]), renewable energy system control (Gayathri and Ezhilarasi
  [20], Paquianadin et al. [44]),
  and evolutionary algorithm enhancement
  (Cuevas et al. [8]). Current research focuses on integrating golden section
  ratio principles with stochastic optimization methods, such as augmenting particle
  swarm algorithms for global search efficiency (Talaei et al. [51]) or
  enhancing neural architecture search (NAS) for model compactness.

Cross-disciplinary studies further validate the golden section
ratio's universality. Architectural marvels from the Parthenon to Le
Corbusier's Modulor system (Le Corbusier [35]), Mondrian's
compositional principles and Apple's product design aesthetics
(Huntley [26]), biological phenomena including DNA helical turns and
phyllotaxis patterns (Livio [38], L$\rm\ddot{u}$ttge and Souza
[42]), all attest to its natural optimization characteristics.
Contemporary research employs fractal geometry and complexity
science to elucidate its role in nonlinear systems, self-organizing
structures, and information encoding, driving innovative
applications in cryptography and quantum computing.

The mathematical concept of the Golden Section, developed over two millennia, sustains an
ongoing interplay between theoretical depth and practical vitality.
 Its developmental trajectory, mirroring its proportional perfection, epitomizes the
 harmonious unification of rational thought and natural law.

\subsection*{1.1.2\ \ Evolution and Generalization of the
Fibonacci Sequence}

The Fibonacci sequence was first introduced by the 13th-century
Italian mathematician Leonardo Fibonacci in his seminal work {\it
Liber Abaci} (Sigler [48]). This sequence was originally formulated
through an idealized rabbit population model. Its standard
definition states: A sequence $\{F_{n}\}$ is called the Fibonacci
sequence if it satisfies the recurrence relation
$$F_{0}=0, F_{1}=1,F_{n+1}=F_{n}+F_{n-1}~~~(n\geq1).$$ The Fibonacci sequence exhibits a profound connection
with the golden section ratio $\frac{1+\sqrt{5}}{2}$, as evidenced
by the limit relationship
$\lim_{n\rightarrow\infty}\frac{F_{n}}{F_{n-1}}=\frac{1+\sqrt{5}}{2}$.
This further reveals a fundamental unity between discrete recursive
systems and continuous proportional relationships, which has driven
extensive investigations into the Fibonacci sequence. For instance,
19th-century French mathematician Lucas [41] not only utilized the
sequence to prove the primality of $2^{127}-1$ but also generalized
it into Lucas sequences and developed associated function theory.
Subsequent research has demonstrated the sequence's unique
properties across geometric topology, algebraic systems, and
probabilistic statistics, while also highlighting its important
value in number theory analysis and operations research and
optimization. Remarkably, this sequence is ubiquitously observed in
natural phenomena: from phyllotactic patterns in plants and
pineapple scale arrangements to patterns of tree branch growth, and
from bees returning home to quasicrystal structures. These
phenomena, with their embedded numerical signatures, have further
ignited profound academic explorations in number theory and
combinatorial mathematics. In 1963, V.E. Hoggatt and others founded
the international authoritative journal {\it Fibonacci Quarterly},
which has continuously promoted the development of the Fibonacci
theoretical system and its interdisciplinary applications.

The research on Fibonacci sequences has been constantly expanding
and shows a diversified trend. A large number of papers on the
generalized forms of Fibonacci sequences and their applications have
emerged (Miles [43], Gabai [18], Ferguson [14], Sun [50], Kessler
and Schiff [31], Howard and Cooper [22], Laipaporn et al. [34], Hare
and Saunders [21], Chen et al. [7]). A important generalization is
the $k$-generalized Fibonacci sequence (Miles [43], Gabai [18],
Kessler and Schiff [31], Howard and Cooper [22], Garc{\'i}a et al.
[19]). The study of $k$-generalized Fibonacci sequences traces back
to Miles' groundbreaking work(Miles [43]), where the sequence
$\{f_{j,k}\}$ is defined as:$$f_{j,k}=0, 0\leq j \leq k-2,
f_{k-1,k}=1, f_{j,k}=\sum_{n=1}^{k}f_{j-n,k}.$$ Miles further proved
that $det A_{j,k}=(-1)^{(2j+k)(k-1)/2}$, where
$A_{j,k}=(a_{mn})_{k\times k}$ with entries
 $a_{mn}=f_{j+m+n-2,k}$ for $1\leq m, n\leq k$. Subsequent scholars--including
 Gabai [18], Kessler and Schiff [31],
 Howard and Cooper [22], and Garc{\'i}a et al. [19]
--have expanded this framework through diverse approaches. It is
worth noting that Kessler and Schiff [31] pointed out that Miles'
paper seems to be the oldest well-known literature on the study of
such generalized Fibonacci sequences.

\subsection*{1.2\ \ Ski Rental Problem and Competitive
Analysis}

The classical ski-rental problem was first introduced in the context
of snoopy caching. The problem admits an optimal 2-competitive
strategy (Karlin et al. [28]). The formal naming of this problem as
the "ski-rental problem" and the establishment of a competitive
analysis framework under this terminology originated from Karp's
paper {\it Online Algorithms Versus Offline Algorithms:
 How Much Is It Worth to Know the Future?} (Karp [29]).
In the model, a person plans to go skiing but does not know the
exact number of days $N$ they intend to ski. The individual must
decide whether to rent or buy skis. Renting skis costs 1 unit per
day, while buying skis costs $s (s > 1)$ units , which eliminates
further rental costs. The dilemma is: Which option minimizes
costs---renting, buying immediately, or renting for a few days and
then buying? Using competitive analysis, Karp proved that the
optimal strategy is to rent skis for the first $s-1$ days and buy
them on the $s$-th day if skiing continues. This yields a
competitive ratio of $2-1/s$ . Notably, in Karp's model, if the
unknown discrete skiing days $N$ are generalized to continuous time
$T$, with rental costs calculated proportionally to the continuous
rental duration, the resulting competitive ratio aligns with the
2-competitive result obtained by Karlin in the context of snoopy
caching. As Lotker [39] pointed out in their article: ``In the rent
option, the cost is proportional to usage time, and there is no
one-time cost. The deterministic solution is straightforward (with
competitive factor 2)."

Competitive analysis, a fundamental methodology, was first
systematically introduced by Sleator and Tarjan [49] in their
seminal work in the 1980s. Since then, it has been widely applied to
the design and analysis of algorithms for various online
optimization problems in computer systems, becoming a cornerstone of
theoretical research in this field. The key metric for evaluating
the performance of online algorithms---``competitive ratio"---was
originally defined by Karlin et al. [28]. In 1994, Karlin et al.
[27] further employed competitive analysis to devise an optimal
randomized strategy for the Ski Rental Problem, achieving a
competitive ratio of $e/(e-1)$. As a classic example of
cost-minimization online optimization problems, the performance of
algorithms for the Ski Rental Problem is typically quantified using
the competitive ratio. Relevant competitive analysis models and
methods can be found in El-Yaniv et al.'s paper {\it Competitive
Optimal On-line Leasing} (El-Yaniv et al. [11]), where the ``{\bf
On-line Leasing Problem}" refers to the ``{\bf Ski Rental Problem}".
To facilitate deeper exploration and broader applications of this
classical problem, we will henceforth uniformly refer to it as the
``{\bf Online Leasing Problem}" or simply ``{\bf Online Leasing}".

The online leasing problem provides a concise and effective model
for the study of online algorithms, while competitive analysis has
laid an important theoretical foundation for the design and analysis
of online algorithms. This model and its solutions can be widely
applied in practical scenarios. Its simplest variant was first
presented by Rudolph [46] as an example of an online problem and
competitive analysis. Lotker et al. [39] [40] studied the
multi-slope online leasing problem with random strategies, based on
the real-world scenario where intermediate choices exist between
pure renting and pure buying. For the additive general model, Lotker
et al. obtained a randomized algorithm with a competitive ratio
bounded by $(e-r_{k}/r_{0})/(e-1)$. Hu and Xu [24] further improved
this result, providing a tighter
 upper bound of $e/(e-r_{k}/r_{0})$. For the non-additive variant, Levi and Patt-Shamir
[36], Chen and Xu [6] conducted in-depth analysis of the non-trivial
 characteristics in the two-option case. Other studies on the multi-slope
 online leasing problem and its variants can be found in references (Fujiwara et al.
 [16], Fujiwara et al. [17], Bohnlein et al. [4]).
 Patt-Shamir and Yadai [45] investigated the Non-linear Online Leasing problem.
 Hu [23] studied the multiple online leasing problem with continuously divisible
 assets and provided optimal algorithms and competitive ratios. In 2018, Hu et al.
 [25] extended the classical (rent-or-buy a pair of skis) online leasing problem to
 study the optimal decision-making for multiple skis with periodic payments
 (discrete time). In addition, many scholars have extended and applied the
 classical Karp model in multiple dimensions based on real-world economic
 contexts. For example, El-Yaniv et al. [11], Xu et al. [57], and Yang et al.
 [58] introduced interest rate factors into their studies; Dai et al. [9] considered the
 scenario where lessees obtain equipment usage rights through financial leasing;
 Al-Binali [1], Dong et al. [10] explored risk factors and established risk-return
 models; Zhang et al. [59] incorporated both market interest rates and equipment
 depreciation; Fujiwara and Iwama [15], Epstein and Zebedat-Haider [12], Wu et al. [55],
 and Wang et al. [54] extended their research from other economic perspectives. Overall,
 the online leasing problem and its variants hold significant value in both theoretical and
 applied aspects, with recent research further expanding their scope of application and
 theoretical depth.

\subsection*{1.3\ \ Double Online Leasing
Problem}

In Section 1.2, we observe that most research on the online leasing
 problem and its variants is based on Karp's single ski-rental model.
 However, in many practical scenarios, investors need to lease or purchase
  multiple units of equipment or assets, rather than a single unit.
  Consider the following extended case:

Suppose two people go to a ski resort to practice skiing. The cost
of purchasing one pair of skis at the resort is $s$ monetary units
$(s\gg 1)$, and the rental fee for one pair per unit time is 1
monetary units. If skis are rented, the fee is calculated
proportionally for fractional time intervals (i.e., the cost scales
continuously with actual usage time when the rental duration is not
an integer multiple of the unit time). Once they purchase a pair of
skis, they no longer need to pay any rental fees for that pair. They
can choose to rent two pairs, buy two pairs, or buy one pair and
rent the other, and so on... Faced with uncertain skiing durations
$T$, what is their optimal strategy?

We can draw an analogy between this case and the optimal rental
strategy for a single pair of skis in the classical continuous-time
model. According to this analogy, the two individuals should
initially rent two pairs of skis until s time units. If they
continue skiing after s time units, they should then purchase two
pairs of skis. The competitive ratio for this strategy is 2.
However, this is not the optimal strategy. In fact, if they choose
to: 1) rent two pairs of skis during the first
$(\sqrt{5}-1)s/2\approx 0.618s$, 2) purchase one pair of skis after
$0.618s$ time units
 and continue renting the other pair until s time units are up, and 3) purchase
 the other pair of skis after $s$ time units, the competitive ratio will be
$(5+\sqrt{5})/4\approx
 1.81$.

The results of this case are indeed surprising! When the number of
skis to be rented or purchased increases from one pair to two pairs,
not only does the classical optimal strategy lose its optimality,
but also in the strategy that achieves the smaller competitive ratio
of $(5+\sqrt{5})/4\approx
 1.81$, a magical number 0.618 appears. To distinguish it from the classical
 single-ski online leasing problem, we refer to problems of the type described
 above as the ``double online leasing problem." Moreover, we call the strategy of
 purchasing one pair of skis at the golden section point $0.618s$ within the time
 interval $[0,s]$ to achieve a competitive ratio of 1.81 the "golden section strategy."
Now, what we need to explore is whether the golden section strategy
is indeed the optimal strategy in the generalized double online
leasing problem. The answer is affirmative. Below, we will first
provide a proof for the double (two-pair ski) online leasing problem
using the game-theoretic analysis between the online decision-maker
and the offline adversary (assuming, for the sake of analysis, that
the minimum rental duration for skis in the leasing option is 1 time
unit). We will then present a generalized proof for the
multiple-pair ski online leasing problem in Section 5.

For the above double online leasing problem, it is not difficult to
see that, facing the unknown skiing duration $T$, there are only the
following six possible options for the two individuals to obtain two
pairs of skis:

$(i)$ Renting two pairs of skis for the entire duration;

$(ii)$ Purchasing two pairs of skis at the outset;

$(iii)$ Purchasing one pair of skis and renting the other until the
end of use;

$(iv)$ Purchasing one pair of skis and renting the other until time
$t(t\geq1)$, then buying it out;

$(v)$ Renting both pairs of skis until time $t(t\geq1)$, and buying
them out if further use is required;

$(vi)$ Renting both pairs of skis until time $t_{1}(t_{1}\geq1)$,
buying out one pair and continuing to use it until time
$t_{2}(t_{1}<t_{2})$, and then buying out the other pair.

Next, we briefly discuss why options $(i)$ to $(v)$ cannot be
optimal strategies.

Option $(i)$~ The offline adversary can set a sufficiently large
$\mu>2$ such that the actual usage duration is $\mu s$. In this
case, the rental cost 2$\mu s$ for the online player will far exceed
twice the offline optimal cost $2s$.

Option $(ii)$~ The offline adversary can simply set the usage
duration $T=1$. In this case, the purchase cost $2s$ for the online
player will exceed twice the actual required rental cost of 2.

Option $(iii)$~ For this strategy, the offline adversary can adopt
two countermeasures: setting $T=1$, which makes the online cost
$s+1$ exceed twice the optimal cost; or setting $T\geq ms$ (when $m$
is sufficiently large), which makes the online cost $s+ms$ far
exceed twice the offline optimal cost $2s$.

Option $(iv)$~ When the offline adversary sets $T=1$, the online
cost $s+1$ will also exceed twice the optimal cost.

Option $(v)$~ This strategy essentially reduces the double online
leasing problem for two pairs of skis to the classical online
leasing problem: the rental cost is 2 (for two pairs of skis), and
the purchase cost is $2s$ (for two pairs of skis). The optimal
strategy is to buy out the skis when the accumulated rental cost
reaches $2s$ (i.e., $t=s$), at which point the competitive ratio is
2.

From the above analysis, it can be seen that only option $(vi)$, due
to its dynamic adjustment characteristics, may potentially be the
optimal strategy. Below, we prove that the best result of option
$(vi)$ is the golden section strategy.

First, we prove that in option $(vi)$, if we take $t_{1}\geq s$,
then we can only obtain a competitive ratio greater than 2. This is
because when$t_{1}\geq s$, it follows that $t_{2}(t_{2}>t_{1})$. The
offline adversary can then set the usage time to $T=t_{2}$,
resulting in the ratio of online cost to offline cost being
\begin{eqnarray*}
\frac{t_{1}+s+t_{2}+s}{2s}>\frac{4s}{2s}=2.
\end{eqnarray*}
Next, we argue that in option $(vi)$, the optimal decision for the
online player is to take $t_{2}=s$. In fact, for any $1\leq
t_{1}<s$, the worst-case scenario for the online player when
purchasing a second skis at time $t_{2}(t_{2}>t_{1})$ is that the
offline adversary stops the skiing activity immediately afterward,
i.e., setting the skiing time to $T=t_{2}$ . In this case, the ratio
of online cost to offline cost is
\begin{eqnarray*}
\alpha=\left\{
               \begin{array}{ll}
               \frac{t_{1}+s+t_{2}+s}{2t_{2}} , & t_{2}\leq s; \\
               \frac{t_{1}+s+t_{2}+s}{2s}, & t_{2}\geq s.
               \end{array}
             \right.
\end{eqnarray*}
From the above expression, we can see that when $t_{2}\leq s$, the
ratio $\alpha$ decreases monotonically with $t_{2}$; while when
$t_{2}\geq s$, the ratio $\alpha$ increases monotonically with
$t_{2}$. The ratio $\alpha$ achieves its minimum value when
$t_{2}=s$, as shown in Figure 2.
%%%%%%%%%%%%%%%%%%%%%%%%%%%%%%%%%%%%%%%%%%%%%%%%%%%%%%%%图象插入%%%%%%%%%%%%%%%%%%
\begin{figure}[H]
  \centering
  \includegraphics[width=\textwidth, trim=5 360 5 280, clip]{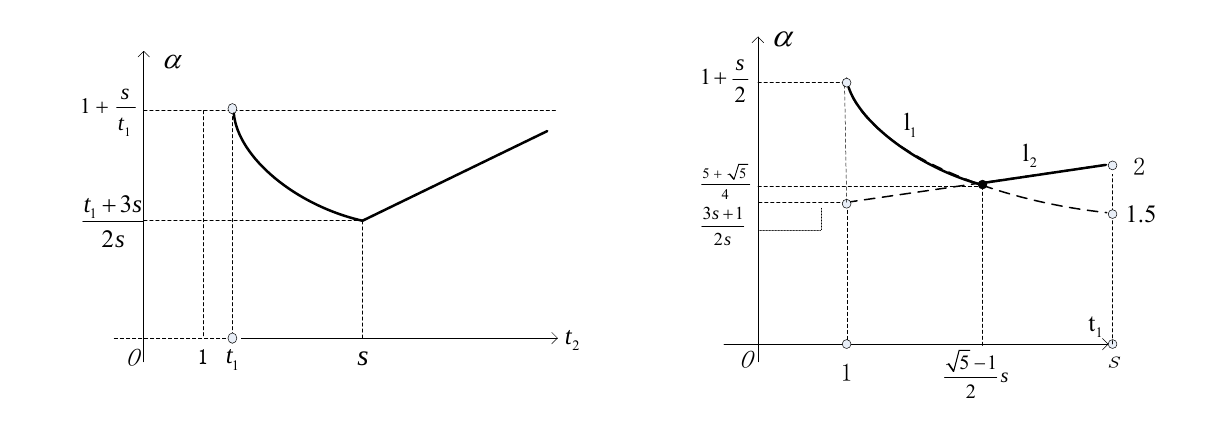}
  {\small Figure 2\hspace{6cm}Figure 3}
  \end{figure}
%%%%%%%%%%%%%%%%%%%%%%%%%%%%%%%%%%%%%%%%%%%%%%%%%%%%%%%%%%%%%%%%%%%%%%%%%%%%%%%%%%%

The analysis above shows that option $(vi)$ can achieve a smaller
competitive ratio only when $1\leq t_{1}<s$ and $t_{2}=s$. In this
case, the quality of the online decision depends on the value of the
switching point within the interval $[1, s)$, where the first pair
of skis transitions from renting to purchasing. Given the switching
point $t_{1}\in[1, s)$ chosen by the online decision-maker for the
first pair of skis and the switching point $t_{2}=s$ for the second
pair, the offline adversary can choose to stop skiing at either
$T=t_{1}$ or $T=s$ to maximize the ratio between the online and
offline costs. If the offline adversary stops the activity at
$T=t_{1}$, the ratio of online to offline costs is
\begin{eqnarray*}
\alpha=\frac{t_{1}+s+t_{1}}{2t_{1}}=1+\frac{s}{2t_{1}}
\end{eqnarray*}
which decreases monotonically with $t_{1}$ in the interval $[1, s)$,
as shown by the curve segment $l_{1}$ in Figure 3. If the offline
adversary stops the game at $T=s$, the ratio of online to offline
costs is
\begin{eqnarray*}
\alpha=\frac{t_{1}+s+s+s}{2s}=\frac{3s+t_{1}}{2s},
\end{eqnarray*}
 which increases monotonically
with $t_{1}$ in the interval $[1, s)$, as shown by the line segment
$l_{2}$ in Figure 3. When the online decision-maker chooses a
smaller $t_{1}$, the offline adversary sets the stopping time
$T=t_{1}$ to increase the ratio $\alpha$; when the online
decision-maker chooses a larger $t_{1}$, the offline adversary sets
$T=s$ to end the activity to increase the ratio $\alpha$. Therefore,
the optimal decision point $t_{1}$  for the online decision-maker
must be the value that makes the two ratios set by the offline
adversary equal, i.e., the optimal decision point $t_{1}$ satisfies
\begin{eqnarray*}
1+\frac{s}{2t_{1}}=\frac{3s+t_{1}}{2s}.
\end{eqnarray*}
From this, we obtain the optimal decision point
$t_{1}=(\sqrt{5}-1)s/2$, which is the golden section point of the
interval $[0,s]$, as shown in Figure 3.

Based on the comprehensive analysis above, for the double online
leasing problem, the optimal online strategy is to initially rent
two pairs of skis until reaching the golden section point of the
interval $[0, s]$, which is at time $(\sqrt{5}-1)s/2$. If further
usage is required beyond this point, one pair should be purchased
while continuing to rent the other until time $s$, at which point
the second pair is purchased and used until the end. Therefore, the
golden section strategy constitutes the optimal solution for the
double online leasing problem. According to the golden section
strategy, for any skiing duration $T\geq1$ set by the offline
adversary, the ratio of online cost to offline cost is given by:
\begin{eqnarray*}
\alpha=\left\{
               \begin{array}{ll}
               \frac{2T}{2T}=1<\frac{5+\sqrt{5}}{4} , & T<\frac{\sqrt{5}-1}{2}s; \\
               \frac{2\frac{\sqrt{5}-1}{2}s+s}{2\frac{\sqrt{5}-1}{2}s}=\frac{5+\sqrt{5}}{4}, & T=\frac{\sqrt{5}-1}{2}s;\\
               \frac{\frac{\sqrt{5}-1}{2}s+s+T}{2T}<\frac{5+\sqrt{5}}{4}, & \frac{\sqrt{5}-1}{2}s<T<s;\\
               \frac{\frac{\sqrt{5}-1}{2}s+s+s+s}{2s}=\frac{5+\sqrt{5}}{4}, & T\geq s.
               \end{array}
             \right.
\end{eqnarray*}
Thus, the competitive ratio of the golden section strategy is
$(5+\sqrt{5})/4\approx1.81$.

Above, we extended the classical online leasing problem of a single
pair of skis to the online leasing of two pairs of skis, thereby
proposing the double online leasing problem. Through competitive
analysis, we have shown that the optimal deterministic online
strategy is the golden section strategy, and we have improved the
optimal competitive ratio of the original classical problem from 2
to $(5+\sqrt{5})/4\approx1.81$. It is worth noting that the optimal
strategy reveals an elegant mathematical property: the optimal
timing for purchasing the first pair of skies coincides exactly with
the golden section point of the interval $[0,s]$. As a universally
optimal proportion in nature, the golden section has profound
manifestations in architecture, art, and combinatorial optimization.
The perfect fit between the double online leasing problem and the
golden section raises a deeper question: Is this connection
universal? To explore this question, we will conduct the following
research:

$\bullet$ Introducing multiple division points within the interval
$[0,1]$ to construct a high-dimensional golden section theory and
the bi-infinite $m$-generalized Fibonacci sequence;

$\bullet$ Considering and investigating the multiple online leasing
problem involving the rental or purchase of multiple pairs of skis.

\subsection*{1.4\ \ Content and Organization }

The main content and structure of this paper are organized as
follows:

$\bullet$ {\bf Section 2} generalizes the classical golden section
to higher dimensions by proposing the geometric-type
higher-dimensional golden section and the harmonic-type
higher-dimensional golden section.

$\bullet$ {\bf Section 3} extends the generalized Fibonacci sequence
to a bi-infinite generalized Fibonacci sequence and investigates the
relationship between higher-dimensional golden sections (points) and
bi-infinite generalized Fibonacci sequences.

$\bullet$ {\bf Section 4} introduces the Fibonacci matrix for
harmonic-type higher-dimensional golden section points and discusses
the properties of its submatrices.

$\bullet$ {\bf Section 5} extends the classical online leasing to a
multiple online leasing model, presenting and proving that its
optimal strategy coincides precisely with the harmonic-type
higher-dimensional golden section method.

$\bullet$ {\bf Section 6} considers the continuous relaxation
problem for multiple online leasing and its optimal solution.

$\bullet$ {\bf Section 7} demonstrates the feasibility of
approximating the optimal solution using relaxed solutions from
multiple online leasing through numerical analysis.

$\bullet$ {\bf Section 8} proposes a correction method for the
relaxed solutions based on the classical golden section---the relaxed
solutions golden-section mantissa weighted method.

$\bullet$ {\bf Section 9} concludes the paper and outlines potential
directions for future research.

\section{ Higher-dimensional Golden Section}

For convenience and consistency, we refer to a complete unit line
segment $AB$ as the one-dimensional golden section, and its endpoint
$B$ as the one-dimensional golden section point; we call the
classical golden section the two-dimensional golden section. In
conjunction with Figure 1, under the definitions, assumptions, and
notations of the two-dimensional golden section in Section 1.1.1, it
is not difficult to prove that the following four propositions are
equivalent:

$(i)$ The point $C$ on the unit line segment $AB$ is the
two-dimensional golden section point;
$$(ii)\hspace{5cm}\frac{AC}{CB}=\frac{1+\sqrt{5}}{2};\hspace{7.2cm}$$
$$(iii)\hspace{3.9cm}\frac{AB}{AC}=\frac{AC}{CB}~~~\hbox{or}~~~\frac{1}{\lambda}
+\frac{1}{\lambda^{2}}=1;\hspace{6cm}$$
 $$(iv)\hspace{3.6cm}\frac{AC}{CB}=1+AC~~~\hbox{or}~~~\frac{1}{\omega}
+\frac{1}{\omega+1}=1.\hspace{5.3cm}$$

We now consider a more general problem: How to partition the unit
line segment $AB$ into $m~(m\geq 2)$ segments ---$AC_{1},
C_{1}C_{2}, \cdots, C_{m-2}C_{m-1}, C_{m-1}B$---such that this
partitioning possesses biological characteristics, aesthetic value,
 architectural significance, and applications in operations research optimization, similar
 to those of the classical golden section?
The most intuitive approach is based on the equivalent condition
$(ii)$ of the two-dimensional golden section, where $m-1$ division
points $C_{1}, C_{2}, $ $\cdots, C_{m-2}, C_{m-1}$
   are inserted into $AB$ so that the ratios of adjacent segments satisfy the golden section ratio:
\begin{eqnarray*}
\frac{AC_{1}}{C_{1}C_{2}}=\frac{C_{1}C_{2}}{C_{2}C_{3}}=\cdots=\frac{C_{m-2}C_{m-1}}{C_{m-1}B}
=\frac{1+\sqrt{5}}{2}.
\end{eqnarray*}
 However, this equal-ratio partitioning only holds when $m=2$.
When $m>2$, the sum of the segment lengths no longer satisfies the
unit length constraint, i.e., $AB+AC_{1}+C_{1}C_{2}+\cdots+
C_{m-2}C_{m-1}+C_{m-1}B\neq1$.

Below, we will propose two generalized forms based on the equivalent
conditions $(iii)$ and $(iv)$ of the two-dimensional golden
section---{\bf geometric-type higher-dimensional golden section} and
{\bf harmonic-type higher-dimensional golden section}.

\subsection*{2.1 \ Geometric-type Higher-dimensional Golden
Section}

\noindent{\bf Definition 2.1}~ Based on the equivalent condition
$(iii)$ of the two-dimensional golden section, insert $m-1$ division
points $A=C_{0}, C_{1}, C_{2}, \cdots, C_{m-2},$ $C_{m-1}, C_{m}=B$
on the unit line segment $AB$ to partition it into $m$ segments
$$C_{0}C_{1},~ C_{1}C_{2},~ C_{2}C_{3},~ \cdots,~
C_{m-2}C_{m-1}, C_{m-1}C_{m}$$ such that
$$\frac{C_{0}C_{m}}{C_{0}C_{1}}=\frac{C_{0}C_{1}}{C_{1}C_{2}}=\frac{C_{1}C_{2}}{C_{2}C_{3}}=
\cdots=\frac{C_{m-2}C_{m-1}}{C_{m-1}C_{m}}.$$
 This is called the
$m$-dimensional geometric-type golden section of the line segment
$AB$, as shown in Figure 4.

%%%%%%%%%%%%%%%%%%%%%%%%%%%%%%%%%%%%%%%%%%%%%%%%%%%%%%%%图象插入%%%%%%%%%%%%%%%%%%
\begin{figure}[H]
  \centering
  \includegraphics[width=\textwidth, trim=5 540 5 250, clip]{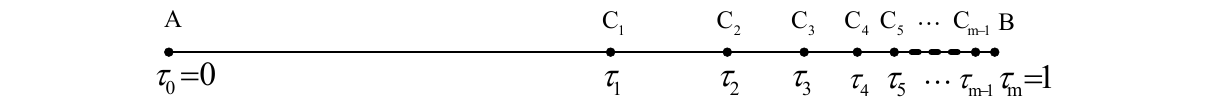}
  {\small Figure 4:~~$m$-dimensional geometric-type
golden section}
  \end{figure}
%%%%%%%%%%%%%%%%%%%%%%%%%%%%%%%%%%%%%%%%%%%%%%%%%%%%%%%%%%%%%%%%%%%%%%%%%%%%%%%%%%%

 The ratio $\frac{C_{0}C_{m}}{C_{0}C_{1}}$
  is called the $m$-dimensional geometric-type golden section ratio, denoted as $\lambda$.
  Its reciprocal $\phi=\frac{1}{\lambda}$
  is called the $m$-dimensional geometric-type golden section rate. The points
  $C_{1}, C_{2}, \cdots, C_{m-1}, C_{m}$
  on the unit line segment $AB$ are called the $m$-dimensional geometric-type
  golden section points.
  Under normal circumstances, the unit line segment $AB$ is regarded as the
  interval $[0,1]$ on the real axis.
  The $m$-dimensional geometric-type golden section points $C_{1}, C_{2}, \cdots, C_{m-1}, C_{m}$
  are represented by their coordinates $\tau_{1},
\tau_{2}, \cdots, \tau_{m-1}, \tau_{m}$. The left and right
endpoints of the interval $[0,1]$ are denoted as $\tau_{0}=0$ and
$\tau_{m}=1$, respectively. The $m$-dimensional geometric-type
golden section on the interval $[0,1]$ is denoted as
$\tau^{(m)}(\tau_{m}, \tau_{m-1}, \cdots, \tau_{2}, \tau_{1} )$,
which is briefly written as $\tau^{(m)}$. The interval $[0,1]$ is
divided into $m$ subintervals, called the section intervals of the
$m$-dimensional
 geometric-type golden section, and these subintervals and their lengths are denoted by $\delta_{1}, ~\delta_{2},
~\cdots, ~\delta_{m}$.

When $m$ tends to infinity, it is referred to as the
infinite-dimensional geometric-type golden section. The
geometric-type golden section in three dimensions and higher is
collectively referred to as the geometric-type higher-dimensional
golden section.

According to Definition 2.1, the $m$-dimensional geometric-type
golden section ratio $\lambda$ on the interval $[0,1]$ satisfies
$$\frac{1}{\lambda}+\frac{1}{\lambda^{2}}+\cdots+\frac{1}{\lambda^{m}}=1~~\hbox{or}~~
\lambda^{m}-\lambda^{m-1}-\cdots-\lambda^{2}-\lambda-1=0.$$ In
 general, we call the equations
\begin{eqnarray}
\frac{1}{x}+\frac{1}{x^{2}}+\cdots+\frac{1}{x^{m}}=1~~\hbox{or}~~
x^{m}-x^{m-1}-\cdots-x^{2}-x-1=0
\end{eqnarray}
the $m$-dimensional geometric-type golden section equation.

\vspace{1.5mm}\noindent{\bf Lemma 2.1} {\it For every $m \in \{2, 3,
\cdots, k, \cdots\}$, there exists a unique real root
$\lambda^{(m)}$ of equation (1) in the interval $(\lambda^{(m-1)},
2)$. Moreover, $f'_{m}(x) > 0$ holds in $(\lambda^{(m-1)}, 2)$,
where $\lambda^{(1)}=1, f_{m}(x)=x^{m}-x^{m-1}-\cdots-x^{2}-x-1$}.

\vspace{1.5mm}{\bf Proof}~ We proceed by induction on $m$.

When $m=2$, the equation $x^{2}-x-1=0$ clearly has a unique real
root $\lambda^{(2)}=(1+\sqrt{5})/2$
  in $(\lambda^{(1)},2)$. Moreover, in $(\lambda^{(1)}, 2)$, $f_{2}'(x)=2x-1>0$.

For $m=3$, we have \begin{eqnarray*}
f_{3}(x)=x^{3}-x^{2}-x-1=x(x^{2}-x-1)-1,
\end{eqnarray*}
which yields $f_{3}(\lambda^{(2)})=-1<0$ ~and~$f_{3}(2)=1>0$. The
derivative satisfies \begin{eqnarray*}
f_{3}'(x)=3x^{2}-2x-1=x(2x-1)+x^{2}-x-1>0
\end{eqnarray*} for $x\in(\lambda^{(2)}, 2)$, establishing the existence of a unique
root $\lambda^{(3)} \in (\lambda^{(2)}, 2)$.

Assume the statement holds for $m=k$: equation (1) has a unique root
$\lambda^{(k)} \in (\lambda^{(k-1)},2)$ with
\begin{eqnarray*}
f_{k}'(x)
= kx^{k-1} - (k-1)x^{k-2} - \cdots - 2x - 1 > 0.
\end{eqnarray*}

For $m=k+1$, observe that
\begin{eqnarray*}
f_{k+1}(x)=x^{k+1}-x^{k}-\cdots-x^{2}-x-1=x(x^{k}-x^{k-1}-\cdots-x^{2}-x-1)-1,
\end{eqnarray*}
By the induction hypothesis, we have $f_{k+1}(\lambda^{(k)})=-1<0$
and
\begin{eqnarray*}
f_{k+1}'(x)&=&(k+1)x^{k}-kx^{k-1}-\cdots-2x-1\\
&=&x[kx^{k-1}-(k-1)x^{k-2}-\cdots-2x-1]+x^{k}-x^{k-1}-\cdots-x^{2}-x-1\\
&>&0.
\end{eqnarray*}
Since $~f_{k+1}(2)=1>0$, equation (1) has a unique real root
$\lambda^{(k+1)}$ in $(\lambda^{(k)},2)$. The induction is complete.

\vspace{1.5mm}\noindent{\bf Theorem 2.1} {\it For the
$m$-dimensional geometric-type golden section equation (1), when
$m=2,3,\cdots,k,\cdots$, there exist the largest real roots $
 \lambda^{(2)},\lambda^{(3)},\cdots,\lambda^{(k)},\cdots,\lambda^{(\infty)}$
  which are strictly monotonically increasing with $m$, and these roots satisfy the following
  inequalities:
\begin{eqnarray*}
\frac{1+\sqrt{5}}{2}=\lambda^{(2)}<\lambda^{(3)}<\cdots<\lambda^{(k)}<\cdots<\lambda^{(\infty)}=2.
\end{eqnarray*}}

 \vspace{1.5mm}{\bf Proof} Since for all
$m\in\{2,3,\cdots,k,\cdots\}$, when $x\geq2$,
\begin{eqnarray*}
f'_{m}(x)=\frac{mx^{m-1}(x-2)}{x-1}+\frac{x^{m}-1}{(x-1)^{2}}>0.
\end{eqnarray*} and
$f_{m}(2)=1>0$, it follows that all real roots of equation (1) are
less than 2. Furthermore, if there exists some $1\leq\lambda'<2$
such that $\lambda$ is the unique real root of equation (1) in the
interval $(\lambda',2)$, then $\lambda$ must be the largest real
root of the equation. By Lemma 2.1, equation (1) has a unique
largest real root $\lambda^{(m)}$ in the interval $[(1+\sqrt{5})/2,
2)$, and the sequence $\{\lambda^{(m)}\}$ (for $m = 2, 3, \cdots$)
is strictly increasing and bounded above by 2.
 Therefore, the $\lim_{m \to \infty} \lambda^{(m)}=\lambda^{(\infty)}$ exists.
Moreover, since $\lambda^{(m)}$ is a root of equation (1), we have
$$(\lambda^{(m)})^{m}-(\lambda^{(m)})^{m-1}-\cdots-(\lambda^{(m)})^{2}
-\lambda^{(m)}-1=0,$$ which can be rewritten as
$(\lambda^{(m)})^{m}=\frac{(\lambda^{(m)})^{m}-1}{\lambda^{(m)}-1}$
or equivalently $\lambda^{(m)}=2-\frac{1}{(\lambda^{(m)})^{m}}.$
Taking the limit as $m \to \infty$ on both sides, we obtain
$\lim_{m\rightarrow\infty}\lambda^{(m)}=\lambda^{(\infty)}=2$.  The
proof of Theorem 2.1 is complete.

According to Theorem 2.1, the $m$-dimensional geometric-type golden
section equation (1) has a unique largest real root in the interval
$[(\sqrt{5}+1)/2, 2]$. Therefore, the $m$-dimensional geometric-type
golden section ratio $\lambda$ defined above specifically refers to
the largest real root of equation (1). For example:

 The classical
golden section ratio (2-dimensional geometric-type golden section
ratio) $\lambda$ is the largest real root of the equation
$x^{2}-x-1=0:\lambda=\frac{1+\sqrt{5}}{2}\approx 1.618033989$, with
the golden section rate $\phi=\frac{\sqrt{5}-1}{2}\approx
0.618033989$, the golden section point
$\tau_{1}=\frac{\sqrt{5}-1}{2}, \tau_{2}=1$; and the interval
$[0,1]$ divided into two subintervals:
$$\delta_{1}=\frac{\sqrt{5}-1}{2},~~\delta_{2}=\frac{3-\sqrt{5}}{2}.$$

The 3-dimensional geometric-type golden section ratio $\lambda$ is
the largest real root of the equation
$x^{3}-x^{2}-x-1=0$:$$\lambda=\frac{1+\sqrt[3]{19+3\sqrt{33}}+\sqrt[3]{19-3\sqrt{33}}}{3}\approx1.839286755,$$
with the golden section rate $\phi\approx0.5436890128,$ the golden
section points
$$\tau_{1}\approx0.5436890128,~~\tau_{2}\approx0.8392867554,~~
\tau_{3}=1$$; and the interval $[0,1]$ divided into three
subintervals:
$$\delta_{1}\approx0.5436890128,~~\delta_{2}\approx0.2955977426,~~\delta_{3}\approx0.1607132446.$$

The infinite-dimensional geometric-type golden section ratio
$\lambda$ is the largest real root of the infinite- dimensional geometric-type golden section equation, which is
$\lambda=2$, with the golden section rate value $\phi=\frac{1}{2},$
the golden section points:
$$\tau_{1}=\frac{1}{2},~\tau_{2}=\frac{3}{4},~
\cdots,~\tau_{n}=\frac{2^{n}-1}{2^{n}},~\cdots ,$$ and the interval
[0,1] divided into infinitely many subintervals:
$$\delta_{1}=\frac{1}{2}, ~\delta_{2}=\frac{1}{2^{2}},
 ~\cdots, ~\delta_{n}=\frac{1}{2^{n}}, ~\cdots~ .$$

The significance of the infinite-dimensional geometric-type golden
section perfectly illustrates a quote from the ancient Chinese
philosopher Zhuang Zhou in  {\it Zhuangzi $\cdot$ The World}
chapter: ``A rod a foot long is cut in half each day, yet it will
never be exhausted for all eternity."

\vspace{1.5mm}\noindent{\bf Theorem 2.2}~{\it Equation (1) has no
repeated roots and possesses no negative real roots less than or
equal to $-1$.}

\vspace{1.5mm}{\bf Proof}~ $(i)$ Since $x=1$ is not a root of
equation (1), the equation can be rewritten as
\begin{eqnarray}
x^{m+1}-2x^m +1=0.
\end{eqnarray}
Except for the single root $x=1$, equation (2) shares identical
roots with equation (1). The Sylvester discriminant of equation (2)
is
\begin{eqnarray*}
D&=&(-1)^{\frac{m(m+1)}{2}}\left|
\begin{array}{ccccccccccccc}
1      & -2    & 0      & 0      & \cdots & 0     & 0      & 1      & 0      & 0      &\cdots &0      & 0    \\
0      &1      & -2     & 0      & \cdots & 0     & 0      & 0      & 1      & 0      &\cdots &0      & 0    \\
0      & 0     &1       & -2     & \cdots & 0     & 0      & 0      & 0      & 1      &\cdots &0      & 0    \\
\vdots &\vdots & \vdots & \vdots & \vdots &\vdots & \vdots & \vdots & \vdots & \vdots &\vdots &\vdots &\vdots\\
0      & 0     & 0      & 0      & \cdots & 1     & -2     & 0      & 0      & 0      &\cdots &1      & 0    \\
0      & 0     & 0      & 0      & \cdots & 0     & 1      & -2     & 0      & 0      &\cdots &0      & 1    \\
m+1    & -2m   & 0      & 0      & \cdots & 0     & 0      & 0      & 0      & 0      &\cdots &0      & 0    \\
0      & m+1   & -2m    & 0      & \cdots & 0     & 0      & 0      & 0      & 0      &\cdots &0      & 0    \\
\vdots &\vdots & \vdots & \vdots & \vdots &\vdots & \vdots & \vdots & \vdots & \vdots &\vdots &\vdots &\vdots \\
0      & 0     & 0      & 0      & \cdots & m+1   & -2m    &  0     & 0      & 0      &\cdots &0      & 0\\
0      & 0     & 0      & 0      & \cdots & 0     & m+1    & -2m    & 0      & 0      &\cdots &0      & 0\\
 \end{array}
\right|\\
\end{eqnarray*}
\begin{eqnarray*}
&=&(-1)^{\frac{m(m+1)}{2}}\left|
   \begin{array}{cccccccc}
   1      & -2    & 0      & 0      & \cdots & 0     & 0      & 1      \\
   m+1    & -2m   & 0      & 0      & \cdots & 0     & 0      & 0      \\
   0      & m+1   & -2m    & 0      & \cdots & 0     & 0      & 0       \\
   \vdots &\vdots & \vdots & \vdots & \vdots &\vdots & \vdots & \vdots  \\
   0      & 0     & 0      & 0      & \cdots & -2m   &  0     & 0       \\
   0      & 0     & 0      & 0      & \cdots & m+1   & -2m    & 0       \\
   0      & 0     & 0      & 0      & \cdots & 0     & m+1    & -2m      \\
    \end{array}
\right|\\
&=&(-1)^{\frac{m(m+1)}{2}+2m(m+1)}(m+1)^{m+1}+(-1)^{\frac{m(m+1)}{2}+2m(m-1)+11}2(2m)^{m}.\\
\end{eqnarray*}
From the above, we obtain
\begin{eqnarray}
D=\left\{
         \begin{array}{ll}
           5, & m=2; \\
           (m+1)^{m+1}-2(2m)^{m}, & m\geq3.
         \end{array}
       \right.
\end{eqnarray}
For $m\geq3$,
\begin{eqnarray*}
 \frac{(m+1)^{m+1}}{2(2m)^{m}}=\frac{m+1}{2^{m+1}}(1+\frac{1}{m})^{m}<\frac{m+1}{2^{m-1}}\cdot
\frac{3}{4}\leq\frac{3}{4}<1,
 \end{eqnarray*}
which implies $D<0$ when $m\geq3$. Therefore, equation (2) has no
repeated roots, and consequently, equation (1) has no repeated
roots. For $m=2, D=5>0$, and equation (2) has three distinct real
roots. Naturally, equation (1) has two distinct real roots, meaning
no repeated roots.

$(ii)$ Rewriting equation (2) gives
\begin{eqnarray}
x = 2 - \frac{1}{x^m}.
\end{eqnarray}

When $m=2k+1~(k\geq1)$ and $x\leq-1$, the left-hand side of equation
(4) satisfies $x\leq-1$, while the right-hand side satisfies
$2-\frac{1}{x^{2k+1}}>2.$

When $m=2k~(k\geq1)$ and $x\leq-1$, the left-hand side of equation
(4) satisfies $x\leq-1$, while the right-hand side satisfies
$2-\frac{1}{x^{2k}}>1.$

Combining these results, equation (1) has no negative real roots
less than or equal to $-1$. The proof of Theorem 2.2 is complete.

\subsection*{2.2\ Harmonic-type Higher-dimensional Golden
Section}

\noindent{\bf Definition 2.2}~Based on the equivalent condition
$(iv)$ of the 2-dimensional golden section, insert $m-1$ points
$A=C_{0}, C_{1}, C_{2}, \cdots, C_{m-2},$ $C_{m-1},  C_{m}=B$ on the
unit line segment $AB$, dividing the unit line segment $AB$ into $m$
segments $$C_{0}C_{1},~ C_{1}C_{2},~ C_{2}C_{3},~ \cdots,~
C_{m-2}C_{m-1},~ C_{m-1}C_{m}$$ such that
\begin{eqnarray}
\frac{C_{0}C_{1}}{C_{1}C_{2}}=1+C_{0}C_{1}, \quad
\frac{C_{1}C_{2}}{C_{2}C_{3}}=1+C_{1}C_{2}, \quad
\frac{C_{2}C_{3}}{C_{3}C_{4}}=1+C_{2}C_{3}, \quad \cdots, \quad
\frac{C_{m-2}C_{m-1}}{C_{m-1}C_{m}}=1+C_{m-2}C_{m-1}.
\end{eqnarray}
This is called the $m$-dimensional harmonic-type golden section of
the line segment $AB$. As shown in Figure 5.

%%%%%%%%%%%%%%%%%%%%%%%%%%%%%%%%%%%%%%%%%%%%%%%%%%%%%%%%图象插入%%%%%%%%%%%%%%%%%%
\begin{figure}[H]
  \centering
  \includegraphics[width=\textwidth, trim=5 510 5 275, clip]{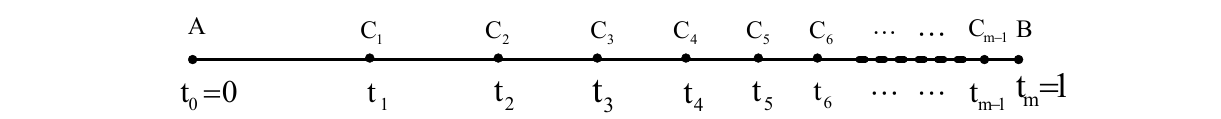}
  {\small Figure 5:~~$m$-dimensional harmonic-type
golden section}
  \end{figure}
%%%%%%%%%%%%%%%%%%%%%%%%%%%%%%%%%%%%%%%%%%%%%%%%%%%%%%%%%%%%%%%%%%%%%%%%%%%%%%%%%%%

 The ratio $\frac{C_{0}C_{m}}{C_{0}C_{1}}$ is called
the section basis of the $m$-dimensional harmonic-type golden
section and is denoted by $\omega$. The section rate of the
$m$-dimensional harmonic-type golden section is defined as
$\varphi=\frac{\omega}{m}$, and the balance ratio or competitive
ratio of the $m$-dimensional harmonic-type golden section is
$\alpha=\varphi+1$. The points $C_{1}, C_{2}, \cdots, C_{m-1}$ and
the right endpoint of the line segment $ C_{m}$
  are called the $m$-dimensional harmonic-type golden section points.
The $m$-dimensional harmonic-type golden section points on the
interval $[0,1]$ are represented by their coordinates $t_{1}, t_{2},
\cdots, t_{m-1}, t_{m}$, while the left and right endpoints of the
interval $[0,1]$ are denoted by $t_{0}=0$ and $t_{m}=1$,
respectively. The $m$-dimensional harmonic-type golden section on
the interval $[0,1]$ is denoted by $T^{(m)}: t_{m}, t_{m-1}, \cdots,
t_{2},  t_{1}$, and is briefly written as $T^{(m)}$. The interval
$[0,1]$ is divided into $m$ subintervals, which are called the
section intervals of the $m$-dimensional harmonic-type golden
section. These subintervals and their lengths are denoted by
$\Delta_{1}, ~\Delta_{2}, ~\cdots, ~\Delta_{m}$.

From $AB=1$ and the equations in (5), we have
$$C_{0}C_{1}=\frac{1}{\omega},
C_{1}C_{2}=\frac{1}{1+\frac{1}{C_{0}C_{1}}},
~~C_{2}C_{3}=\frac{1}{2+\frac{1}{C_{0}C_{1}}},~~ \cdots,~~
C_{m-1}B=\frac{1}{m-1+\frac{1}{C_{0}C_{1}}}.$$ By Definition 2.2,
the section basis $\omega$ of the $m$-dimensional harmonic-type
golden section on the interval $[0,1]$ satisfies
\begin{eqnarray*}
\frac{1}{\omega}+\frac{1}{\omega+1}+\frac{1}{\omega+2}+\cdots+\frac{1}{\omega+m-1}=1.
\end{eqnarray*}
It is evident that when $m=2$, $\lambda=\omega=(1+\sqrt{5})/2,$
which means that the 2-dimensional geometric-type golden section and
the 2-dimensional harmonic-type golden section are both the
classical golden section.

In general, we refer to the following equation (6) as the
$m$-dimensional harmonic-type golden section equation.
\begin{eqnarray}
\frac{1}{x}+\frac{1}{x+1}+\frac{1}{x+2} +\cdots+\frac{1}{x+m-1}=1
\end{eqnarray}

The harmonic-type golden sections of dimension 3 and higher are
collectively referred to as the harmonic-type higher-dimensional
golden sections.

\vspace{2mm}\noindent{\bf Theorem 2.3}~{\it For any positive integer
$m\geq 2$, the $m$-dimensional harmonic-type golden section on the
interval $[0,1]$ exists and is unique.}

\vspace{1.5mm}{\bf Proof}~By Definition 2.2, the validity of this
theorem implies that the $m$-dimensional harmonic-type golden
section equation (6) has a unique positive real solution $\omega$.
Thus, it suffices to prove that the equation in terms of $\varphi$
\begin{eqnarray}
\frac{1}{m\varphi}+\frac{1}{m\varphi+1}+\frac{1}{m\varphi+2}+\cdots+\frac{1}{m\varphi+m-1}=1
\end{eqnarray}
has a unique positive real solution. In fact, equation (7) has a
unique positive real solution in the interval $(1/(e-1),1)$. Let
\begin{eqnarray*}
f(\varphi)=\frac{1}{m\varphi}+\frac{1}{m\varphi+1}+\cdots+\frac{1}{m\varphi+m-1}-1,
\end{eqnarray*}
Clearly, for any fixed $m\in \mathbb{Z}^{+}$
  with $m\geq 2$, $f(\varphi)$ is a continuous function of $\varphi$ on $(0, +\infty)$, and we have
\begin{eqnarray*}
f(1)&=&\frac{1}{m}+\frac{1}{m+1}+\cdots+\frac{1}{m+m-1}-1\\
&<&\frac{m}{m}-1=0,
\end{eqnarray*}
and
\begin{eqnarray*}
f\left(\frac{1}{e-1}\right)&=&\frac{1}{m\frac{1}{e-1}}+\frac{1}{m\frac{1}{e-1}+1}+\cdots+\frac{1}{m\frac{1}{e-1}+m-1}-1\\
&>&\ln\left(1+\frac{1}{m\frac{1}{e-1}}\right)+\ln\left(1+\frac{1}{m\frac{1}{e-1}+1}\right)+\cdots+
\ln\left(1+\frac{1}{m\frac{1}{e-1}+m-1}\right)-1\\
&=&\ln\left[\frac{m+e-1}{m}\cdot\frac{m+2(e-1)}{m+e-1}\cdot\cdots\cdot\frac{m+m(e-1)}
{m+(m-1)(e-1)}\right]-1\\ &=&0.
\end{eqnarray*}
By the Intermediate Value Theorem for continuous functions, equation
(7) has at least one real root in the interval $(1/(e-1),1)$.
Moreover, since on $(0,\infty)$
\begin{eqnarray*}
f'(\varphi)=-\frac{m}{m^{2}\varphi^{2}}-\frac{m}{[m\varphi+1]^{2}}-\cdots-\frac{m}{[m\varphi+m-1]^{2}}<0,
\end{eqnarray*}
$f(\varphi)$ is a strictly decreasing function on $(0,\infty)$.
Therefore, equation (7) has a unique positive real solution in the
interval $(1/(e-1),1)$. Consequently, equation (6) has a unique
positive real solution in the interval $(\frac{m}{e-1},m)$. The
proof of Theorem 2.3 is complete.

The section basis $\omega$ of the 3-dimensional harmonic-type golden
section is the unique positive real root of the equation
$$\frac{1}{x}+\frac{1}{x+1}+\frac{1}{x+2}=1,$$ which is
$x^{3}-4x-2=0$:$$\omega=\frac{4\sqrt{3}}{3}\cos(\frac{1}{3}\arctan\frac{\sqrt{111}}{9})\approx2.214319743.$$
The section rate $\varphi\approx0.738106581,$ the balance ratio or
competitive ratio $\alpha\approx1.738106581,$ The section points
$t_{1}\approx0.451605963,t_{2}\approx0.7627137805, t_{3}=1.$ The
section intervals of the interval $[0,1]$ are:
$\Delta_{1}\approx0.451605963,\Delta_{2}\approx0.311107817,\Delta_{3}\approx0.237286220.$

Since the harmonic-type golden section equation in dimensions 5 and
higher may not have explicit solutions, below, we will perform some
transformations on equation (6) and conduct analysis and discussions
on the properties of its roots.

The function
\begin{eqnarray*}
f(x) = \frac{1}{x} + \frac{1}{x+1} + \frac{1}{x+2} + \cdots +
\frac{1}{x+m-1} - 1
\end{eqnarray*}
is piecewise continuous on the intervals
$(-\infty,-m+1)\bigcup(-m+1,-m+2)\bigcup(-m+1,-m+2)\bigcup\cdots\bigcup(-1,0)$,
with\begin{eqnarray*}
f_{-}(0)=f_{-}(-1)=\cdots=f_{-}(-m+1)=-\infty,\\
f_{+}(-1)=f_{+}(-2)=\cdots=f_{+}(-m+2)=+\infty.
\end{eqnarray*}
Moreover, on each interval, the derivative satisfies
\begin{eqnarray*}
f'(x)=-\frac{1}{x^{2}}-\frac{1}{(x+1)^{2}}-\cdots-
\frac{1}{(x+m-1)^{2}} < 0.
\end{eqnarray*}
Thus, equation (6) has one negative real root in each of the
intervals $(-m+1,-m+2),~(-m+1,-m+2),~\cdots,~(-2,-1),~(-1,0)$,
yielding a sequence of negative roots:
$$-m+1<x_{m-1}<x_{m-2}<\cdots<x_{2}<x_{1}<0.$$ For $x<-m+1$,
$f(x)<0$ always holds.

 Combining the theorem 2.3 with the above
analysis, we conclude: The $m$-dimensional harmonic-type golden
section equation (6) has no complex roots---all roots are simple and
real. Apart from the unique positive real root $\omega$, the
remaining $m-1$ roots are negative real roots:
$0>x_{1}>x_{2}>\cdots>x_{m-1}>-m+1.$

\vspace{1.5mm}\noindent{\bf Theorem 2.4} {\it The section basis
$\omega$ of the $m$-dimensional harmonic-type golden section
strictly increases without upper bound as $m$ increases, and the
section rate and the balance ratio or competitive ratio $\varphi,
\alpha$ both strictly decrease as $m$ increases. Moreover,
$$\lim_{m\rightarrow\infty}\varphi=\frac{1}{e-1},~~\lim_{m\rightarrow\infty}\alpha
=\frac{e}{e-1}.$$}

 \vspace{1.5mm}{\bf Proof}~ It is obvious that $\omega$
strictly increases without upper bound as $m$ increases. Below, we
prove that $\varphi$ strictly decreases as $m$ increases. For the
time being, denote $\varphi$ in equation (7) as $\varphi_{m}$.

We first prove the following inequality:
\begin{eqnarray}
\sum^{m-1}_{i=0}\frac{1}{m\varphi_{m}+i} >
\sum^{m}_{i=0}\frac{1}{(m+1)\varphi_{m}+i}.
\end{eqnarray}

Let
\begin{eqnarray*}
X_{i} = \frac{1}{m\varphi_{m}+i} - \frac{1}{(m+1)\varphi_{m}+i},
\quad i=0,1,\cdots,m-1, \quad Z_{m} = \frac{1}{(m+1)\varphi_{m}+m},
\end{eqnarray*}
then the problem reduces to proving
\begin{eqnarray*}
X_{0} + X_{1} + \cdots + X_{m-1} > Z_{m}.
\end{eqnarray*}

Assume, for contradiction, that
\begin{eqnarray*}
X_{0} + X_{1} + \cdots + X_{m-1} \leq Z_{m}.
\end{eqnarray*}
Then,
\begin{eqnarray}
\ln(1 + X_{0} + X_{1} + \cdots + X_{m-1}) - \ln(1 + Z_{m}) \leq 0.
\end{eqnarray}
On the other hand, since the function $f(x)=\ln(1+x)$ satisfies the
conditions of the Lagrange Mean Value Theorem on the interval
$[0,X_{0}+X_{1}+\cdots+X_{m-1}]$, there exists
$\xi\in(0,X_{0}+X_{1}+\cdots+X_{m-1})$, i.e.,
$0<\xi<\frac{1}{(m+1)\varphi_{m}+m}$, such that
\begin{eqnarray}
\ln(1 + X_{0} + X_{1} + \cdots + X_{m-1}) &=& \ln(1 + X_{0} + X_{1}
+ \cdots + X_{m-1}) - \ln(1 + 0) \nonumber \\
&=& \frac{1}{1 + \xi}(X_{0} + X_{1} + \cdots + X_{m-1}).
\end{eqnarray}
 Furthermore, it can be
verified that $f(x)=\ln(1+x)$ satisfies the Lagrange Mean Value
Theorem conditions on each interval
$$I_{i}=\Big[\frac{1}{(m+1)\varphi_{m}+i},\frac{1}{m\varphi_{m}+i}\Big]
\ \hbox{for all} \ i=0,1,\cdots,m-1,$$ Thus, there exist
$$\xi_{i}\in\Big(\frac{1}{(m+1)\varphi_{m}+i},\frac{1}{m\varphi_{m}+i}\Big),$$
such that
\begin{eqnarray}
\ln\left(1 + \frac{1}{m\varphi_{m} + i}\right) - \ln\left(1 +
\frac{1}{(m+1)\varphi_{m} + i}\right) = \frac{1}{1 + \xi_{i}} X_{i}.
\end{eqnarray}
Since for all $i=0,1,\cdots,m-1$,
$$\xi_{i}>\frac{1}{(m+1)\varphi_{m}+i}>\frac{1}{(m+1)\varphi_{m}+m}>\xi,$$
it follows that
\begin{eqnarray}
\frac{1}{1 + \xi} X_{i} > \frac{1}{1 + \xi_{i}} X_{i}.
\end{eqnarray}
From equations (10), (11), and (12), we obtain
\begin{eqnarray*}
&& \ln(1 + X_{0} + X_{1} + \cdots + X_{m-1}) \nonumber \\
&>&\sum^{m-1}_{i=0} \left[ \ln\left(1 + \frac{1}{m\varphi_{m} +
i}\right) - \ln\left(1 + \frac{1}{(m+1)\varphi_{m} + i}\right)
\right] \nonumber \\
&=& \sum^{m-1}_{i=0} \ln \left[ \frac{m\varphi_{m} + i +
1}{m\varphi_{m} + i} \cdot \frac{(m+1)\varphi_{m} +
i}{(m+1)\varphi_{m} + i + 1} \right] \nonumber \\
&=& \ln\prod^{m-1}_{i=0} \left[ \frac{m\varphi_{m} + i +
1}{m\varphi_{m} + i} \cdot \frac{(m+1)\varphi_{m} +
i}{(m+1)\varphi_{m} + i + 1} \right] \nonumber \\
&=& \ln
\frac{(m+1)(\varphi_{m} + 1)}{(m+1)(\varphi_{m} + 1) - 1}.
\end{eqnarray*}
Meanwhile,
\begin{eqnarray*}
\ln \frac{(m+1)\varphi_{m} + 1}{(m+1)(\varphi_{m} + 1) - 1} - \ln(1
+ Z_{m}) = \ln \frac{(m+1)(\varphi_{m} + 1)}{(m+1)(\varphi_{m} + 1)
- 1} - \ln\left(1 + \frac{1}{(m+1)\varphi_{m} + m}\right) = 0.
\end{eqnarray*}
Therefore,
\begin{eqnarray*}
\ln(1 + X_{0} + X_{1} + \cdots + X_{m-1}) - \ln(1 + Z_{m}) > 0,
\end{eqnarray*}
which contradicts (9). This proves inequality (8).

Since
\begin{eqnarray*}
\frac{1}{m\varphi_{m}} + \frac{1}{m\varphi_{m} + 1} + \cdots +
\frac{1}{m\varphi_{m} + m - 1} = 1
\end{eqnarray*}
and
\begin{eqnarray*}
\frac{1}{(m+1)\varphi_{m+1}} + \frac{1}{(m+1)\varphi_{m+1} + 1} +
\cdots + \frac{1}{(m+1)\varphi_{m+1} + m} = 1,
\end{eqnarray*}
inequality (8) implies that $\varphi_{m+1}<\varphi_{m}$. This proves
that the sequence $\{\varphi_{m}\}$ is strictly decreasing with
respect to $m$.

From the proof of Theorem 2.3, we know that the sequence
$\{\varphi_{m}\}$ has a lower bound $\frac{1}{e-1}$. By the monotone
convergence theorem, the sequence $\{\varphi_{m}\}$ converges. Let
$\lim_{m\rightarrow \infty}\varphi_{m}=\varphi_{0}$. Now, denote the
left-hand side of equation (7) as
\begin{eqnarray*}
G(m)=\frac{1}{m\varphi_{m}}+\frac{1}{m\varphi_{m}+1}+\cdots+\frac{1}{m\varphi_{m}+m-1},
\end{eqnarray*}
then for $m>2$, by the logarithmic inequality
$$\frac{1}{1+x}<\ln(1+\frac{1}{x})<\frac{1}{x}~~\hbox{for}~~x\geq 1,$$ we
have
\begin{eqnarray*}
G(m)&<&\ln\Big(1+\frac{1}{m\varphi_{m}-1}\Big)+\ln\Big(1+\frac{1}{m\varphi_{m}}\Big)+\cdots+
\ln\Big(1+\frac{1}{m\varphi_{m}+m-2}\Big)\\
&=&\ln\big(1+\frac{1}{\varphi_{m}-\frac{1}{m}}\big)
\end{eqnarray*}
and
\begin{eqnarray*}
G(m)&>&\ln\Big(1+\frac{1}{m\varphi_{m}}\Big)+\ln\Big(1+\frac{1}{m\varphi_{m}+1}\Big)+\cdots+
\ln\Big(1+\frac{1}{m\varphi_{m}+m-1}\Big)\\
&=&\ln\big(1+\frac{1}{\varphi_{m}}\big).
\end{eqnarray*}
Moreover,
\begin{eqnarray*}
\lim_{m\rightarrow
\infty}\ln\big(1+\frac{1}{\varphi_{m}-\frac{1}{m}}\big)=\lim_{m\rightarrow
\infty}\ln\big(1+\frac{1}{\varphi_{m}}\big)=\ln\big(1+\frac{1}{\varphi_{0}}\big).
\end{eqnarray*}
Thus, we have
\begin{eqnarray*}
\lim_{m\rightarrow \infty}G(m)=\ln\big(1+\frac{1}{\varphi_{0}}\big).
\end{eqnarray*}

Taking the limit of both sides of equation (7) as $m \to \infty$, we
solve for
$\lim_{m\rightarrow\infty}\varphi=\varphi_{0}=\frac{1}{e-1}$.
Furthermore, it is easy to derive that
$\lim_{m\rightarrow\infty}\alpha=\frac{e}{e-1}.$ The proof is
complete.

By Vieta's theorem and the definition of the Stirling numbers of the
first kind $s(m,k)$:
$$\prod^{m-1}_{k=0}(x-k)=\sum^{m}_{k=0}(-1)^{n-k}s(m,k)x^{k},$$ we
obtain the following Lemma 2.2.

\vspace{1.5mm}\noindent{\bf Lemma 2.2}~ 1)~ {\it Let
$I_{m}=\{1,2,\cdots,m-1\},$ and let $i_{1:k} \in I_{m}$ denote $i_1,
i_2, \dots, i_k \in I_{m}$ with $i_{1}<i_{2}<\cdots<i_{k}$, then
\begin{eqnarray}
\prod^{m-1}_{k=0}(x+k)=x^{m}+\Big(\sum^{m-1}_{i=1}{i}\Big)x^{m-1}+
              \Big(\sum_{i_{1:2}\in I_{m}}\prod^{2}_{j=1}i_{j}\Big)x^{m-2}+\cdots
               +\Big(\sum_{i_{1:m-2}\in I_{m}}\prod^{m-2}_{j=1}i_{j}\Big)x^{2}+
                (m-1)!~x.
 \end{eqnarray}

\begin{eqnarray*}
2)\hspace{0.5cm}
\prod^{m-1}_{k=0}(x+k)=s(m,m)x^{m}+s(m,m-1)x^{m-1}+s(m,m-2)x^{m-2}+\cdots+s(m,2)x^{2}+s(m,1)x.
 \end{eqnarray*}}

\vspace{1.5mm}\noindent{\bf Theorem  2.5}
 {\it The $m$-dimensional harmonic-type golden section equation(6)
\begin{eqnarray*}
\frac{1}{x}+\frac{1}{x+1}+\frac{1}{x+2}+\cdots+\frac{1}{x+m-1}=1
\end{eqnarray*}
and
\begin{eqnarray}
&&x^{m}+\big(\sum^{m-1}_{i=1}i-m\big)x^{m-1}+
 \Big(\sum_{i_{1:2}\in I_{m}}\prod^{2}_{j=1}i_{j}-(m-1)\sum^{m-1}_{i=1}i\Big)x^{m-2}+\cdots+\nonumber\\
     &&\Big(\sum_{i_{1:m-2}\in I_{m}}\prod^{m-2}_{j=1}i_{j}-3
\sum_{i_{1:m-3}\in I_{m}}\prod^{m-3}_{j=1}i_{j}\Big)x^{2}
+\Big((m-1)!-2\sum_{i_{1:m-2}\in
I_{m}}\prod^{m-2}_{j=1}i_{j}\Big)x-(m-1)!=0
\end{eqnarray}
have the same solutions.}

\vspace{1.5mm}{\bf Proof} When $x\neq 0, -1, -2, \cdots, -(m-1)$,
equation (6) can be rewritten as
\begin{eqnarray}
 P_{0}(x)+P_{1}(x)+P_{2}(x)+\cdots+P_{m-1}(x)=\prod^{m-1}_{k=0}(x+k),
 \end{eqnarray}
where $P_{k}(x)=x(x-1)\cdots[x-(k-1)][x-(k+1)]\cdots[x-(m-1)],
(k=0,1,2,\cdots,m-1). $

Since the derivative of the product is
$$\Big(\prod^{m-1}_{k=0}(x+k)\Big)'=\sum^{m-1}_{k=0}P_{k}(x)$$
equation (15) is equivalent to
\begin{eqnarray}\Big(\prod^{m-1}_{k=0}(x+k)\Big)'=\prod^{m-1}_{k=0}(x+k).\end{eqnarray}
Furthermore, when substituting $x=0, -1, -2, \cdots, -(m-1)$, the
left-hand side of (15), $\sum^{m-1}_{k=0}P_{k}(x)\neq0,$ is never
zero, while the right-hand side, $\prod^{m-1}_{k=0}(x+k)=0$, becomes
zero. Therefore, equation (6) and equation (16) have the same
solutions.

By Lemma 2.2, we have
\begin{eqnarray}
\Big(\prod^{m-1}_{k=0}(x+k)\Big)'&=&mx^{m-1}+(m-1)\Big(\sum^{m-1}_{i=1}{i}\Big)x^{m-2}+
              (m-2)\Big(\sum_{i_{1:2}\in I_{m}}\prod^{2}_{j=1}i_{j}\Big)x^{m-3}+\cdots\nonumber\\
               &+&2\Big(\sum_{i_{1:m-2}\in I_{m}}\prod^{m-2}_{j=1}i_{j}\Big)x+
                (m-1)! .
 \end{eqnarray}
Substituting (13) and (17) into equation (16) and simplifying, we
conclude that equation (6) is indeed equivalent to equation (14).

The proof is complete.

Using the Stirling numbers of the first kind $s(m,k)$ as
coefficients and noting that $s(m,m)=1$ and $s(m,1)=(m-1)!$,
equation (14) can alternatively be expressed as:
$$x^{m}+[s(m,m-1)-m]x^{m-1}+\cdots+[s(m,2)-3s(m,3)]x^{2}+[s(m,1)-2s(m,2)]x-(m-1)!=0.$$

For example, the equation
$$\frac{1}{x}+\frac{1}{x+1}+\frac{1}{x+2}+\frac{1}{x+3}+\frac{1}{x+4}+\frac{1}{x+5}=1$$
takes the following form when expressed with Stirling number
coefficients:
$$x^{6}+9x^{5}+10x^{4}-115x^{3}-401x^{2}-428x-120=0.$$

\section{ Properties of Higher-dimensional Golden Section}

To better characterize the properties of the higher-dimensional
golden section, we supplement the definition of the $k$-generalized
Fibonacci numbers sequence as defined by Miles [43] as follows:

\vspace{1.5mm} \noindent\textbf{Definition 3.1} Let $m\geq2$ be a
positive integer, and let
$$F^{(m)}_{0}=F^{(m)}_{1}=\cdots=F^{(m)}_{m-3}=F^{(m)}_{m-2}=0,
  F^{(m)}_{m-1}=1, $$
$$F^{(m)}_{n}=F^{(m)}_{n-1}+F^{(m)}_{n-2}+\cdots+F^{(m)}_{n-m+1}+F^{(m)}_{n-m}~~~(n=m, m+1, \cdots).$$
$$F^{(m)}_{n}=F^{(m)}_{n+m}-F^{(m)}_{n+m-1}-\cdots-F^{(m)}_{n+2}-F^{(m)}_{n+1}~~~(n=-1, -2, \cdots).$$
Then the bi-infinite sequence $\{F^{(m)}_{n}\}_{n\in \mathbb{Z}} $
  is called the $m$-generalized Fibonacci sequence or $m$-dimensional Fibonacci sequence, where $m$
  is its dimension.

When $m=2$, this reduces to the classical Fibonacci sequence
extended to include terms with negative indices, which aligns with
the extended Fibonacci sequence by Shi [47]. Moreover, Shi [47]
proved that $F^{(2)}_{-(2k-1)}=F^{(2)}_{2k-1},
F^{(2)}_{-2k}=-F^{(2)}_{2k}, k \in \mathbb{Z}^{+} $, i.e.,
$F^{(2)}_{-n}=(-1)^{n-1}F^{(2)}_{n}.$

For $m>2$, calculating the terms of the $m$-generalized Fibonacci
sequence using the recursive formula in Definition 3.1 is relatively
slow. However, the original recursive formula easily yields the
following recursive relations.

\vspace{1.5mm}\noindent\textbf{Theorem 3.1} For the $m$-generalized
Fibonacci sequence $\{F^{(m)}_{n}\}$, $\forall ~n\in \mathbb{Z}$,
$$F^{(m)}_{n}=2F^{(m)}_{n-1}-F^{(m)}_{n-m-1},~ F^{(m)}_{n}=2F^{(m)}_{n+m}-F^{(m)}_{n+m+1}.$$

The formula $F^{(m)}_{n}=2F^{(m)}_{n+m}-F^{(m)}_{n+m+1}$
  is primarily used to compute the values of the $m$-generalized Fibonacci sequence
  $\{F^{(m)}_{n}\}$ for negative indices. To this end, we provide the initial
values for calculating the negative-indexed terms of
$\{F^{(m)}_{n}\}$: $$F^{(m)}_{0}=0, ~F^{(m)}_{-1}=1,~
F^{(m)}_{-2}=-1,~~
F^{(m)}_{-3}=\cdots=F^{(m)}_{-m+1}=F^{(m)}_{-m}=0.$$

\vspace{1.5mm}\noindent\textbf{Corollary 3.1} For the
$m$-generalized Fibonacci sequence $\{F^{(m)}_{n}\}$, $n\in
\mathbb{Z},$ and $\forall ~k\in \mathbb{Z}^{+}$, we have
$$F^{(m)}_{n}=2^{k}F^{(m)}_{n-k}-2^{k-1}F^{(m)}_{n-k-m}-2^{k-2}F^{(m)}_{n-(k-1)-m}-\cdots-2^{2}
F^{(m)}_{n-3-m}-2F^{(m)}_{n-2-m}-F^{(m)}_{n-1-m},$$
$$F^{(m)}_{n}=2^{k}F^{(m)}_{n+km}-2^{k-1}F^{(m)}_{n+km+1}-2^{k-2}F^{(m)}_{n+(k-1)m+1}-\cdots-2^{2}
F^{(m)}_{n+3m+1}-2F^{(m)}_{n+2m+1}-F^{(m)}_{n+m+1}.$$

\subsection*{3.1\ Properties of Geometric-type
Higher-dimensional Golden Section}

By the definition of geometric-type higher-dimensional golden
section, the $m$-dimensional geometric-type golden section has the
following properties.

\subsection*{3.1.1\ Basic Relationships among Section Points,
Section Intervals, and the Section Ratio}

~~~~$\bullet$~The relationship between section points and section
intervals:

$$\delta_{k}=\tau_{k}-\tau_{k-1}; ~ \hbox{conversely,}~
\tau_{k}=\sum_{i=1}^{k}\delta_{i}, ~\hbox{where}~k=1, 2, \cdots, m;
\tau_{0}=0, \tau_{m}=1.$$

$\bullet$~ The $m$-dimensional geometric-type golden section ratio
satisfies
$$\lambda=\frac{1}{\delta_{1}}=\frac{\delta_{1}}{\delta_{2}}=\cdots=\frac{\delta_{m-1}}
{\delta_{m}};~ \hbox{conversely,}~~ \lambda^{-1}=\delta_{1},
\lambda^{-2}=\delta_{2}, \cdots, \lambda^{-m}=\delta_{m}.$$

$\bullet$~ The section points of the $m$-dimensional geometric-type
golden section can be expressed in terms of the section ratio as:
$$\tau_{1}=\lambda^{-1},\tau_{2}=\lambda^{-1}+\lambda^{-2},\cdots,\tau_{m-1}=\lambda^{-1}+
\lambda^{-2}+\cdots+\lambda^{-(m-1)},
\tau_{m}=\lambda^{-1}+\lambda^{-2}+\cdots+\lambda^{-m}=1.$$

\subsection*{3.1.2\ The Relationship between the Section Ratio
$\lambda$ and $\{F^{(m)}_{n}\}$}

We are familiar with the fact that the classical (2-dimensional)
golden section has a `beautiful' relationship with the Fibonacci
sequence:
$(1+\sqrt{5})/2=\lim_{n\rightarrow\infty}{F^{(2)}_{n+1}}/{F^{(2)}_{n}}.
$ Delightfully, the $m$-dimensional geometric-type golden section
also has a similar relationship with the $m$-generalized Fibonacci
sequence:

\vspace{1.5mm}\noindent{\bf Theorem 3.2} Let $\lambda_{1},
\lambda_{2}, \cdots, \lambda_{m-1}, \lambda$ be the distinct roots
of the equation $x^{m}-x^{m-1}-\cdots-x-1=0$ ~$(m\in \mathbb{Z}^{+}$
and $m\geq2)$. For any integer $n$, we have
\begin{eqnarray*}
 \lambda=\lim_{n\rightarrow\infty}\frac{F^{(m)}_{n+1}}{F^{(m)}_{n}}.
\end{eqnarray*}

\vspace{1.5mm}{\bf Proof} Since the recurrence relation of the
$m$-generalized Fibonacci sequence $\{F^{(m)}_{n}\}$ is
$$F^{(m)}_{n}=F^{(m)}_{n-1}+F^{(m)}_{n-2}+\cdots+F^{(m)}_{n-m+1}+F^{(m)}_{n-m}~~~(n=1, 2, \cdots),$$
the corresponding characteristic equation is precisely the
$m$-dimensional geometric-type golden section equation
$$x^{m}-x^{m-1}-\cdots-x^{2}-x-1=0~(m\geq2).$$
By Theorem 2.1, the $m$-dimensional geometric-type golden section
ratio $\lambda$ is the unique largest positive real root of this
equation. Let
\begin{eqnarray*}
c_{i}=(\lambda_{i}-\lambda)\prod_{j\neq
i}^{m-1}(\lambda_{i}-\lambda_{j})~ (i=1,2,\cdots,m-1), ~
c_{m}=\prod_{j=1}^{m-1}(\lambda-\lambda_{j}).
\end{eqnarray*}

From Theorem 2.2 and the general term formula in reference (Sun
[50], the general term of the $m$-generalized Fibonacci sequence
$\{F^{(m)}_{n}\}$ is
$$
  {F^{(m)}_{n}}=c_{1}(\lambda_{1})^{n+m-1}+c_{2}(\lambda_{2})^{n+m-1}+\cdots+
c_{m-1}(\lambda_{m-1})^{n+m-1}+c_{m}(\lambda)^{n+m-1} . $$ From
Theorem 2.1 and Theorem 2.2, for any $i\in \{1,2,\cdots,m-1\}$, we
have $\mid\lambda_{i}\mid<\lambda$. Therefore,
\begin{eqnarray*}
\lim_{n\rightarrow\infty}\frac{F^{(m)}_{n+1}}{F^{(m)}_{n}}&=&
\lim_{n\rightarrow\infty}\frac{ c_{1}(\lambda_{1})^{n+m}+\cdots+
c_{m-1}(\lambda_{m-1})^{n+m}+c_{m}(\lambda)^{n+m}
}{c_{1}(\lambda_{1})^{n+m-1}+\cdots+
c_{m-1}(\lambda_{m-1})^{n+m-1}+c_{m}(\lambda)^{n+m-1}}\\
&=&\lim_{n\rightarrow\infty}\frac{
c_{1}\lambda_{1}\left(\frac{\lambda_{1}}{\lambda}\right)^{n+m-1}+\cdots+
c_{m-1}\lambda_{m-1}\left(\frac{\lambda_{m-1}}{\lambda}\right)^{n+m-1}+c_{m}\lambda}
{c_{1}\left(\frac{\lambda_{1}}{\lambda}\right)^{n+m-1}+\cdots+
c_{m-1}\left(\frac{\lambda_{m-1}}{\lambda}\right)^{n+m-1}+c_{m}}\\
&=&\lambda.
\end{eqnarray*}
The proof is complete.

 Theorem 3.2 shows that the geometric-type higher-dimensional golden section
ratio $\lambda$ can be approximated by the ratio of consecutive
terms $F_{n+1}^{(m)}/F_{n}^{(m)}$
  in the $m$-generalized Fibonacci sequence $F_{n}^{(m)}$. This provides a method for
  numerically solving the largest root of the geometric-type higher-dimensional golden section equation
$x^{m}-x^{m-1}-\cdots-x^{2}-x-1=0$.

For example, consider the 3-generalized Fibonacci sequence
$\{F_{n}^{(3)}\}$. When $n= 11, 12, 13, 14, 15, 16$, the
corresponding terms are: $274, ~504,~927,~1705,~3136,~5768.274, 504,
927, 1705, 3136, 5768.$ The ratios of consecutive terms are
approximately:
\begin{eqnarray*}
1.839416058,~1.839285714,~1.839266451,~1.839296188,~1.839285714.
\end{eqnarray*}
Using $F_{16}^{(3)}/F_{15}^{(3)}\approx 1.839285714$ as an
approximation for the 3-dimensional geometric-type golden section
ratio $\lambda$, the resulting error compared to the previously
derived formula is $-1.04\times^{-6}.$

\subsection*{3.1.3\ Generating Sequence of the Section Ratio
$\lambda$}

In 1876, Lucas discovered that any linear combination of the
classical golden section ratio $\lambda=(\sqrt{5}+1)/2$ and
$-\lambda^{-1}=(1-\sqrt{5})/2$ satisfies the recurrence relation
$F_{n+1}^{(2)}=F_{n}^{(2)}+F_{n-1}^{(2)}$ (Wu [56]). Reference (Shi
[47]) defined and studied the golden sequence using powers of
$\lambda=(\sqrt{5}+1)/2$ as follows:
$$\{G^{\pm n}\}_{n=0}^{+\infty}:~ \cdots, (\frac{\sqrt{5}+1}{2})^{-n},
\cdots, (\frac{\sqrt{5}+1}{2})^{-1}, ~1, ~\frac{\sqrt{5}+1}{2},
\cdots, (\frac{\sqrt{5}+1}{2})^{n}, \cdots.$$ They proved the
relationship between this sequence and the Fibonacci sequence
$\{F^{(2)}_{n}\}$: $G^{n}=F^{(2)}_{n+1}+G^{-1}F^{(2)}_{n}.$
Similarly, we define the following generating sequence for the
geometric-type higher-dimensional golden section ratio.

\vspace{1.5mm} \noindent{\bf Definition 3.2} Let $\lambda$ be the
$m$-dimensional geometric-type higher-dimensional golden section
ratio satisfying the equation
$$\frac{1}{\lambda}+\frac{1}{\lambda^{2}}+\cdots+\frac{1}{\lambda^{m}}=1.$$
Then the sequence $$\lambda^{-m}, \cdots, \lambda^{-2},
\lambda^{-1},~1,~ \lambda, \lambda^{2}, \cdots, \lambda^{m}, \cdots,
\lambda^{n },  \cdots, $$ is called the generating sequence of the
$m$-dimensional geometric-type higher-dimensional golden section
ratio, denoted as $\{\lambda^{n}\}_{n=-m}^{+\infty}$.

When $m=2$, $\lambda=(\sqrt{5}+1)/2$. The generating sequence of the
2-dimensional golden section ratio we defined:
$$(\frac{\sqrt{5}+1}{2})^{-2}, (\frac{\sqrt{5}+1}{2})^{-1}, ~1, ~\frac{\sqrt{5}+1}{2},
\cdots, (\frac{\sqrt{5}+1}{2})^{n}, \cdots$$ represents a special
case of the golden sequence starting from $n=-2$ as described in
reference (Shi [47]). This sequence satisfies the relation:
$\lambda^{n}=F^{(2)}_{n+1}+\lambda^{-1}F^{(2)}_{n}.$

As $m$ approaches infinity, the infinite-dimensional geometric-type
golden section ratio converges to $\lambda=2$. The corresponding
generating sequence becomes:
$$\cdots,~\frac{1}{2^{n}}, ~\cdots, ~\frac{1}{2^{2}}, ~\frac{1}{2}, ~1,~ 2,~ 2^{2},~\cdots,~2^{n},~\cdots.$$
The left portion of this generating sequence forms the series:
$1,\frac{1}{2}, ~\frac{1}{2^{2}}, ~\cdots, \frac{1}{2^{n}},
~\cdots$. This demonstrates an interesting connection: the bisection
search algorithm used for finding numerical solutions to nonlinear
algebraic equations is fundamentally related to the
infinite-dimensional geometric-type golden section ratio concept.

Below, let's look at an example of the application of the generating
sequence of the infinite-dimensional geometric-type golden section
ratio in the field of robot algorithms. Suppose that the robot needs
to find some distinguished point on a line. Assume the robot neither
knows the actual distance to the point nor its left or right. What
is the minimum number of steps it must make to find the point?
Baeza-Yates et al. [2] proposed the linear spiral search: Execute
cycles of steps where the function determining the number of steps
to walk before the $i^{th}$  turn from the origin is
$f(i)=2^{i}~\forall i\geq 1.$ They also proved that the linear
spiral algorithm is the optimal online algorithm, with a competitive
ratio of 9. The number of steps taken in each search of this optimal
algorithm is precisely the sequence on the right of the
infinite-dimensional geometric-type golden section ratio:  $ ~1, ~
2,~ 2^{2},~ \cdots, ~2^{n}, ~\cdots.$

When $m=3$, the three-dimensional geometric-type golden section
ratio
$$\lambda=\frac{1+\sqrt[3]{19+3\sqrt{33}}+\sqrt[3]{19-3\sqrt{33}}}{3}
\approx1.839286755 $$ satisfies the equation
$$\frac{1}{\lambda}+\frac{1}{\lambda^{2}}+\frac{1}{\lambda^{3}}=1.$$
The generating sequence of the three-dimensional geometric golden
section ratio we constructed $$\lambda^{-3}, \lambda^{-2},
\lambda^{-1}, ~1,~ \lambda, \lambda^{2}, \lambda^{3}, \cdots,
\lambda^{n}, \cdots$$ exhibits characteristics of the 3-generalized
Fibonacci sequence $\{F_{n}^{(3)}\}$:
$\lambda^{n}=\lambda^{n-1}+\lambda^{n-2}+\lambda^{n-3}~(n\geq0)$,
and can be expressed using the $3$-generalized Fibonacci sequence
$\{F_{n}^{(3)}\}$:
$$\lambda^{n}=F_{n+2}^{(3)}+(\lambda^{-1}+\lambda^{-2})F_{n+1}^{(3)}
+\lambda^{-1}F_{n}^{(3)}, ~(n\geq -3) .$$

\vspace{1.5mm}{\bf Proof} When $n=-3$, since $ F_{-1}^{(3)}=1,
F_{-2}^{(3)}=-1, F_{-3}^{(3)}=0, $ we have:
$$\lambda^
{-2}=F_{0}^{(3)}+(\lambda^{-1}+\lambda^{-2})F_{-1}^{(3)}+\lambda^{-1}F_{-2}^{(3)}.$$
When $n=-2$, since $F_{0}^{(3)}=0, F_{-1}^{(3)}=1, F_{-2}^{(3)}=-1,$
we have: $$\lambda^
{-2}=F_{0}^{(3)}+(\lambda^{-1}+\lambda^{-2})F_{-1}^{(3)}+\lambda^{-1}F_{-2}^{(3)}.$$
When $n=-1$, since $F_{1}^{(3)}=0, F_{0}^{(3)}=0,$ and
$F_{-1}^{(3)}=1,$ we have: $$\lambda^ {-1}=F_{1}^{(3)}
+(\lambda^{-1}+\lambda^{-2})F_{0}^{(3)}+\lambda^{-1}F_{-1}^{(3)}.$$
Assume the proposition holds for all $0\leq n\leq k~(k\geq -1)$.
Then, for $n=k+1$, since $\lambda^{k+1}=\lambda^{k}
+\lambda^{k-1}+\lambda^{k-2},$ by the induction hypothesis, we have:
$$\lambda^{k}=F_{k+2}^{(3)}+(\lambda^{-1}+\lambda^{-2})F_{k+1}^{(3)}
+\lambda^{-1}F_{k}^{(3)},$$
$$\lambda^{k-1}=F_{k+1}^{(3)}+(\lambda^{-1}+\lambda^{-2})F_{k}^{(3)}
+\lambda^{-1}F_{k-1}^{(3)},$$
$$\lambda^{k-2}=F_{k}^{(3)}+(\lambda^{-1}+\lambda^{-2})F_{k-1}^{(3)}
+\lambda^{-1}F_{k-2 }^{(3)}.$$ Adding these three equations yields:
$$\lambda^{k+1}=F_{k+3}^{(3)}+(\lambda^{-1}+\lambda^{-2})F_{k+2}^{(3)}
+\lambda^{-1}F_{k+1 }^{(3)}.$$ By the principle of mathematical
induction, the proposition holds for all integers $n\geq -3$.

The above conclusions and it's proofs similarly apply to the
following general case:

 \vspace{1.5mm}  \noindent{\bf Theorem ~3.3} The
generating sequence $\{\lambda^{n}\}_{n=-m}^{+\infty}$ of the
$m$-dimensional geometric-type golden section ratio can be expressed
in terms of the $m$-generalized Fibonacci sequence $\{F_{n}^{(m
)}\}$ as:
$$\lambda^{n}=F_{n+m-1}^{(m)}+(\lambda^{-1}+\lambda^{-2}+\cdots+\lambda^{-(m-1)})F_{n+m-2}^{(m)}
+(\lambda^{-1}+\lambda^{-2}+\cdots+\lambda^{-(m-2)})F_{n+m-3}^{(m)}
+\cdots+\lambda^{-1}F_{n}^{(m)}.$$

From Section 3.1.1 and Theorem ~3.3, the generating sequence
$\{\lambda^{n}\}_{n=-m}^{+\infty}$
  can alternatively be represented using the section intervals
$\delta_k$
  and the $m$-generalized Fibonacci sequence:
$$\lambda^{n}=F_{n+m-1}^{(m)}+(\delta_{1}+\delta_{2}+\cdots+\delta_{m-1})F_{n+m-2}^{(m)}+
(\delta_{1}+\delta_{2}+\cdots+\delta_{m-2})F_{n+m-3}^{(m)}
+\cdots+\delta_{1}F_{n}^{(m)}. $$
And the general term of the
generating sequence for the $m$-dimensional geometric-type golden
section ratio can alternatively be expressed using these golden
section points as:
$$\lambda^{n}=F_{n+m-1}^{(m)}\tau_{m}+F_{n+m-2}^{(m)}\tau_{m-1}+F_{n+m-3}^{(m)}
\tau_{m-2}\cdots+F_{n+1}^{(m)}\tau_{2}+ +
F_{n}^{(m)}\tau_{1}~~(n\geq -m) .$$

\subsection*{3.2\ Properties of Harmonic-type
Higher-dimensional Golden Section}

According to the definition of the harmonic-type higher-dimensional
golden section, the $m$-dimensional harmonic-type golden section
$T^{(m)}: t_{m}, t_{m-1}, \cdots,t_{2},  t_{1}$ has the following
properties.

\subsection*{3.2.1 Basic Relationships among Section Points,
Section Intervals, and the Section Basis}

~~~~~$\bullet$~$\Delta_{k}=t_{k}-t_{k-1}; $ conversely,
$t_{k}=\sum_{i=1}^{k}\Delta_{i}$,~ where $k=1, 2, \cdots, m;
t_{0}=0, t_{m}=1.$

$\bullet$~The ratio of the lengths of any two adjacent subintervals:
$\frac{\Delta_{k}}{\Delta_{k+1}}=1+\Delta_{k}$, where $k=1, 2,
\cdots, m-1$.

$\bullet$~The lengths of the subintervals can be expressed in terms
of the section basis as: $$\Delta_{1}=\frac{1}{\omega},
\Delta_{2}=\frac{1}{\omega+1}, \cdots,
\Delta_{m}=\frac{1}{\omega+m-1}.$$

$\bullet$~The coordinates of the section points can be expressed in
terms of the section basis as:
$$t_{k}=\sum_{i=1}^{k}\frac{1}{\omega+i-1},~\hbox{where}~k=1, 2, \cdots,
m;  t_{m}=1.$$

$\bullet$~ The sum of the right endpoints of all subintervals is the
root of the equation (6).

\vspace{1.5mm}{\bf Proof}~The
equation$$\frac{1}{\omega}+\frac{1}{\omega+1}+\cdots+
  \frac{1}{\omega+m-1}=1$$
is multiplied on both sides by $\omega+m$ to obtain
$$\frac{\omega+m}{\omega}+\frac{\omega+m}{\omega+1}+\cdots+\frac{\omega+m}{\omega+m-1}=\omega+m.$$
Rearranging this gives
$$\sum^{m-1}_{k=0}(\frac{1}{\omega}+\frac{1}{\omega+1}+\cdots+\frac{1}{\omega+k})=\omega. $$
That is, $\sum_{k=1}^{m}t_{k}=\omega .$

\subsection*{3.2.2\  The Generating Sequence of the
$m$-dimensional Harmonic-type Golden Section
$T^{(m)}$}

For the classical 2-dimensional golden section $T^{(2)}:t_{2}=1,
t_{1}=\frac{\sqrt{5}-1}{2},$ if we define $T^{(2)}_{-2}=t_{2}-t_{1},
T^{(2)}_{-1}=t_{1}, T^{(2)}_{0}=t_{2},$ then the 2-dimensional
harmonic-type golden section equation
$\frac{1}{\omega}+\frac{1}{\omega+1}=1$ can be expressed as
$$T^{(2)}_{-2}+T^{(2)}_{-1}=T^{(2)}_{0}.$$
This exhibits a Fibonacci-like property, demonstrating the close
relationship between the golden section and the Fibonacci sequence.
If we take $T^{(2)}_{-2}=t_{2}-t_{1}, T^{(2)}_{-1}=t_{1}$ as the
initial values and recursively apply
$T^{(2)}_{n}=T^{(2)}_{n-1}+T^{(2)}_{n-2}$, we can generate a
sequence composed of the section points $t_{2}=1,
t_{1}=\frac{\sqrt{5}-1}{2}$ of $T^{(2)}$. More generally, we define
as follows.

\vspace{1.5mm}\noindent{\bf Definition 3.3} For the $m$-dimensional
harmonic-type golden section $T^{(m)}: t_{m},~ t_{m-1},
~\cdots,~t_{2}, ~ t_{1}$, we define the
sequence$$T^{(m)}_{-m}=t_{m}-t_{m-1},
T^{(m)}_{-m+1}=t_{m-1}-t_{m-2}, \cdots, T^{(m)}_{-2}=t_{2}-t_{1},
T^{(m)}_{-1}=t_{1}, T^{(m)}_{0}=t_{m},$$
$$T^{(m)}_{n}=T^{(m)}_{n-1}+T^{(m)}_{n-2}+\cdots+
T^{(m)}_{n-m},~~(n\geq0), $$ as the generating sequence of the
$m$-dimensional harmonic-type golden section $T^{(m)}$, denoted by
$\{T^{(m)}_{n}\}_{n=-m}^{\infty}$.

It is not difficult to prove that the generating sequence
$\{T^{(m)}_{n}\}_{n=-m}^{\infty}$ of the $m$-dimensional harmonic
golden section $T^{(m)}$ and the $m$-generalized Fibonacci sequence
$\{F^{(m)}_{n}\}$ satisfy the following relationship:

\vspace{1.5mm}\noindent{\bf Theorem 3.4}~ For any integer $n\geq
-m$, we have
$$T^{(m)}_{n}=F^{(m)}_{n+m-1}t_{m}+F^{(m)}_{n+m-2}t_{m-1}+F^{(m)}_{n+m-3}t_{m-2}+
\cdots+F^{(m)}_{n+1}t_{2}+F^{(m)}_{n}t_{1},$$ where $t_{k}
~(k=1,2,\dots, m)$ are the section points of the $m$-dimensional
harmonic-type golden section.

From Theorem 3.4, when $m=2$, the generating sequence of the
classical (2-dimensional) golden section $T^{(2)}$ is
$$\{T^{(2)}_{n}\}_{n=-2}^{\infty}: ~t_{2}-t_{1}, ~~t_{1},
~~t_{2},~~t_{2}+t_{1},~~2t_{2}+t_{1}, ~\cdots, ~F^{(2)}_{n+1}t_{2}+
F^{(2)}_{n}t_{1},~ \cdots. $$ Since $t_{2}=1,$ and $
t_{1}=(\sqrt{5}-1)/2=\omega^{-1}=\lambda^{-1},$  the generating
sequence of the 2-dimensional harmonic-type golden section $T^{(2)}$
can also be expressed as
$$\{T^{(2)}_{n}\}_{n=-2}^{\infty}:
~\omega^{-2}, ~~\omega^{-1},
~~\omega^{0},~~\omega^{1},~~\omega^{2},~~\cdots, ~\omega^{n},~
\cdots. $$ This is entirely consistent with the generating sequence
of the classical (2-dimensional) golden section when viewed as a
geometric-type golden section.

\section{ The Fibonacci Matrix of Harmonic-type
Higher-dimensional Golden Section}

In Sections 3.1 and 3.2, we discussed the generating sequences
$\{\lambda^{n}\}_{n=-m}^{+\infty}$ of the geometric-type
higher-dimensional golden section ratio $\lambda$ and the generating
sequences $\{T^{(m)}_{n}\}_{n=-m}^{\infty}$ of the harmonic-type
higher-dimensional golden section (point) $T^{(m)}$, respectively.
For the geometric-type higher-dimensional golden section, since
$\delta_{k}$ decreases geometrically with ratio $\lambda$, the
general term of its generating sequence
$\{\lambda^{n}\}_{n=-m}^{+\infty}$ can be expressed as a linear
combination of the section points $\tau_{m},~ \tau_{m-1},
~\cdots,~\tau_{2}, ~ \tau_{1}$ and the $m$-generalized Fibonacci
sequence:
\begin{eqnarray*}
\lambda^{n}=F_{n+m-1}^{(m)}\tau_{m}+F_{n+m-2}^{(m)}\tau_{m-1}+\cdots+F_{n+1}^{(m)}\tau_{2}+F_{n}^{(m)}\tau_{1}.
\end{eqnarray*}

However, for the harmonic-type higher-dimensional golden section,
$\Delta_{k}$ decreases harmonically as $\frac{1}{\omega+k-1}$, and
its section basis $\omega$ cannot be directly expressed as in the
above formula. Nevertheless, the generating sequence
$\{T^{(m)}_{n}\}_{n=-m}^{\infty}$ of the section points can still be
represented as a linear combination of the division points $t_{m},~
t_{m-1}, ~\cdots,~t_{2}, ~ t_{1}$ and the $m$-generalized Fibonacci
sequence:
\begin{eqnarray*}
T^{(m)}_{n}=F^{(m)}_{n+m-1}t_{m}+F^{(m)}_{n+m-2}t_{m-1}+F^{(m)}_{n+m-3}t_{m-2}+
\cdots+F^{(m)}_{n+1}t_{2}+F^{(m)}_{n}t_{1}.
\end{eqnarray*}

Based on this expression, we further study the generating sequence
$\{T^{(m)}_{n}\}_{n=-m}^{\infty}$ of the harmonic-type
higher-dimensional golden section $T^{(m)}$, and introduce the
concepts of the Fibonacci matrix, row-column index system, and
various submatrices for the harmonic-type higher-dimensional golden
section. At the same time, we explore the properties of the
Fibonacci matrix and its submatrices.

\subsection*{4.1 The Fibonacci Matrix of Harmonic-type
Higher-dimensional Golden Section: Row-column Index
System and Submatrix Classification }

To better characterize the $m$-dimensional harmonic-type golden
section $T^{(m)}$, we consider the general form of its generating
sequence.

For the generating sequence $\{T^{(2)}_{n}\}_{n=-2}^{\infty}$ of
$T^{(2)}$, if we let $T^{(2)}_{-n-1}=T^{(2)}_{-n+1}-
T^{(2)}_{-n},~n=2, 3, 4, \cdots$, then the generating sequence of
$T^{(2)}$ can be extended to a bidirectional infinite sequence
$\{T^{(2)}_{n}\}_{-\infty}^{\infty}:$
$$\cdots,~~F^{(2)}_{-n+1}t_{2}+F^{(2)}_{-n}t_{1},~~\cdots, ~~-t_{2}+2t_{1}, ~~t_{2}-t_{1},
~~t_{1}, ~~t_{2},~~t_{2}+t_{1},~~2t_{2}+t_{1}, ~\cdots,
~F^{(2)}_{n+1}t_{2}+F^{(2)}_{n}t_{1},~ \cdots. $$ Thus, the
generating sequence $\{T^{(2)}_{n}\}_{-\infty}^{\infty}$ of
$T^{(2)}$ corresponds to the golden sequence in ([47]): $$\cdots,
~\omega^{-n},~ \cdots,~~\omega^{-2}, ~~\omega^{-1},
~~\omega^{0},~~\omega^{1},~~ \omega^{2},~~\cdots, ~\omega^{n},~
\cdots. $$

When $m \geq 3$, if we disregard the geometric meaning of the
$m$-dimensional harmonic-type golden section $T^{(m)}$ and allow
$T^{(m)}$ to generate negative terms, we can define the generating
sequence of $T^{(m)}$ in the following bidirectional infinite form:
$$T^{(m)}_{-m}=t_{m}-t_{m-1}, T^{(m)}_{-m+1}=t_{m-1}-t_{m-2},\cdots,T^{(m)}_{-2}=t_{2}-t_{1},
T^{(m)}_{-1}=t_{1}, T^{(m)}_{0}=t_{m},$$
$$T^{(m)}_{n}=T^{(m)}_{n-1}+T^{(m)}_{n-2}+\cdots+T^{(m)}_{n-m},~~(n\geq1),$$
$$T^{(m)}_{-n}=T^{(m)}_{-n+m}-T^{(m)}_{-n+m-1}-\cdots-T^{(m)}_{-n+1},~~(n\geq m+1),$$
denoted as $\{T^{(m)}_{n}\}_{-\infty}^{\infty}$. It can be easily
verified that for all $n \in \mathbb{Z}$ and $~m\geq 2$,
$$T^{(m)}_{n}=F^{(m)}_{n+m-1}t_{m}+F^{(m)}_{n+m-2}t_{m-1}+F^{(m)}_{n+m-3}t_{m-2}+\cdots
+F^{(m)}_{n+1}t_{2}+F^{(m)}_{n}t_{1}.$$ For all $n \in \mathbb{Z}$,
if we let $$\textbf{h}^{(m)}_{n}=(F^{(m)}_{n+m-1}, F^{(m)}_{n+m-2},
F^{(m)}_{n+m-3}, \cdots, F^{(m)}_{n+1}, F^{(m)}_{n} )^{\top},$$
thus, the generating sequence $\{T^{(m)}_{n}\}_{-\infty}^{\infty}$
of the $m$-dimensional harmonic-type golden section $T^{(m)}$
corresponds one-to-one with the following $m\times\infty$-order
matrix $\mathbb{H}^{(m)}$:
\begin{eqnarray*}
\mathbb{H}^{(m)}&=&\left(
\begin{array}{ccccccccccccccccccc}
\hspace{0.35cm}\cdots\hspace{0.4cm} \textbf{h}^{(m)}_{-n}~
\hspace{0.3cm}\cdots\hspace{0.3cm}
\textbf{h}^{(m)}_{-m+1}~\hspace{0.3cm}\cdots~\hspace{0.3cm}
\textbf{h}^{(m)}_{-1}~\hspace{0.3cm}\textbf{h}^{(m)}_{0}~\hspace{0.4cm}\textbf{h}^{(m)}_{1}~\hspace{0.3cm}\cdots~
\hspace{0.3cm}
\textbf{h}^{(m)}_{m-1}~ \hspace{0.3cm}\cdots ~\hspace{0.4cm}\textbf{h}^{(m)}_{n}~\hspace{0.25cm} \cdots~~~\\
\end{array}
\right)\\
&=&\left(
\begin{array}{ccccccccccccccccccc}
\cdots & F^{(m)}_{-n+m-1}&\cdots &F^{(m)}_{0}  &\cdots &F^{(m)}_{m-2}  &F^{(m)}_{m-1} &F^{(m)}_{m}  &\cdots &F^{(m)}_{2m-2} &\cdots&F^{(m)}_{n+m-1}&\cdots \\
\cdots &F^{(m)}_{-n+m-2}&\cdots &F^{(m)}_{-1} &\cdots &F^{(m)}_{m-3}  &F^{(m)}_{m-2} &F^{(m)}_{m-1}&\cdots &F^{(m)}_{2m-3} &\cdots &F^{(m)}_{n+m-2} &\cdots \\
\cdots &\vdots         &\cdots &\vdots        &\vdots&\vdots         &\vdots       &\vdots        &\vdots   &\vdots         &\vdots  &\vdots     &\vdots \\
\cdots &F^{(m)}_{-n+2} &\cdots &F^{(m)}_{-m+3}&\cdots&F^{(m)}_{1}    &F^{(m)}_{2}  &F^{(m)}_{3}   &\cdots &F^{(m)}_{m+1}  &\cdots &F^{(m)}_{n+2} &\cdots \\
\cdots &F^{(m)}_{-n+1} &\cdots &F^{(m)}_{-m+2}&\cdots&F^{(m)}_{0}    &F^{(m)}_{1}  &F^{(m)}_{2}   &\cdots &F^{(m)}_{m}   &\cdots&F^{(m)}_{n+1}&\cdots \\
\cdots &F^{(m)}_{-n}  &\cdots &F^{(m)}_{-m+1} &\cdots&F^{(m)}_{-1}   &F^{(m)}_{0}  &F^{(m)}_{1}   &\cdots &F^{(m)}_{m-1} &\cdots&F^{(m)}_{n}  &\cdots\\
 \end{array}
\right).\\
\end{eqnarray*}

Since $\{T^{(m)}_{n}\}_{-\infty}^{\infty}$ and $\mathbb{H}^{(m)}$
are in a one-to-one correspondence, we can study the generating
sequence $\{T^{(m)}_{n}\}_{-\infty}^{\infty}$ of the $m$-dimensional
harmonic-type golden section $T^{(m)}$ through $\mathbb{H}^{(m)}$.
To this end, we further generalize $\mathbb{H}^{(m)}$ by adding
infinitely many rows both above and below $\mathbb{H}^{(m)}$,
thereby obtaining an $\infty\times\infty$ matrix $\mathcal{H}^{(m)}$
composed of $m$-generalized Fibonacci numbers:
\begin{eqnarray*}
\left(
\begin{array}{ccccccccccccccccccc}
\vdots  &\vdots             &\vdots &\vdots           &\vdots            &\uparrow \bf R            &\vdots         &\vdots         &\vdots  &\vdots     &\vdots \\
\cdots  &F^{(m)}_{0}        &\cdots &F^{(m)}_{n-2}    &F^{(m)}_{n-1}     &{\bf F}^{(m)}_{n} &F^{(m)}_{n+1}  &F^{(m)}_{n+2}  &\cdots  &F^{(m)}_{2n} &\cdots \\
\vdots  &\vdots             &\vdots&\vdots            &\vdots           &\vdots            &\vdots         &\vdots          &\vdots  &\vdots     &\vdots \\
\cdots  &F^{(m)}_{-n+2}     &\cdots&F^{(m)}_{0}       &F^{(m)}_{1}       &{\bf F}^{(m)}_{2} &F^{(m)}_{3}   &F^{(m)}_{4}     &\cdots &F^{(m)}_{n+2} &\cdots \\
\cdots  &F^{(m)}_{-n+1}     &\cdots&F^{(m)}_{-1}      &F^{(m)}_{0}       &{\bf F}^{(m)}_{1} &F^{(m)}_{2}   &F^{(m)}_{3}     &\cdots &F^{(m)}_{n+1} &\cdots \\
\cdots  &{\bf F}^{(m)}_{-n} &\cdots&{\bf F}^{(m)}_{-2}&{\bf F}^{(m)}_{-1}&\bf F^{(m)}_{0}&{\bf F}^{(m)}_{1}&{\bf F}^{(m)}_{2}&\cdots&{\bf F}^{(m)}_{n}&\rightarrow \bf C\\
\cdots  &F^{(m)}_{-n-1}     &\cdots&F^{(m)}_{-3}      &F^{(m)}_{-2}      &{\bf F}^{(m)}_{-1}&F^{(m)}_{0}   &F^{(m)}_{1}     &\cdots &F^{(m)}_{n-1}&\cdots \\
\cdots  &F^{(m)}_{-n-2}     &\cdots&F^{(m)}_{-4}      &F^{(m)}_{-3}     &{\bf F}^{(m)}_{-2} &F^{(m)}_{-1}  &F^{(m)}_{0}     &\cdots &F^{(m)}_{n-2} &\cdots \\
\vdots  &\vdots             &\vdots&\vdots            &\vdots           &\vdots             &\vdots         &\vdots         &\vdots &\vdots     &\vdots \\
\cdots  &F^{(m)}_{-2n}     &\cdots &F^{(m)}_{-n-2}   &F^{(m)}_{-n-1}    &{\bf F}^{(m)}_{-n} &F^{(m)}_{-n+1}&F^{(m)}_{-n+2}  &\cdots &F^{(m)}_{0} &\cdots \\
\vdots  &\vdots            &\vdots &\vdots           &\vdots            &\vdots            &\vdots         &\vdots         &\vdots  &\vdots     &\vdots \\
 \end{array}
\right).\\
\end{eqnarray*}

\vspace{1.5mm}\noindent{\bf Definition 4.1} For all $n \in
\mathbb{Z}$ and $~m\geq 2$, the $\infty\times\infty$-order matrix
$\mathcal{H}^{(m)}$ is called the Fibonacci matrix of the
$m$-dimensional harmonic-type golden section $T^{(m)}$.

We select an element $F^{(m)}_{0}$ in $\mathcal{H}^{(m)}$ as the
(virtual) origin, denoted as $\bf{F}^{(m)}_{0}$. Two perpendicular
dashed lines intersecting vertically at the origin
$\bf{F}^{(m)}_{0}$--one extending upward and the other to the
right--are then established, referred to as the row index axis
($R$-axis) and column index axis ($C$-axis), respectively. The
subscripts $i$ and $j~(i, j\in \mathbb{Z} )$ of the elements
$F^{(m)}_{i}$ and $F^{(m)}_{j}$ on the $R$-axis and $C$-axis are
called their row index and column index, respectively. The row
containing $F^{(m)}_{i}$ is designated as the $i-$th row of
$\mathcal{H}^{(m)}$, with the elements on the $C$-axis specifically
assigned to the 0-th row. Similarly, the column containing
$F^{(m)}_{j}$ is designated as the $j$-th column, with the elements
on the $R$-axis assigned to the 0-th column. The system of
perpendicularly intersecting row and column index axes established
on the matrix $\mathcal{H}^{(m)}$ as described above is referred to
as the (virtual) {\bf row-column index system} of
$\mathcal{H}^{(m)}$. The matrix $\mathcal{H}^{(m)}$ is thereby
divided into four blocks, corresponding to the four quadrants of a
Cartesian coordinate system, termed Block I, II, III, and IV of
$\mathcal{H}^{(m)}$.

The row-column index system of the matrix $\mathcal{H}^{(m)}$
differs from the Cartesian coordinate system. Since the origin must
be established at $F^{(m)}_{0}$, the row index axis and column index
axis are mutually determined:

 $\bullet$ If the row index axis is
selected first, the column index axis is consequently determined.
Such a system is called a row-primary index system.

 $\bullet$
Conversely, if the column index axis is chosen first, the row index
axis follows accordingly. This is referred to as a column-primary
index system.

 Moreover, the establishment of the row-column index
system should be flexibly adapted to the specific problem at hand.
For example, to prove
\begin{eqnarray}
F^{(2)}_{n}F^{(2)}_{n+5}+2F^{(2)}_{n+1}F^{(2)}_{n+4}=3F^{(2)}_{n+2}F^{(2)}_{n+3},
\end{eqnarray}
we consider the following two second-order matrices:
\begin{eqnarray*}
\mathbf{H}_{1}=\mathbf{H}^{(2)}_{n+1}(1, 2)=\left(
\begin{array}{cc}
F^{(2)}_{n+2}     &F^{(2)}_{n+4}    \\
F^{(2)}_{n+1}      &F^{(2)}_{n+3}     \\
 \end{array}
\right),~~~~~~~\mathbf{H}_{2}=\mathbf{H}^{(2)}_{n}(2, 3)=\left(
\begin{array}{cc}
 F^{(2)}_{n+2}       &F^{(2)}_{n+5}      \\
F^{(2)}_{n}      &F^{(2)}_{n+3} \\
 \end{array}
\right).
\end{eqnarray*}

For simplicity, we establish a column-primary index system for both
matrices by setting their bottom rows as the $C$-axis (column index
axis) in $\mathcal{H}^{(2)}$. Under this system:

$\bullet$ For $\mathbf{H}^{(2)}_{n+1}(1, 2)$: Rows: Row 0 and Row 1,
Columns: Column $n+1$ and Column $n+3$.

$\bullet$ For $\mathbf{H}^{(2)}_{n}(2, 3)$: Rows: Row 0 and Row 2,
Columns:
Column $n$ and Column $n+3$.\\
Alternatively, a row-primary index system can be established by
designating the leftmost column as the $R$-axis (row index axis). In
this case:

$\bullet$ For $\mathbf{H}_{1}$: Columns: Column 0 and Column 2,
Rows: Row $n+1$ and Row $n+2$.

$\bullet$ For $\mathbf{H}_{2}$: Columns: Column 0 and Column 3,
Rows: Row $n$
and Row $n+2$.\\
From Corollary 4.1 (see below), we derive:
\begin{eqnarray*}
3F^{(2)}_{n+2}F^{(2)}_{n+3}-F^{(2)}_{n}F^{(2)}_{n+5}-2F^{(2)}_{n+1}F^{(2)}_{n+4}=
2det\mathbf{H}_{1}+det\mathbf{H}_{2}=2(-1)^{n+1}F^{(2)}_{1}F^{(2)}_{2}+(-1)^{n}F^{(2)}_{2}F^{(2)}_{3}=0.
\end{eqnarray*}

Unless otherwise specified, the following discussions will be
conducted under the column-primary index system, where a selected
row of the matrix in question serves as the column ($C$) index axis.

If the row indices and column indices of a submatrix of
$\mathcal{H}^{(m)}$ are consecutive, then this submatrix is called a
{\bf consecutive submatrix} of $\mathcal{H}^{(m)}$. Otherwise, it is
referred to as a {\bf non-consecutive submatrix} of
$\mathcal{H}^{(m)}$. An $m$-th order square submatrix of
$\mathcal{H}^{(m)}$ is called a {\bf full-dimensional submatrix} of
$\mathcal{H}^{(m)}$, denoted by $\bf H^{(m)}$. If a full-dimensional
submatrix $\mathbf{H}^{(m)}$ of $\mathcal{H}^{(m)}$ contains at
least one consecutive $(m-1)$-th order square submatrix, then
$\mathbf{H}^{(m)}$ is said to be $(m-1)$-consecutive. Moreover, the
full-dimensional submatrix $\bf H^{(m)}$ of $\mathcal{H}^{(m)}$ can
be categorized into: {\bf consecutive full-dimensional
submatrices}$^{1}$ \footnote{ It is not difficult to see that:
if $n\in \mathbb{N}$ is restricted, then the consecutive
full-dimensional submatrix $\bf H^{(m)}$ is the transpose of Miles'
matrix $A_{j,k}$ (Miles [43]).} and {\bf non-consecutive
full-dimensional submatrices}. A consecutive full-dimensional
submatrix $\mathbf{H}^{(m)}$ of $\mathcal{H}^{(m)}$ is also called
$m$-consecutive. It can be observed that: If a submatrix of the
Fibonacci matrix is $m$-consecutive, then it must also be
$(m-1)$-consecutive. However, the converse does not necessarily
hold.

Let $k\in~\mathbb{Z}^{+}$ and $1\leq k\leq m-1$. An $(m-k)$-th order
square submatrix of $\mathcal{H}^{(m)}$ is called a {\bf
dimension-reduced $k$-submatrix} of $\mathcal{H}^{(m)}$, denoted by
$\bf H^{(m_{-k})}$. For all $1\leq k\leq m-1$, the collection of all
dimension-reduced $k$-submatrices is referred to as the {\bf
dimension-reduced submatrices} of $\mathcal{H}^{(m)}$, denoted by
$\bf H^{(m_{-})}$. Similarly, there are {\bf dimension-reduced
consecutive submatrices} and {\bf dimension-reduced non-consecutive
submatrices}.

Let $k\in~\mathbb{Z}^{+}$ and $k\geq 1$. An $(m+k)$-th order
(square) submatrix of $\mathcal{H}^{(m)}$ is called a {\bf
dimension-augmented $k$-submatrix} of $\mathcal{H}^{(m)}$, denoted
by $\bf H^{(m_{+k})}$. For all $k\in~\mathbb{Z}^{+}$, the collection
of all dimension-augmented $k$-submatrices is referred to as the
{\bf dimension-augmented submatrices} of $\mathcal{H}^{(m)}$,
denoted by $\bf H^{(m_{+})}$. Similarly, $\bf H^{(m_{+})}$ can be
classified into {\bf dimension-augmented consecutive submatrices}
and {\bf dimension-augmented non-consecutive submatrices}.

\subsection*{4.2\ $(m-1)$-Consecutive Full-dimensional
Submatrices in Fibonacci Matrices}

We investigate the following four types of $(m-1)$-consecutive
full-dimensional submatrices in the Fibonacci matrix
$\mathcal{H}^{(m)}$.

$(i)$~~Let $n\in \mathbb{Z}$, $u, r \in \mathbb{N}$. Consider the
$m$-th order square matrix formed by elements at the intersections
of:

$\bullet$ Rows: consecutive $0$-th, $1$-th, $\cdots$, $(m-2)$-th
rows plus the $u$-th row;

$\bullet$ Columns: consecutive $n$-th, $(n+1)$-th, $\cdots$,
$(n+m-2)$-th columns plus the $(n+r)$-th column.\\
This matrix is denoted by $\mathcal{H}^{(m)}_{n}(u,r)$ and called
the $\urcorner$-shaped full-dimensional submatrix of
$\mathcal{H}^{(m)}$. Special case: When $u=m-1$ and $r=m-1$,
$\textbf{H}^{(m)}_{n}(m-1,m-1)$ becomes the full-dimensional
consecutive submatrix of $\mathcal{H}^{(m)}$, simply denoted as
$\textbf{H}^{(m)}_{n}$.

When $0\leq u\leq m-2$ or $0\leq r \leq m-2$, the matrix
$\textbf{H}^{(m)}_{n}(u,r)$ contains two identical rows or two
identical columns. In such cases, we call such a matrix a {\bf
quasi-full-dimensional submatrix}$^{2}$ \footnote{ For example,
in $\mathcal{H}^{(2)}$ with $u=0$ and $r =0$, the matrix
$\textbf{H}^{(2)}_{n}(0,0)$ consists of four $F^{(2)}_{n}$ elements
located at the intersections of two identical 0-rows and two
identical $n$-columns. Therefore, the second-order
quasi-full-dimensional submatrix of $\mathcal{H}^{(2)}$ composed of
four $F^{(2)}_{n}$ elements is given by
$\textbf{H}^{(2)}_{n}(0,0)=\left(
\begin{array}{cc}
F^{(2)}_{n}     &F^{(2)}_{n}    \\
F^{(2)}_{n}      &F^{(2)}_{n}     \\
 \end{array}
\right).$}of $\mathcal{H}^{(m)}$.

$(ii)$~~Let $n \in \mathbb{Z}$ and $u \in \mathbb{N},~l \in
\mathbb{Z}^{+}$. The $m$-th order square matrix formed by the
elements at the intersections of:

$\bullet$ Rows: The consecutive rows 0, 1, $\dots$, $m-2$ and the
$u$-th row of $\mathcal{H}^{(m)}$;

$\bullet$ Columns: The $(n-l)$-th column and the consecutive columns
$n, n+1, \dots, n+m-2$,\\
arranged in their original order, is denoted as
$\textbf{H}^{(m)}_{n}(u,-l)$ and called an $\ulcorner$-shaped
full-dimensional submatrix of $\mathcal{H}^{(m)}$. Special case:
$\textbf{H}^{(m)}_{n}(m-1,-1)$ is a full-dimensional consecutive
submatrix, the notation can be simplified to
$\textbf{H}^{(m)}_{n-1}$.

$(iii)$~~Let $ n \in \mathbb{Z}$, $  d,~l \in \mathbb{Z}^{+}$. The
$m$-th order square matrix formed by the elements at the
intersections of:

$\bullet$ Rows: The $-d$-th row and the consecutive rows $0, 1,
\dots, m-2$ of $\mathcal{H}^{(m)}$;

$\bullet$ Columns: The $(n-l)$-th column and the consecutive columns
$n, n+1, \dots, n+m-2$,\\
arranged in their original order, is denoted as
$\textbf{H}^{(m)}_{n}(-d, -l)$ and called an $\llcorner$-shaped
full-dimensional submatrix of $\mathcal{H}^{(m)}$. Special case:
$\textbf{H}^{(m)}_{n}(-1,-1)$ is the full-dimensional consecutive
submatrix, and in this case, the notation simplifies to
$\textbf{H}^{(m)}_{n-2}$.

$(iv)$~~Let $ n \in \mathbb{Z}$, $ d \in \mathbb{Z}^{+}$ and $r \in
\mathbb{N}$. The $m$-th order square matrix formed by the elements
at the intersections of:

$\bullet$ Rows: The $-d$-th row and consecutive rows $0, 1, \dots,
m-2$ of $\mathcal{H}^{(m)}$;

$\bullet$ Columns: Consecutive columns $n, n+1, \dots, n+m-2$ and
the $(n-r)$-th column,\\
arranged in their original order, is denoted as
$\textbf{H}^{(m)}_{n}(-d,r)$ and called a $\lrcorner$-shaped
full-dimensional submatrix of $\mathcal{H}^{(m)}$. Special case:
$\textbf{H}^{(m)}_{n}(-1,m-1)$ is the full-dimensional consecutive
submatrix, with simplified notation $\textbf{H}^{(m)}_{n-1}$.

All four of these matrices are $(m-1)$-consecutive full-dimensional
submatrices of the Fibonacci matrix $\mathcal{H}^{(m)}$. They share
a common $(m-1)$-th order consecutive submatrix, is called their
principal submatrix. Formed by the elements at the intersection of
consecutive rows (0-th, 1-st, ..., $(m-2)$-th) and consecutive
columns ($n$-th, $(n+1)$-th, ..., $(n+m-2)$-th) of
$\mathcal{H}^{(m)}$, arranged in their original order:
$$\left(
\begin{array}{ccccccccccccccccc}
F^{(m)}_{n+m-2} &F^{(m)}_{n+m-1}&F^{(m)}_{n+m}&\cdots &F^{(m)}_{n+2m-5}&F^{(m)}_{n+2m-4} \\
F^{(m)}_{n+m-3} &F^{(m)}_{n+m-2}&F^{(m)}_{n+m-1}&\cdots &F^{(m)}_{n+2m-4}&F^{(m)}_{n+2m-3} \\
\vdots          &\vdots          &\vdots        &\vdots         &\vdots        &\vdots\\
F^{(m)}_{n+2}   &F^{(m)}_{n+3} &F^{(m)}_{n+4}  &\cdots &F^{(m)}_{n+m-1}  &F^{(m)}_{n+m}\\
F^{(m)}_{n+1}  &F^{(m)}_{n+2}  &F^{(m)}_{n+3}  &\cdots &F^{(m)}_{n+m-2}  &F^{(m)}_{n+m-1}  \\
F^{(m)}_{n}    &F^{(m)}_{n+1}  &F^{(m)}_{n+2}  &\cdots  &F^{(m)}_{n+m-3}  &F^{(m)}_{n+m-2}  \\
 \end{array}
\right).$$ In fact, this principal submatrix is the
dimension-reduced by 1 consecutive submatrix of $\mathcal{H}^{(m)}$
as defined above, denoted as $\textbf{H}^{(m_{-1})}_{n}$. Here, $u$
and $d$ are called the upper row distance and lower row distance of
the principal submatrix, respectively, while $r$ and $l$ are called
the right column distance and left column distance of the principal
submatrix, respectively.

Below, we will discuss the properties of their determinants and
their applications separately according to the above four cases.
\subsection*{4.2.1 Properties and Applications of $\det\textbf{H}^{(m)}_{n}(u,r)$}
When $0\leq u \leq m-2$ or $0\leq r\leq m-2$,
$\textbf{H}^{(m)}_{n}(u,r)$ is a quasi-full-dimensional submatrix,
and $det\textbf{H}^{(m)}_{n}(u,r)=0$. When $ u= m-1$ and $ r= m-1$,
$\textbf{H}^{(m)}_{n}(u,r)$ is a consecutive full-dimensional
submatrix; when $ u\geq m-1$ and $ r\geq m-1$ but not simultaneously
equal, $\textbf{H}^{(m)}_{n}(u,r)$ is a non-consecutive
full-dimensional
 submatrix. In all cases, we will derive and prove a
unified expression for $det\textbf{H}^{(m)}_{n}(u,r)$ in the
following text, and apply it to the two-dimensional case to prove
some classical Fibonacci sequence identities.

\vspace{1.5mm}\noindent{\bf Lemma 4.1}~ {\it When $n=-1, 0, 1$, the
full-dimensional submatrices} $\textbf{H}^{(m)}_{-1}(m-1, m-1)$,
$\textbf{H}^{(m)}_{0}(m-1, m-1)$, $\textbf{H}^{(m)}_{1}(m-1, m-1)$,
{\it i.e., the consecutive full-dimensional submatrices}
$\textbf{H}^{(m)}_{-1}$, $\textbf{H}^{(m)}_{0}$,
$\textbf{H}^{(m)}_{1}$, {\it have determinants given by:}

$(i)~
det\textbf{H}^{(m)}_{-1}=(-1)^{-1(m-1)}F^{(m)}_{m-1}F^{(m)}_{m-1}$;\hspace{2cm}
$(ii)~det\textbf{H}^{(m)}_{0}=(-1)^{0(m-1)}F^{(m)}_{m-1}F^{(m)}_{m-1}$;

$(iii)~
det\textbf{H}^{(m)}_{1}=(-1)^{1(m-1)}F^{(m)}_{m-1}F^{(m)}_{m-1}$.

\vspace{1.5mm}{\bf Proof ~} $(i)$ Since $F^{(m)}_{-1}=1,
F^{(m)}_{0}=F^{(m)}_{1}=\cdots=F^{(m)}_{m-2}=0, F^{(m)}_{m-1}=1$, it
follows that
\begin{eqnarray*}
det\textbf{H}^{(m)}_{-1} &=&det\left(
\begin{array}{ccccccccccccccccc}
F^{(m)}_{m-2} &\cdots&F^{(m)}_{2m-6}&F^{(m)}_{2m-5} &F^{(m)}_{2m-4}&F^{(m)}_{2m-3} \\
F^{(m)}_{m-3} &\cdots&F^{(m)}_{2m-7}&F^{(m)}_{2m-6} &F^{(m)}_{2m-5}&F^{(m)}_{2m-4} \\
\vdots        &\vdots&\vdots        &\vdots         &\vdots        &\vdots\\
F^{(m)}_{1}   &\cdots&F^{(m)}_{m-3} &F^{(m)}_{m-2} &F^{(m)}_{m-1}  &F^{(m)}_{m}\\
F^{(m)}_{0}   &\cdots&F^{(m)}_{m-4} &F^{(m)}_{m-3} &F^{(m)}_{m-2}  &F^{(m)}_{m-1}  \\
F^{(m)}_{-1}   &\cdots&F^{(m)}_{m-5} &F^{(m)}_{m-4} &F_{m-3}^{(m)}&F^{(m)}_{m-2}  \\
 \end{array}
\right)\\
\end{eqnarray*}
\begin{eqnarray*}
&=&det\left(
\begin{array}{ccccccccccccccccc}
0     &\cdots&F^{(m)}_{2m-6}&F^{(m)}_{2m-5} &F^{(m)}_{2m-4}&F^{(m)}_{2m-3} \\
0     &\cdots&F^{(m)}_{2m-7}&F^{(m)}_{2m-6} &F^{(m)}_{2m-5}&F^{(m)}_{2m-4} \\
\vdots&\vdots&\vdots        &\vdots         &\vdots        &\vdots\\
0     &\cdots&0             &0              &1             &F^{(m)}_{m}\\
0     &\cdots&0             &0              &0             &1  \\
1     &\cdots&0             &0              &0             &0  \\
 \end{array}
\right)\\&=&(-1)^{m+1}=(-1)^{-1(m-1)}F^{(m)}_{m-1}F^{(m)}_{m-1}.
\end{eqnarray*}

$(ii)$ \begin{eqnarray*} det\textbf{H}^{(m)}_{0} &=&det\left(
\begin{array}{ccccccccccccccccc}
F^{(m)}_{m-1} &\cdots&F^{(m)}_{2m-5}&F^{(m)}_{2m-4} &F^{(m)}_{2m-3}&F^{(m)}_{2m-2} \\
F^{(m)}_{m-2} &\cdots&F^{(m)}_{2m-6}&F^{(m)}_{2m-5} &F^{(m)}_{2m-4}&F^{(m)}_{2m-3} \\
\vdots        &\vdots&\vdots        &\vdots         &\vdots        &\vdots\\
F^{(m)}_{2}   &\cdots&F^{(m)}_{m-2} &F^{(m)}_{m-1} &F^{(m)}_{m}  &F^{(m)}_{m+1}\\
F^{(m)}_{1}   &\cdots&F^{(m)}_{m-3} &F^{(m)}_{m-2} &F^{(m)}_{m-1}  &F^{(m)}_{m}  \\
F^{(m)}_{0}   &\cdots&F^{(m)}_{m-4} &F^{(m)}_{m-3} &F_{m-2}^{(m)}&F^{(m)}_{m-1}  \\
 \end{array}
\right)\\
&=&det\left(
\begin{array}{ccccccccccccccccc}
F^{(m)}_{m-1}     &\cdots&F^{(m)}_{2m-5}&F^{(m)}_{2m-4} &F^{(m)}_{2m-3}&F^{(m)}_{2m-2} \\
0     &\cdots&F^{(m)}_{2m-6}&F^{(m)}_{2m-5} &F^{(m)}_{2m-4}&F^{(m)}_{2m-3} \\
\vdots&\vdots&\vdots        &\vdots         &\vdots        &\vdots\\
0     &\cdots&0             &F^{(m)}_{m-1}  &F^{(m)}_{m} &F^{(m)}_{m+1}\\
0     &\cdots&0             &0              &F^{(m)}_{m-1} &F^{(m)}_{m}  \\
0     &\cdots&0             &0              &0             &F^{(m)}_{m-1}  \\
 \end{array}
\right)\\&=&(F^{(m)}_{m-1})^{m}=1=(-1)^{0(m-1)}F^{(m)}_{m-1}F^{(m)}_{m-1}.
\end{eqnarray*}

$(iii)$ Expanding the rightmost column of $det\textbf{H}^{(m)}_{1}$
using the recurrence relation
$F^{(m)}_{n}=\sum^{m}_{i=1}F^{(m)}_{n-i}$, then subtracting
multiples of the left columns (each multiplied by -1) from the
rightmost column, and finally expanding along the rightmost column,
we obtain
\begin{eqnarray*} det\textbf{H}^{(m)}_{1}
&=&det\left(
\begin{array}{ccccccccccccccccc}
F^{(m)}_{m}   &\cdots&F^{(m)}_{2m-4}&F^{(m)}_{2m-3} &F^{(m)}_{2m-2}&F^{(m)}_{2m-1} \\
F^{(m)}_{m-1} &\cdots&F^{(m)}_{2m-5}&F^{(m)}_{2m-4} &F^{(m)}_{2m-3}&F^{(m)}_{2m-2} \\
\vdots        &\vdots&\vdots        &\vdots         &\vdots        &\vdots\\
F^{(m)}_{3}   &\cdots&F^{(m)}_{m-1} &F^{(m)}_{m} &F^{(m)}_{m+1}  &F^{(m)}_{m+2}\\
F^{(m)}_{2}   &\cdots&F^{(m)}_{m-2} &F^{(m)}_{m-1} &F^{(m)}_{m}  &F^{(m)}_{m+1}  \\
F^{(m)}_{1}   &\cdots&F^{(m)}_{m-3} &F^{(m)}_{m-2} &F_{m-1}^{(m)}&F^{(m)}_{m}  \\
 \end{array}
\right)\\
\end{eqnarray*}
\begin{eqnarray*}
&=&det\left(
\begin{array}{ccccccccccccccccc}
F^{(m)}_{m}   &\cdots&F^{(m)}_{2m-4}&F^{(m)}_{2m-3} &F^{(m)}_{2m-2}&\sum^{m}_{i=1}F^{(m)}_{i+m-2} \\
F^{(m)}_{m-1} &\cdots&F^{(m)}_{2m-5}&F^{(m)}_{2m-4} &F^{(m)}_{2m-3}&\sum^{m}_{i=1}F^{(m)}_{i+m-3} \\
\vdots        &\vdots&\vdots        &\vdots         &\vdots        &\vdots\\
F^{(m)}_{3}   &\cdots&F^{(m)}_{m-1} &F^{(m)}_{m} &F^{(m)}_{m+1}  &\sum^{m}_{i=1}F^{(m)}_{i+1}\\
F^{(m)}_{2}   &\cdots&F^{(m)}_{m-2} &F^{(m)}_{m-1} &F^{(m)}_{m}  &\sum^{m}_{i=1}F^{(m)}_{i} \\
F^{(m)}_{1}   &\cdots&F^{(m)}_{m-3} &F^{(m)}_{m-2} &F_{m-1}^{(m)}&\sum^{m}_{i=1}F^{(m)}_{i-1} \\
 \end{array}
\right)\\
&=&det\left(
\begin{array}{ccccccccccccccccc}
F^{(m)}_{m}   &\cdots&F^{(m)}_{2m-4}&F^{(m)}_{2m-3} &F^{(m)}_{2m-2}&F^{(m)}_{m-1} \\
F^{(m)}_{m-1} &\cdots&F^{(m)}_{2m-5}&F^{(m)}_{2m-4} &F^{(m)}_{2m-3}&F^{(m)}_{m-2} \\
\vdots        &\vdots&\vdots        &\vdots         &\vdots        &\vdots\\
F^{(m)}_{3}   &\cdots&F^{(m)}_{m-1} &F^{(m)}_{m} &F^{(m)}_{m+1}  &F^{(m)}_{2}\\
F^{(m)}_{2}   &\cdots&F^{(m)}_{m-2} &F^{(m)}_{m-1} &F^{(m)}_{m}  &F^{(m)}_{1} \\
F^{(m)}_{1}   &\cdots&F^{(m)}_{m-3} &F^{(m)}_{m-2} &F_{m-1}^{(m)}&F^{(m)}_{0} \\
 \end{array}
\right)\\
&=&det\left(
\begin{array}{ccccccccccccccccc}
1     &\cdots&F^{(m)}_{2m-4}&F^{(m)}_{2m-3} &F^{(m)}_{2m-2}&1 \\
1     &\cdots&F^{(m)}_{2m-5}&F^{(m)}_{2m-4} &F^{(m)}_{2m-3}&0 \\
\vdots&\vdots&\vdots        &\vdots         &\vdots        &\vdots\\
0     &\cdots&1             &1              &2             &0\\
0     &\cdots&0             &1              &1             &0  \\
0     &\cdots&0             &0              &1             &0  \\
 \end{array}
\right)\\&=&(-1)^{m+1}=(-1)^{1(m-1)}F^{(m)}_{m-1}F^{(m)}_{m-1}.
\end{eqnarray*}

\vspace{1.5mm}\noindent{\bf Lemma 4.2}~ {\it For all $n \in
\mathbb{Z}$ and $u=m-1, r = m-1$, the following holds:}
$$det\textbf{H}^{(m)}_{n}(m-1,m-1)=det\textbf{H}^{(m)}_{n}=(-1)^{n(m-1)}F^{(m)}_{m-1}F^{(m)}_{m-1}.$$

\vspace{1.5mm}{\bf Proof} Given $u=m-1$ and $r=m-1$, the matrix
$\textbf{H}^{(m)}_{n}(m-1,m-1)$ is a consecutive full-dimensional
submatrix with $F^{(m)}_{n}$ as its bottom-left element. Thus,
$det\textbf{H}^{(m)}_{n}(m-1,m-1)=det\textbf{H}^{(m)}_{n}$. We
proceed to prove $det\textbf{H}^{(m)}_{n}=
(-1)^{n(m-1)}F^{(m)}_{m-1}F^{(m)}_{m-1}$ by mathematical induction,
distinguishing between the cases $n\in \mathbb{Z}^+$ (positive
integers) and $n \in \mathbb{Z}^-$ (negative integers).

When $n \in \mathbb{Z}^{+}$, By Lemma 4.1, the statement holds for
$n =1$. Assume the statement holds for $n=k$. For $n=k+1$, we expand
the rightmost column of $det\textbf{H}^{(m)}_{k+1}$ using the
recurrence relation $F^{(m)}_{n}=\sum^{m}_{i=1}F^{(m)}_{n-i}$, and
then add -1 times the first $m-1$ columns to the rightmost column,
yielding \begin{eqnarray*}
det\textbf{H}^{(m)}_{k+1}&=&det\left(\begin{array}{ccccccccccccccccc}
F^{(m)}_{k+m}   &\cdots&F^{(m)}_{k+2m-4}&F^{(m)}_{k+2m-3} &F^{(m)}_{k+2m-2}&F^{(m)}_{k+2m-1} \\
F^{(m)}_{k+m-1} &\cdots&F^{(m)}_{k+2m-5}&F^{(m)}_{k+2m-4} &F^{(m)}_{k+2m-3}&F^{(m)}_{k+2m-2+} \\
\vdots        &\vdots&\vdots        &\vdots         &\vdots        &\vdots\\
F^{(m)}_{k+3}   &\cdots&F^{(m)}_{k+m-1} &F^{(m)}_{k+m} &F^{(m)}_{k+m+1}  &F^{(m)}_{k+m+2}\\
F^{(m)}_{k+2}   &\cdots&F^{(m)}_{k+m-2} &F^{(m)}_{k+m-1} &F^{(m)}_{k+m}  &F^{(m)}_{k+m+1}  \\
F^{(m)}_{k+1}   &\cdots&F^{(m)}_{k+m-3} &F^{(m)}_{k+m-2} &F_{k+m-1}^{(m)}&F^{(m)}_{k+m}  \\
 \end{array}
\right)\\
\end{eqnarray*}
\begin{eqnarray*}
&=&det\left(
\begin{array}{ccccccccccccccccc}
F^{(m)}_{k+m}   &\cdots&F^{(m)}_{k+2m-4}&F^{(m)}_{k+2m-3} &F^{(m)}_{k+2m-2}&\sum^{m-1}_{i=0}F^{(m)}_{k+m-1+i} \\
F^{(m)}_{k+m-1} &\cdots&F^{(m)}_{k+2m-5}&F^{(m)}_{k+2m-4} &F^{(m)}_{k+2m-3}&\sum^{m-1}_{i=0}F^{(m)}_{k+m-2+i} \\
\vdots        &\vdots&\vdots        &\vdots         &\vdots        &\vdots\\
F^{(m)}_{k+3}   &\cdots&F^{(m)}_{k+m-1} &F^{(m)}_{k+m} &F^{(m)}_{k+m+1}  &\sum^{m-1}_{i=0}F^{(m)}_{k+2+i}\\
F^{(m)}_{k+2}   &\cdots&F^{(m)}_{k+m-2} &F^{(m)}_{k+m-1} &F^{(m)}_{k+m}  &\sum^{m-1}_{i=0}F^{(m)}_{k+1+i} \\
F^{(m)}_{k+1}   &\cdots&F^{(m)}_{k+m-3} &F^{(m)}_{k+m-2} &F_{k+m-1}^{(m)}&\sum^{m-1}_{i=0}F^{(m)}_{k+i} \\
 \end{array}
\right)\\&=&det\left(
\begin{array}{ccccccccccccccccc}
F^{(m)}_{k+m}   &\cdots&F^{(m)}_{k+2m-4}&F^{(m)}_{k+2m-3} &F^{(m)}_{k+2m-2}&F^{(m)}_{k+m-1} \\
F^{(m)}_{k+m-1} &\cdots&F^{(m)}_{k+2m-5}&F^{(m)}_{k+2m-4} &F^{(m)}_{k+2m-3}&F^{(m)}_{k+m-2} \\
\vdots        &\vdots&\vdots        &\vdots         &\vdots        &\vdots\\
F^{(m)}_{k+3}   &\cdots&F^{(m)}_{k+m-1} &F^{(m)}_{k+m} &F^{(m)}_{k+m+1} &F^{(m)}_{k+2}\\
F^{(m)}_{k+2}   &\cdots&F^{(m)}_{k+m-2} &F^{(m)}_{k+m-1} &F^{(m)}_{k+m}  &F^{(m)}_{k+1} \\
F^{(m)}_{k+1}   &\cdots&F^{(m)}_{k+m-3} &F^{(m)}_{k+m-2} &F_{k+m-1}^{(m)}&F^{(m)}_{k} \\
 \end{array}
\right).\\
\end{eqnarray*}
By successively swapping the rightmost column with the preceding
$m-1$ columns until it becomes the leftmost column, and applying the
induction hypothesis, we obtain \begin{eqnarray*}
det\textbf{H}^{(m)}_{k+1}&=&(-1)^{m-1}det\left(
\begin{array}{ccccccccccccccccccc}
F^{(m)}_{k+m-1}&F^{(m)}_{k+m}   &\cdots&F^{(m)}_{k+2m-4}&F^{(m)}_{k+2m-3} &F^{(m)}_{k+2m-2} \\
F^{(m)}_{k+m-2}&F^{(m)}_{k+m-1} &\cdots&F^{(m)}_{k+2m-5}&F^{(m)}_{k+2m-4} &F^{(m)}_{k+2m-3} \\
\vdots         &\vdots          &\vdots        &\vdots  &\vdots        &\vdots\\
F^{(m)}_{k+2}  &F^{(m)}_{k+3}   &\cdots&F^{(m)}_{k+m-1} &F^{(m)}_{k+m} &F^{(m)}_{k+m+1} \\
F^{(m)}_{k+1}  &F^{(m)}_{k+2}   &\cdots&F^{(m)}_{k+m-2} &F^{(m)}_{k+m-1} &F^{(m)}_{k+m}  \\
F^{(m)}_{k}    &F^{(m)}_{k+1}   &\cdots&F^{(m)}_{k+m-3} &F^{(m)}_{k+m-2} &F_{k+m-1}^{(m)} \\
 \end{array}
\right)\\&=&(-1)^{m-1}det\textbf{H}^{(m)}_{k}=(-1)^{m-1}(-1)^{k(m-1)}F^{(m)}_{m-1}F^{(m)}_{m-1}\\
&=&(-1)^{(k+1)(m-1)}F^{(m)}_{m-1}F^{(m)}_{m-1}.
\end{eqnarray*}
Thus, by the principle of mathematical induction, the statement
holds for all $n \in \mathbb{Z}^{+}$.

When $n \in \mathbb{Z}^{-}$, For negative integers, we observe the
recurrence $$F^{(m)}_{-(k+1)}=F^{(m)}_{-k+m-1}-F^{(m)}_{-k+m-2}-
\cdots-F^{(m)}_{-k+1}-F^{(m)}_{-k},$$ which allows us to extend the
proof to $n \in \mathbb{Z}^{-}$. Combining this with Lemma 4.1
(which ensures the statement holds for $n=0$), we conclude that the
result holds for all $n \in \mathbb{Z}$.

By the principle of induction, the lemma is proved.

We observe that if we restrict $n \in \mathbb{N}$ and combine it
with Footnote 1, then the result of Lemma 4.2 coincides with Miles'
result $\det A_{j,k} = (-1)^{(2j+k)(k-1)/2}$ (Miles [43]).

\vspace{1.5mm}\noindent{\bf Lemma 4.3}~ {\it Let $n \in \mathbb{Z}$.
If $0 \leq u \leq m-1$ and $0 \leq r \leq m-1$, then}
$$det\textbf{H}^{(m)}_{n}(u,r)=(-1)^{n(m-1)}F^{(m)}_{u}F^{(m)}_{r}.$$

\vspace{1.5mm}{\bf Proof} When $u=m-1$ and $r=m-1$, This has already
been proven in Lemma 4.2.

When $0\leq u\leq m-2$~ or~$ 0\leq  r\leq m-2$, The matrix
$\textbf{H}^{(m)}_{n}(u, u)$ must have either two identical rows or
two identical columns, making it a quasi-full-dimensional submatrix
of $\mathcal{H}^{(m)}H (m)$. Thus, $det\textbf{H}^{(m)}_{n}(u,r)=0$.
Additionally, since $F^{(m)}_{0}=F^{(m)}_{1}=\cdots=F^{(m)}_{m-2}=0,
F^{(m)}_{m-1}=1$, it follows that for $0\leq u\leq m-2$~ or~$ 0\leq
r\leq m-2$, either $F^{(m)}_{u}=0$ or $F^{(m)}_{r}=0$. Therefore,
$$det\textbf{H}^{(m)}_{n}(u,r)=(-1)^{n(m-1)}F^{(m)}_{u}F^{(m)}_{r}.$$
Which completes the proof.

\vspace{1.5mm}\noindent{\bf Lemma 4.4}~ {\it Let $n$ be any given
integer. Then for all positive integer $r \geq m-1$, we have}
$$det\textbf{H}^{(m)}_{n}( m-1, r)=
(-1)^{n(m-1)}F^{(m)}_{m-1}F^{(m)}_{r}.$$

\vspace{1.5mm}{\bf Proof} We proceed by mathematical induction.

When $r=m-1$, the proposition holds by Lemma 4.2. Assume the
proposition holds for all integers $m-1\leq r\leq k~(k\geq m-1)$.
For $r = k+1$, expanding the rightmost column of
$det\textbf{H}^{(m)}_{n}( m-1, k+1)$ using the recurrence relation
$F^{(m)}_{n}=\sum^{m}_{i=1}F^{(m)}_{n-i}$, and then applying the
properties of determinants along with Lemma 4.3, we obtain

\begin{eqnarray*}
det\textbf{H}^{(m)}_{n}( m-1, k+1) &=&det\left(
\begin{array}{ccccccccccccccccc}
F^{(m)}_{n+m-1}  &F^{(m)}_{n+m}    &\cdots &F^{(m)}_{n+2m-3}   & F^{(m)}_{n+m+k} \\
F^{(m)}_{n+m-2}  &F^{(m)}_{n+m-1}   &\cdots&F^{(m)}_{n+2m-4}   & F^{(m)}_{n+m+k-1} \\
        \vdots   &\vdots            &\cdots &      \vdots      &\vdots\\
F^{(m)}_{n+2}    &F^{(m)}_{n+3}    &\cdots &F^{(m)}_{n+m}     & F^{(m)}_{n+k+3} \\
F^{(m)}_{n+1}    &F^{(m)}_{n+2}    &\cdots &F^{(m)}_{n+m-1}   & F^{(m)}_{n+k+2}\\
F^{(m)}_{n}      &F^{(m)}_{n+1}     &\cdots&F^{(m)}_{n+m-2}   & F^{(m)}_{n+k+1}\\
 \end{array}
\right)\\&=&det\left(
\begin{array}{ccccccccccccccccc}
F^{(m)}_{n+m-1}  &F^{(m)}_{n+m}    &\cdots &F^{(m)}_{n+2m-3}   & \sum^{m-1}_{i=0}F^{(m)}_{n+k+i} \\
F^{(m)}_{n+m-2}  &F^{(m)}_{n+m-1}   &\cdots&F^{(m)}_{n+2m-4}   & \sum^{m-1}_{i=0}F^{(m)}_{n+k-1+i} \\
        \vdots   &\vdots            &\cdots &      \vdots      &\vdots\\
F^{(m)}_{n+2}    &F^{(m)}_{n+3}    &\cdots &F^{(m)}_{n+m}     & \sum^{m-1}_{i=0}F^{(m)}_{n+k-m+3+i} \\
F^{(m)}_{n+1}    &F^{(m)}_{n+2}    &\cdots &F^{(m)}_{n+m-1}   & \sum^{m-1}_{i=0}F^{(m)}_{n+k-m+2+i}\\
F^{(m)}_{n}      &F^{(m)}_{n+1}     &\cdots&F^{(m)}_{n+m-2}   & \sum^{m-1}_{i=0}F^{(m)}_{n+k-m+1+i}\\
 \end{array}
\right)\\&=& \sum^{m-1}_{i=0}det\textbf{H}^{(m)}_{n}(m-1,
k-m+1+i)=\sum^{m-1}_{i=0}\Big((-1)^{n(m-1)}
F^{(m)}_{m-1}F^{(m)}_{k-m+1+i}\Big)\\
&=&(-1)^{n(m-1)}F^{(m)}_{m-1}\sum^{m-1}_{i=0}F^{(m)}_{k-m+1+i}
=(-1)^{n(m-1)}F^{(m)}_{m-1}F^{(m)}_{k+1}.
\end{eqnarray*}

Using a method similar to Lemma 4.4, where induction is applied to
the right column offset $r \geq m-1$, we employ mathematical
induction on the upper row offset $u \geq m-1$ to derive the
following Lemma 4.5.

\vspace{1.5mm}\noindent{\bf Lemma 4.5}~ {\it Let $n \in \mathbb{Z}$.
Then for all positive integers $u \geq m-1$ and $r \geq m-1$, we
have}
$$det\textbf{H}^{(m)}_{n}(u,r)=(-1)^{n(m-1)}F^{(m)}_{u}F^{(m)}_{r}.$$

Combining Lemma 4.3 and Lemma 4.5, we obtain the following Theorem
4.1.

\vspace{1.5mm}\noindent{\bf Theorem 4.1}~ {\it Let $n \in
\mathbb{Z}$. Then for all $u, r \in \mathbb{N}$, we have}
$$det\textbf{H}^{(m)}_{n}(u,r)=(-1)^{n(m-1)}F^{(m)}_{u}F^{(m)}_{r}.$$

Using Theorem 4.1, we can derive various relations for
$m$-generalized Fibonacci sequences. For example, to find the sum of
the first $n$ positive-indexed terms of $\{F^{(3)}_{n}\}$, i.e.,
$\sum^{n}_{i=1}F^{(3)}_{i}$: On one hand, by Theorem 4.1, the
determinant \begin{eqnarray*} &&\left|
\begin{array}{ccccccccccccccccc}
F^{(3)}_{-2}    &F^{(3)}_{-1}        & \sum^{n}_{i=1}F^{(3)}_{i+2} \\
F^{(3)}_{-3}    &F^{(3)}_{-2}        & \sum^{n}_{i=1}F^{(3)}_{i+1}\\
F^{(3)}_{-4}    &F^{(3)}_{-3}        & \sum^{n}_{i=1}F^{(3)}_{i}\\
 \end{array}
\right|=
\sum^{n}_{i=1}det\textbf{H}^{(3)}_{-4}(2, i+4)=\sum^{n}_{i=1}\Big((-1)^{-4(3-1)}F^{(3)}_{2}F^{(3)}_{i+4}\Big)\\
&&=\sum^{n}_{i=1}F^{(3)}_{i+4}=\sum^{n}_{i=1}F^{(3)}_{i}+F^{(3)}_{n+1}+F^{(3)}_{n+2}+F^{(3)}_{n+3}+
F^{(3)}_{n+4}-F^{(3)}_{1}-F^{(3)}_{2}-F^{(3)}_{3}-F^{(3)}_{4}.
\end{eqnarray*}
On the other hand, from the properties of determinants, we have
\begin{eqnarray*}
\left|
\begin{array}{ccccccccccccccccc}
F^{(3)}_{-2}    &F^{(3)}_{-1}        & \sum^{n}_{i=1}F^{(3)}_{i+2} \\
F^{(3)}_{-3}    &F^{(3)}_{-2}        & \sum^{n}_{i=1}F^{(3)}_{i+1}\\
F^{(3)}_{-4}    &F^{(3)}_{-3}        & \sum^{n}_{i=1}F^{(3)}_{i}\\
 \end{array}
\right|=\left|
\begin{array}{ccccccccccccccccc}
-1   &1        & \sum^{n}_{i=1}F^{(3)}_{i+2} \\
0    &-1        & \sum^{n}_{i=1}F^{(3)}_{i+1}\\
2    &0        & \sum^{n}_{i=1}F^{(3)}_{i}\\
 \end{array}
\right|=5\sum^{n}_{i=1}F^{(3)}_{i}+4F^{(3)}_{n+1}+2F^{(3)}_{n+2}-2F^{(3)}_{2}-4F^{(3)}_{1}.
\end{eqnarray*} Therefore, it follows that
$$\sum^{n}_{i=1}F^{(3)}_{i}=\frac{1}{2}\big(F^{(3)}_{n+3}+F^{(3)}_{n+1}-1\big).$$

\vspace{1.5mm}\noindent{\bf Corollary 4.1}~ {\it Let $n \in
\mathbb{Z}$. Then for all $u, r \in \mathbb{N}$, we have}
$$det\textbf{H}^{(2)}_{n}(u,r)=(-1)^{n}F^{(2)}_{u}F^{(2)}_{r}.$$

From Corollary 4.1, we can easily derive some important properties
and formulas for the classical Fibonacci sequence, such as the
following formulas $1^{\circ}--16^{\circ}$ (Wu [56]). Since the
subsequent derivations apply Corollary 4.1 to the classical
Fibonacci sequence, we simplify the notation $F^{(2)}_{n}$ to
$F_{n}$ for brevity.

\begin{eqnarray*}
\begin{array}{ll}
 1^{\circ}~F^{2}_{1}+F^{2}_{2n+1}+F^{2}_{2n+3}=3F_{1}F_{2n+1}F_{2n+3}    &2^{\circ}~ F_{n+m}=F_{n-1}F_{m}+F_{n}F_{m+1} \\
3^{\circ}~F_{m-1}F_{n}-F_{m}F_{n-1}=(-1)^{m}F_{n-m}                      &4^{\circ}~F_{n+1}F_{n-1}-F^{2}_{n}=(-1)^{n}F_{1} \\
5^{\circ}~F_{m+k}F_{n-k}-F_{m}F_{n}=(-1)^{n}F_{m-n-k}F_{k}              &6^{\circ}~\sum^{n}_{r=1}F_{r}=F_{n+2}-1 \\
 7^{\circ}~F_{n+4}F_{n-2}-F_{n+3}F_{n-1}=(-1)^{n-1}5                     &8^{\circ}~ F^{2}_{n+1}-F^{2}_{n-1}=F_{2n} \\
9^{\circ}~F_{n}F_{n+4}-F_{n+1}F_{n+3}=2(-1)^{n-1}                        &10^{\circ}~F_{n}F_{n+4}+F_{n+1}F_{n+3}=2F^{2}_{n+2} \\
 11^{\circ}~ F_{n}F_{n+6}+5F_{n+1}F_{n+5}+F_{n+2}F_{n+4}=7F^{2}_{n+3}   &12^{\circ}~F^{2}_{n}+F^{2}_{n+1}=F_{2n+1}\\
 13^{\circ}~F^{2}_{n-1}+F_{n}F_{n-1}-F^{2}_{n}=(-1)^{n}                  &14^{\circ}~F^{2}_{n+1}-F_{n+4}F_{n-2}=(-1)^{n}4\\
 \end{array}
\end{eqnarray*}
\begin{eqnarray*}
15^{\circ} &&F_{n}F_{n+14}+3F_{n+1}F_{n+13}+F_{n+2}F_{n+12}+F_{n+3}F_{n+11}
+F_{n+4}F_{n+10}+2F_{n+5}F_{n+9}+5F_{n+6}F_{n+8}\\
&&=14F^{2}_{n+7}
\end{eqnarray*}
\begin{eqnarray*} 16^{\circ}\hspace{-0.6cm}&&F_{n}F_{n+22}+3F_{n+1}F_{n+21}+F_{n+2}F_{n+20}+F_{n+3}F_{n+19}+3F_{n+4}F_{n+18}+F_{n+5}F_{n+17}+F_{n+6}F_{n+16}\hspace{0.6cm}\\&&
+4F_{n+7}F_{n+15}+2F_{n+8}F_{n+14}+2F_{n+9}F_{n+13}+F_{n+10}F_{n+12}=20F^{2}_{n+11}\end{eqnarray*}

For formula $1^{\circ}~$, it suffices to note that \begin{eqnarray*}
&&F_{2n+1}F_{2n+3}-F^{2}_{2n+1}+F_{2n+1}F_{2n+3}-F^{2}_{2n+3}+F_{2n+1}F_{2n+3}\\
&=&F_{2n+1}F_{2n+2}-F_{2n+2}F_{2n+3}+F_{2n+1}F_{2n+3}=-F^{2}_{2n+2}+F_{2n+1}F_{2n+3}=-\left|
\begin{array}{cc}
F_{2n+2}     &F_{2n+3}   \\
F_{2n+1}      &F_{2n+2}     \\
  \end{array}
  \right|.
 \end{eqnarray*}

For formula $2^{\circ}~$, one only needs to observe
\begin{eqnarray*} (-1)^{-n}F_{n+m}F_{1}&=&\left|
\begin{array}{cc}
F_{-(n-1)}     &F_{m+1}    \\
F_{-n}      &F_{m}     \\
 \end{array}
\right|=F_{-(n-1)}F_{m}-F_{-n}F_{m+1}\\
&=&(-1)^{n}F_{(n-1)}F_{m}+(-1)^{n+2}F_{n}F_{m+1}.\\
\end{eqnarray*}

Formulas $3^{\circ},~4^{\circ},~5^{\circ},~14^{\circ}~$
  are direct results of Corollary 4.1.

For formula $6^{\circ}~$, it can be directly derived using Corollary
4.1:
 \begin{eqnarray*}
1-F_{n+2}=\left|
\begin{array}{cc}
F_{2}     &\sum^{n+2}_{i=3}F_{i}    \\
F_{1}      &\sum^{n+1}_{i=2}F_{i}     \\
 \end{array}
\right|=\sum^{n}_{r=1}\left|
\begin{array}{cc}
 F_{2}       &F_{2+r}      \\
F_{1}      &F_{1+r} \\
 \end{array}
\right|=(-1)^{1}F_{1}\sum^{n}_{r=1}F_{r}.
\end{eqnarray*}

For formula $7^{\circ}~$, simply add $F^{2}_{n+1}-F^{2}_{n+1}$
  to the left-hand side and apply Corollary 4.1.
Formula $9^{\circ}$ can be obtained in the same manner.

For formula $13^{\circ}$, note that the left-hand side equals
$F_{n-1}F_{n+1}-F^{2}_{n}$.

Formula $12^{\circ}$ can be derived as follows, and formula
$8^{\circ}$ is similar.
\begin{eqnarray*}
(-1)^{-n-1}F_{1}F_{2n+1}=\left|
\begin{array}{cc}
F_{-n}     &F_{n+1}    \\
F_{-n-1}      &F_{n}     \\
 \end{array}
\right|=F_{-n}F_{n}-F_{-n-1}F_{n+1}=(-1)^{-n-1}(F^{2}_{n}+F^{2}_{n+1}).
\end{eqnarray*}

Formula $10^{\circ}$
  (as well as $9^{\circ}$ ) is one of several formulas proved by Bertucci [3] in 1991 concerning
 the relationship between Fibonacci numbers and Pythagorean triples. It can be
 easily proven using Corollary 4.1 as follows:
 \begin{eqnarray*}
2F^{2}_{n+2}-F_{n}F_{n+4}-F_{n+1}F_{n+3}&=&F^{2}_{n+2}-F_{n}F_{n+4}+F^{2}_{n+2}-F_{n+1}F_{n+3}\\
&=&\left|
\begin{array}{cc}
F_{n+2}     &F_{n+4}    \\
F_{n}      &F_{n+2}     \\
 \end{array}
\right|+\left|
\begin{array}{cc}
F_{n+2}     &F_{n+3}    \\
F_{n+1}      &F_{n+2}     \\
 \end{array}
\right|\\
&=&(-1)^{n}F^{2}_{2}+(-1)^{n+1}F^{2}_{1}=0.
\end{eqnarray*}

Liu [37], considering the characteristics of formulas $9^{\circ}$
and $10^{\circ}$, proposed and proved (18), $11^{\circ}~$,
$15^{\circ}~$, and$16^{\circ}~$, among other relations involving
consecutive $k~ (5\leq k\leq 24)$ Fibonacci numbers. These formulas
can all be easily verified using Corollary 4.1 in the manner
described above. Liu [37] raised the question of whether a general
formula exists for arbitrary $k\geq 5$. They believed that such a
universal formula does exist, but it has not yet been discovered.

First, we note that such a general formula is not unique. To this
end, we present and prove the following relation involving the
squares of Fibonacci numbers: \begin{eqnarray} \left\{
\begin{array}{l}
  F^{2}_{2}=F^{2}_{1};  \\
F^{2}_{3}=3F^{2}_{2}+F^{2}_{1}, \\
     \end{array}
     \right.
      F^{2}_{n}=\left\{
                 \begin{array}{ll}
                 2F^{2}_{2k-1}+\sum^{2k-2}_{i=2}F^{2}_{i}; & n=2k ~(k\geq 2); \\
                 2F^{2}_{2k}+\sum^{2k-1}_{i=1}F^{2}_{i}+F^{2}_{2}, & n=2k+1.\\
                 \end{array}
                 \right.
\end{eqnarray}

\vspace{1.5mm}{\bf Proof} by Mathematical Induction.

For $k=2$, $F^{2}_{4}=2F^{2}_{3}+F^{2}_{2},
F^{2}_{5}=2F^{2}_{4}+F^{2}_{3}+2F^{2}_{2}+F^{2}_{1},$ the
proposition holds.

Assume the proposition holds for $k$. Now consider the case for
$k+1$:
$$F^{2}_{2(k+1)}=(F_{2k+1}+F_{2k})^{2}=F^{2}_{2k+1}+F^{2}_{2k}+2F_{2k+1}F_{2k}$$
$$F^{2}_{2k}=(F_{2k+1}-F_{2k-1})^{2}=F^{2}_{2k+1}+F^{2}_{2k-1}-2F_{2k+1}F_{2k-1}$$
Adding the two equations yields
\begin{eqnarray*}F^{2}_{2(k+1)}+F^{2}_{2k}&=&2F^{2}_{2k+1}+F^{2}_{2k}+F^{2}_{2k-1}+2F_{2k+1}(F_{2k}-F_{2k-1})\\
&=&2F^{2}_{2k+1}+F^{2}_{2k}+F^{2}_{2k-1}+2F_{2k+1}(F_{2k}-F_{2k-1})\\
&=&2F^{2}_{2k+1}+F^{2}_{2k}+F^{2}_{2k-1}+2F^{2}_{2k}-2F^{2}_{2k-1}\\
&=&2F^{2}_{2k+1}+2F^{2}_{2k}+F^{2}_{2k-1}+F^{2}_{2k}-2F^{2}_{2k-1}\\
&=&2F^{2}_{2k+1}+2F^{2}_{2k}+F^{2}_{2k-1}+(2F^{2}_{2k-1}+\sum^{2k-2}_{i=2}F^{2}_{i})-2F^{2}_{2k-1}\end{eqnarray*}
Using the inductive hypothesis and rearranging terms, we obtain
\begin{eqnarray*}
F^{2}_{2(k+1)}=2F^{2}_{2k+1}+\sum^{2k}_{i=2}F^{2}_{i}.
\end{eqnarray*}
Similarly, \begin{eqnarray*}
F^{2}_{2(k+1)+1}=2F^{2}_{2(k+1)}+\sum^{2(k+1)-1}_{i=1}F^{2}_{i}+F^{2}_{2}.
\end{eqnarray*}
 This proves equation (19).

Since $n+3$ and $n+1$ have the same parity, from (19), we derive
\begin{eqnarray*}
F^{2}_{n+3}-F^{2}_{n+1}=2F^{2}_{n+2}+F^{2}_{n+1}+F^{2}_{n}-2F^{2}_{n}.
\end{eqnarray*}
This leads to the well-known relation for four consecutive Fibonacci
squares (valid for any integer $n$):
\begin{eqnarray}F^{2}_{n+3}-2F^{2}_{n+2}-2F^{2}_{n+1}+F^{2}_{n}=0~~ \hbox{or}~~
F^{2}_{n+3}+F^{2}_{n}=2F^{2}_{n+2}+2F^{2}_{n+1}.\end{eqnarray}
Further, we can derive relations for six, eight, etc., consecutive
Fibonacci squares:
\begin{eqnarray}F^{2}_{n+5}+F^{2}_{n}=2F^{2}_{n+4}+F^{2}_{n+3}+F^{2}_{n+2}+2F^{2}_{n+1},\end{eqnarray}
\begin{eqnarray}
F^{2}_{n+7}+F^{2}_{n}=2F^{2}_{n+6}+F^{2}_{n+5}+F^{2}_{n+4}+F^{2}_{n+3}+F^{2}_{n+2}+2F^{2}_{n+1}.
\end{eqnarray}

Next, we demonstrate that $15^{\circ}$ and $16^{\circ}$  are not
unique. Using Corollary 4.1, formula $15^{\circ}$ can be rewritten
as:\begin{eqnarray*}
F^{2}_{7}+F^{2}_{5}+F^{2}_{3}+5F^{2}_{1}=3F^{2}_{6}+F^{2}_{4}+2F^{2}_{2}.
\end{eqnarray*}
It is easy to verify that this equality holds, confirming the
validity of $15^{\circ}$. Now, if we add $F^{2}_{3}$ to both sides
and apply formula (20), letting
$F^{2}_{3}=2F^{2}_{2}+2F^{2}_{1}-F^{2}_{0},$ on the right-hand side,
and then use $F^{2}_{2}=F^{2}_{1}$ to eliminate $2F^{2}_{1}$ and
$2F^{2}_{2}$, we obtain:
\begin{eqnarray*}
F^{2}_{7}+F^{2}_{5}+2F^{2}_{3}+F^{2}_{1}=3F^{2}_{6}+F^{2}_{4}+2F^{2}_{2}.
\end{eqnarray*}
Thus, $15^{\circ}$ can be refined into $15^{\circ}_{+}$, where the
coefficients on the right-hand side are smaller:
\begin{eqnarray*}
15^{\circ}_{+}&&
F_{n}F_{n+14}+3F_{n+1}F_{n+13}+F_{n+2}F_{n+12}+F_{n+3}F_{n+11}
+2F_{n+4}F_{n+10}+2F_{n+5}F_{n+9}+F_{n+6}F_{n+8}\\
&&=11F^{2}_{n+7}.
\end{eqnarray*}
 Similarly, $16^{\circ}$ can be transformed into:
\begin{eqnarray*}
F^{2}_{11}+F^{2}_{9}+3F^{2}_{7}+F^{2}_{5}+2F^{2}_{3}+F^{2}_{1}=3F^{2}_{10}+F^{2}_{8}+F^{2}_{6}+4F^{2}_{4}+2F^{2}_{2}.
\end{eqnarray*}
By adding $F^{2}_{6}$ and $F^{2}_{9}$ to both sides and applying
formula (20), letting $F^{2}_{6}=2F^{2}_{5}+2F^{2}_{4}-F^{2}_{3},$
on the left and $F^{2}_{9}=2F^{2}_{8}+2F^{2}_{7}-F^{2}_{6},$ on the
right, we simplify to:
\begin{eqnarray*}
F^{2}_{11}+2F^{2}_{9}+F^{2}_{7}+3F^{2}_{5}+F^{2}_{3}+F^{2}_{1}=3F^{2}_{10}+3F^{2}_{8}+F^{2}_{6}+2F^{2}_{4}+2F^{2}_{2}.
\end{eqnarray*} Thus,
$16^{\circ}$ can also be expressed as $16^{\circ}_{+}$:
\begin{eqnarray*}16^{\circ}_{+} \hspace{-0.6cm}&&F_{n}F_{n+22}+3F_{n+1}F_{n+21}+2F_{n+2}F_{n+20}+
3F_{n+3}F_{n+19}+F_{n+4}F_{n+18}+F_{n+5}F_{n+17}+3F_{n+6}F_{n+16}\hspace{0.6cm}\\&&
+2F_{n+7}F_{n+15}+F_{n+8}F_{n+14}+2F_{n+9}F_{n+13}+F_{n+10}F_{n+12}=20F^{2}_{n+11}.\end{eqnarray*}

 When the number
of consecutive Fibonacci numbers $k$ (odd) is large, the
manipulation space expands further. In addition to formula (20),
formulas (21), (22), and even (19) can be used, leading to countless
variations. Similarly, when $k$ is even (as in formula (18)), the
form of the relation is also non-unique. Therefore, we conclude that
such relations can be extended indefinitely. However, since these
relations are not unique, deriving a unified expression for general
$k\geq 5$ requires first establishing a rule and then seeking a
general formula under that rule.
\subsection*{4.2.2 \ Properties of
$det\textbf{H}^{(m)}_{n}(u,-l)$}
For $u \in \mathbb{N}$ and $l\in \mathbb{Z}^{+}$, when $0 \leq u
\leq m-2$, $\textbf{H}^{(m)}_{n}(u,-l)$ has two identical rows.
Similar to Section 4.2.1, we also refer to it as a
quasi-full-dimensional submatrix. In this case,
$det\textbf{H}^{(m)}_{n}(u,-l)=0$. Below, we consider the general
case.

\vspace{1.5mm}\noindent\textbf{Lemma 4.6} ~ {\it For all $n \in
\mathbb{Z}$ and $u \in \mathbb{N}$, we have:}
$$(i)~det\textbf{H}^{(m)}_{n}(u,-1)=(-1)^{(n-1)(m-1)}F^{(m)}_{u}F^{(m)}_{-1},~ (ii)~
det\textbf{H}^{(m)}_{n}(u,-2)=(-1)^{(n-1)(m-1)}F^{(m)}_{u}F^{(m)}_{-2}.$$

\vspace{1.5mm}\textbf{Proof}~ $(i)$~ If we shift the original row
index axis $R$ one column to the left, then in the new principal
column index system, the original matrix becomes
$\textbf{H}^{(m)}_{n-1}(u,m-1)$. By Theorem 4.1, we have:
$$det\textbf{H}^{(m)}_{n-1}(u,m-1)=(-1)^{(n-1)(m-1)}F^{(m)}_{u}F^{(m)}_{m-1}=
(-1)^{(n-1)(m-1)}F^{(m)}_{u}F^{(m)}_{-1}.$$ Therefore,
$$det\textbf{H}^{(m)}_{n}(u,-1)=(-1)^{(n-1)(m-1)}F^{(m)}_{u}F^{(m)}_{-1}.$$

$(ii)$~ Expand the leftmost column of $\textbf{H}^{(m)}_{n}(u,-2)$
using$F^{(m)}_{n}=F^{(m)}_{n+m}-\sum^{m-1}_{i=1}F^{(m)}_{n+i}$, then
add the $n$-th, $(n+1)$-th, $\cdots$, $(n+m-3)$-th columns to the
leftmost column. Finally, multiply the $(n+m-2)$-th column by -1 and
add it to the leftmost column, yielding \begin{eqnarray*}
det\textbf{H}^{(m)}_{n}(u,-2) &=&det\left(
\begin{array}{ccccccccccccccccc}
F^{(m)}_{n-2+u}&F^{(m)}_{n+u}  &F^{(m)}_{n+1+u}   &\cdots&F^{(m)}_{n+m-2+u}    \\
        \vdots   &\vdots            &\vdots &      \vdots      &\vdots\\
 F^{(m)}_{n+1}&F^{(m)}_{n+3}    &F^{(m)}_{n+4}    &\cdots &F^{(m)}_{n+m+1}      \\
 F^{(m)}_{n}&F^{(m)}_{n+2}    &F^{(m)}_{n+3}    &\cdots &F^{(m)}_{n+m}   \\
 F^{(m)}_{n-1}&F^{(m)}_{n+1}    &F^{(m)}_{n+2}     &\cdots&F^{(m)}_{n+m-1}   \\
 F^{(m)}_{n-2} &F^{(m)}_{n}  &F^{(m)}_{n+1}    &\cdots &F^{(m)}_{n+m-2}    \\
 \end{array}
\right)\\
\end{eqnarray*}
\begin{eqnarray*}
&=&det\left(
\begin{array}{ccccccccccccccccc}
 F^{(m)}_{n-2+m+u}-\sum^{m-1}_{i=1}F^{(m)}_{n-2+u+i}&F^{(m)}_{n+u}  &F^{(m)}_{n+1+u}   &\cdots&F^{(m)}_{n+m-2+u}    \\
        \vdots   &\vdots            &\vdots &      \vdots      &\vdots\\
 F^{(m)}_{n+2+m}-\sum^{m-1}_{i=1}F^{(m)}_{n+1+i}&F^{(m)}_{n+3}    &F^{(m)}_{n+4}    &\cdots &F^{(m)}_{n+m+1}      \\
 F^{(m)}_{n+1+m}-\sum^{m-1}_{i=1}F^{(m)}_{n+i}&F^{(m)}_{n+2}    &F^{(m)}_{n+3}    &\cdots &F^{(m)}_{n+m}   \\
 F^{(m)}_{n-1+m}-\sum^{m-1}_{i=1}F^{(m)}_{n-1+i}&F^{(m)}_{n+1}    &F^{(m)}_{n+2}     &\cdots&F^{(m)}_{n+m-1}   \\
 F^{(m)}_{n-2+m}-\sum^{m-1}_{i=1}F^{(m)}_{n-2+i}&F^{(m)}_{n}  &F^{(m)}_{n+1}    &\cdots &F^{(m)}_{n+m-2}    \\
 \end{array}
\right)\\
\end{eqnarray*}
\begin{eqnarray*}
&=&-det\left(
\begin{array}{ccccccccccccccccc}
 F^{(m)}_{n-1+u}&F^{(m)}_{n+u}  &F^{(m)}_{n+1+u}   &\cdots&F^{(m)}_{n+m-2+u}    \\
        \vdots   &\vdots            &\vdots &      \vdots      &\vdots\\
 F^{(m)}_{n+2}&F^{(m)}_{n+3}    &F^{(m)}_{n+4}    &\cdots &F^{(m)}_{n+m+1}      \\
 F^{(m)}_{n+1}&F^{(m)}_{n+2}    &F^{(m)}_{n+3}    &\cdots &F^{(m)}_{n+m}   \\
 F^{(m)}_{n}&F^{(m)}_{n+1}    &F^{(m)}_{n+2}     &\cdots&F^{(m)}_{n+m-1}   \\
F^{(m)}_{n-1}&F^{(m)}_{n}  &F^{(m)}_{n+1}    &\cdots &F^{(m)}_{n+m-2}    \\
 \end{array}
\right)\\&=&-det\textbf{H}^{(m)}_{n}(u,-1)=(-1)^{(n-1)(m-1)}F^{(m)}_{u}(-F^{(m)}_{-1})=
(-1)^{(n-1)(m-1)}F^{(m)}_{u}F^{(m)}_{-2}.\\
\end{eqnarray*}

\vspace{1.5mm}\noindent{\bf Lemma 4.7}~ {\it For all $ n \in
\mathbb{Z}$, $u \in \mathbb{N}$, and $l \in \mathbb{Z}^{+}$, if $3
\leq l \leq m$, then}
$$det\textbf{H}^{(m)}_{n}(u,-l)=0=(-1)^{(n-1)(m-1)}F^{(m)}_{u}F^{(m)}_{-l}.$$

\vspace{1.5mm}{\bf Proof}~ Expand the leftmost column of
$\textbf{H}^{(m)}_{n}(u,-l)$ using
$F^{(m)}_{n}=2F^{(m)}_{n+m}-F^{(m)}_{n+m+1}$. Since $3\leq l\leq m$,
the leftmost column must be equal to twice one of the right columns
minus another column. Therefore, we have \begin{eqnarray*}
det\textbf{H}^{(m)}_{n}(u,-l) &=&det\left(
\begin{array}{ccccccccccccccccc}
F^{(m)}_{n-l+u}&F^{(m)}_{n+u}  &F^{(m)}_{n+1+u}   &\cdots&F^{(m)}_{n+m-2+u}    \\
        \vdots   &\vdots            &\vdots &      \vdots      &\vdots\\
 F^{(m)}_{n-l+3}&F^{(m)}_{n+3}    &F^{(m)}_{n+4}    &\cdots &F^{(m)}_{n+m+1}      \\
 F^{(m)}_{n-l+2}&F^{(m)}_{n+2}    &F^{(m)}_{n+3}    &\cdots &F^{(m)}_{n+m}   \\
 F^{(m)}_{n-l+1}&F^{(m)}_{n+1}    &F^{(m)}_{n+2}     &\cdots&F^{(m)}_{n+m-1}   \\
 F^{(m)}_{n-l}&F^{(m)}_{n}  &F^{(m)}_{n+1}    &\cdots &F^{(m)}_{n+m-2}    \\
 \end{array}
\right)\\&=&det\left(
\begin{array}{ccccccccccccccccc}
2F^{(m)}_{n-l+u+m}-F^{(m)}_{n-l+u+m+1}&F^{(m)}_{n+u}  &F^{(m)}_{n+1+u}   &\cdots&F^{(m)}_{n+m-2+u}    \\
        \vdots   &\vdots            &\vdots &      \vdots      &\vdots\\
 2F^{(m)}_{n-l+3+m}-F^{(m)}_{n-l+3+m+1}&F^{(m)}_{n+3}    &F^{(m)}_{n+4}    &\cdots &F^{(m)}_{n+m+1}      \\
 2F^{(m)}_{n-l+2+m}-F^{(m)}_{n-l+2+m+1}&F^{(m)}_{n+2}    &F^{(m)}_{n+3}    &\cdots &F^{(m)}_{n+m}   \\
 2F^{(m)}_{n-l+1+m}-F^{(m)}_{n-l+1+m+1}&F^{(m)}_{n+1}    &F^{(m)}_{n+2}     &\cdots&F^{(m)}_{n+m-1}   \\
 2F^{(m)}_{n-l+m}-F^{(m)}_{n-l+m+1}&F^{(m)}_{n}  &F^{(m)}_{n+1}    &\cdots &F^{(m)}_{n+m-2}    \\
 \end{array}
\right)\\&=&det\left(
\begin{array}{ccccccccccccccccc}
 0&F^{(m)}_{n+u}  &F^{(m)}_{n+m}   &\cdots&F^{(m)}_{n+2m-3}    \\
        \vdots   &\vdots            &\vdots &      \vdots      &\vdots\\
 0&F^{(m)}_{n+3}    &F^{(m)}_{n+4}    &\cdots &F^{(m)}_{n+m+1}      \\
 0&F^{(m)}_{n+2}    &F^{(m)}_{n+3}    &\cdots &F^{(m)}_{n+m}   \\
 0&F^{(m)}_{n+1}    &F^{(m)}_{n+2}     &\cdots&F^{(m)}_{n+m-1}   \\
 0&F^{(m)}_{n}  &F^{(m)}_{n+1}    &\cdots &F^{(m)}_{n+m-2}    \\
 \end{array}
\right)=0.
\end{eqnarray*}
Moreover, for $m\geq 3$,
$F^{(m)}_{-3}=\cdots=F^{(m)}_{-m+1}=F^{(m)}_{-m}=0$, so when $3\leq
l\leq m$,
$$det\textbf{H}^{(m)}_{n}(u,-l)=0=(-1)^{(n-1)(m-1)}F^{(m)}_{u}F^{(m)}_{-l}.$$

\vspace{1.5mm}\noindent{\bf Theorem 4.2}~ {\it Let $n \in
\mathbb{Z}$ and $u \in \mathbb{N}$. Then for all $l \in
\mathbb{Z}^{+}$, we have}
$$det\textbf{H}^{(m)}_{n}(u,-l)=(-1)^{(n-1)(m-1)}F^{(m)}_{u}F^{(m)}_{-l}.$$

\vspace{1.5mm}{\bf Proof}~ We prove this by mathematical induction.

For $l = 1, 2, \dots, m$, the statement holds by Lemma 4.6 and Lemma
4.7. Now assume the statement holds for $l = k \geq m$. For $l =
k+1$, we have
\begin{eqnarray*}
det\textbf{H}^{(m)}_{n}\big(u,-(k+1)\big) &=&det\left(
\begin{array}{ccccccccccccccccc}
F^{(m)}_{n-(k+1)+u}&F^{(m)}_{n+u}  &F^{(m)}_{n+1+u}   &\cdots&F^{(m)}_{n+m-2+u}    \\
F^{(m)}_{n-(k+1)+m-2}&F^{(m)}_{n+m-2}    &F^{(m)}_{n+m-1}    &\cdots &F^{(m)}_{n+2m-4}      \\
        \vdots   &\vdots            &\vdots &      \vdots      &\vdots\\
 F^{(m)}_{n-(k+1)+2}&F^{(m)}_{n+2}    &F^{(m)}_{n+3}    &\cdots &F^{(m)}_{n+m}   \\
 F^{(m)}_{n-(k+1)+1}&F^{(m)}_{n+1}    &F^{(m)}_{n+2}     &\cdots&F^{(m)}_{n+m-1}   \\
 F^{(m)}_{n-(k+1)}&F^{(m)}_{n}  &F^{(m)}_{n+1}    &\cdots &F^{(m)}_{n+m-2}    \\
 \end{array}
\right)\\
\end{eqnarray*}
\begin{eqnarray*}
&=&det\left(
\begin{array}{ccccccccccccccccc}
 F^{(m)}_{n-k-1+u+m}-\sum^{m-1}_{i=1}F^{(m)}_{n-k-1+u+i}&F^{(m)}_{n+u}  &F^{(m)}_{n+1+u}   &\cdots&F^{(m)}_{n+m-2+u}    \\
 F^{(m)}_{n-k-3+2m}-\sum^{m-1}_{i=1}F^{(m)}_{n-k-3+m+i}&F^{(m)}_{n+m-2}    &F^{(m)}_{n+m-1}    &\cdots &F^{(m)}_{n+2m-4}\\
        \vdots   &\vdots            &\vdots &      \vdots      &\vdots\\
 F^{(m)}_{n-k+1+m}-\sum^{m-1}_{i=1}F^{(m)}_{n-k+1+i}&F^{(m)}_{n+2}    &F^{(m)}_{n+3}    &\cdots &F^{(m)}_{n+m}   \\
 F^{(m)}_{n-k+m}-\sum^{m-1}_{i=1}F^{(m)}_{n-k+i}&F^{(m)}_{n+1}    &F^{(m)}_{n+2}     &\cdots&F^{(m)}_{n+m-1}   \\
 F^{(m)}_{n-k-1+m}-\sum^{m-1}_{i=1}F^{(m)}_{n-k-1+i}&F^{(m)}_{n}  &F^{(m)}_{n+1}    &\cdots &F^{(m)}_{n+m-2}    \\
 \end{array}
\right)\\
\end{eqnarray*}
\begin{eqnarray*}
&=&det\left(
\begin{array}{ccccccccccccccccc}
 F^{(m)}_{n-k-1+u+m}&F^{(m)}_{n+u}  &F^{(m)}_{n+1+u}   &\cdots&F^{(m)}_{n+m-2+u}    \\
 F^{(m)}_{n-k-3+2m}&F^{(m)}_{n+m-2}    &F^{(m)}_{n+m-1}    &\cdots &F^{(m)}_{n+2m-4}\\
        \vdots   &\vdots            &\vdots &      \vdots      &\vdots\\
 F^{(m)}_{n-k+1+m}&F^{(m)}_{n+2}    &F^{(m)}_{n+3}    &\cdots &F^{(m)}_{n+m}   \\
 F^{(m)}_{n-k+m}&F^{(m)}_{n+1}    &F^{(m)}_{n+2}     &\cdots&F^{(m)}_{n+m-1}   \\
 F^{(m)}_{n-k-1+m}&F^{(m)}_{n}  &F^{(m)}_{n+1}    &\cdots &F^{(m)}_{n+m-2}    \\
 \end{array}
\right)\\&-&det\left(
\begin{array}{ccccccccccccccccc}
 \sum^{m-1}_{i=1}F^{(m)}_{n-k-1+u+i}&F^{(m)}_{n+u}  &F^{(m)}_{n+1+u}   &\cdots&F^{(m)}_{n+m-2+u}    \\
 \sum^{m-1}_{i=1}F^{(m)}_{n-k-3+m+i}&F^{(m)}_{n+m-2}    &F^{(m)}_{n+m-1}    &\cdots &F^{(m)}_{n+2m-4}\\
        \vdots   &\vdots            &\vdots &      \vdots      &\vdots\\
 \sum^{m-1}_{i=1}F^{(m)}_{n-k+1+i}&F^{(m)}_{n+2}    &F^{(m)}_{n+3}    &\cdots &F^{(m)}_{n+m}   \\
 \sum^{m-1}_{i=1}F^{(m)}_{n-k+i}&F^{(m)}_{n+1}    &F^{(m)}_{n+2}     &\cdots&F^{(m)}_{n+m-1}   \\
 \sum^{m-1}_{i=1}F^{(m)}_{n-k-1+i}&F^{(m)}_{n}  &F^{(m)}_{n+1}    &\cdots &F^{(m)}_{n+m-2}    \\
 \end{array}
\right)\\&=&det\textbf{H}^{(m)}_{n}\big(u,-(k+1-m)\big)-
\sum^{m-1}_{i=1}det\textbf{H}^{(m)}_{n}(u,-(k+1-i))\\
&=&(-1)^{n(m-1)}F^{(m)}_{u}F^{(m)}_{-(k+1-m)}-\sum^{m-1}_{i=0}\Big((-1)^{n(m-1)}F^{(m)}_{u}F^{(m)}_{-(k+1-i)}\Big)\\
&=&(-1)^{n(m-1)}F^{(m)}_{u}\Big(F^{(m)}_{-(k+1-m)}-\sum^{m-1}_{i=1}F^{(m)}_{-(k+1-i)}\Big)
=(-1)^{n(m-1)}F^{(m)}_{u}F^{(m)}_{-(k+1)}.\\
\end{eqnarray*}

\subsection*{4.2.3 \ Properties of
$det\textbf{H}^{(m)}_{n}(-d,-l)$}

\noindent {\bf Lemma 4.8}~ {\it For all $n \in \mathbb{Z}$, we
have:}
$$(i)~det\textbf{H}^{(m)}_{n}(-1,-1)=(-1)^{n(m-1)}F^{(m)}_{-1}F^{(m)}_{-1},~~(ii)~~
det\textbf{H}^{(m)}_{n}(-2,-1)=(-1)^{n(m-1)}F^{(m)}_{-2}F^{(m)}_{-1}
.$$

\vspace{1.5mm}{\bf Proof}~ $(i)$~ By shifting the original row index
axis $R$ downward by one row, the matrix transforms into
$\textbf{H}^{(m)}_{n-1}(m-1,-1)$ in the new principal column index
system. Given that $F^{(m)}_{m-1}=F^{(m)}_{-1}=1$, Theorem 4.2
yields:
$$det\textbf{H}^{(m)}_{n-1}(m-1,-1)=(-1)^{(n-2)(m-1)}F^{(m)}_{m-1}F^{(m)}_{-1}=
(-1)^{n(m-1)}=(-1)^{n(m-1)}F^{(m)}_{-1}F^{(m)}_{-1}.$$
Consequently,$$det\textbf{H}^{(m)}_{n}(-1,-1)=(-1)^{n(m-1)}F^{(m)}_{-1}F^{(m)}_{-1}.$$

$(ii)$~ For $\textbf{H}^{(m)}_{n}(-2,-1)$, we expand the bottom row
using the recurrence relation:
$F^{(m)}_{n}=F^{(m)}_{n+m}-\sum^{m-1}_{i=1}F^{(m)}_{n+i}$.

We then perform the following operations:

$\bullet$ Add rows 0 through $m-3$ to the bottom row

$\bullet$ Multiply row $m-2$ by -1 and add it to the bottom row

From the properties of determinants and the relations
$F^{(m)}_{-1}=-F^{(m)}_{-2}=1$, we derive
\begin{eqnarray*}
&&det\textbf{H}^{(m)}_{n}(-2,-1) =det\left(
\begin{array}{ccccccccccccccccc}
F^{(m)}_{n+m-3}&F^{(m)}_{n+m-2}  &F^{(m)}_{n+m-1}   &\cdots&F^{(m)}_{n+2m-4}    \\
        \vdots   &\vdots            &\vdots &      \vdots      &\vdots\\
 F^{(m)}_{n}&F^{(m)}_{n+1}    &F^{(m)}_{n+2}    &\cdots &F^{(m)}_{n+m-1}   \\
 F^{(m)}_{n-1}&F^{(m)}_{n}    &F^{(m)}_{n+1}     &\cdots&F^{(m)}_{n+m-2}   \\
 F^{(m)}_{n-3} &F^{(m)}_{n-2}  &F^{(m)}_{n-1}    &\cdots &F^{(m)}_{n+m-4}    \\
 \end{array}
\right)\\
\end{eqnarray*}
\begin{eqnarray*}
&=&det\left(
\begin{array}{ccccccccccccccccc}
F^{(m)}_{n+m-3}&F^{(m)}_{n+m-2}     &\cdots&F^{(m)}_{n+2m-4}    \\
        \vdots   &\vdots            &      \vdots      &\vdots\\
 F^{(m)}_{n}&F^{(m)}_{n+1}        &\cdots &F^{(m)}_{n+m-1}   \\
 F^{(m)}_{n-1}&F^{(m)}_{n}         &\cdots&F^{(m)}_{n+m-2}   \\
 F^{(m)}_{n-3+m}-\sum^{m-1}_{i=1}F^{(m)}_{n-3+i} &F^{(m)}_{n-2+m}-\sum^{m-1}_{i=1}F^{(m)}_{n-2+i}
   &\cdots &F^{(m)}_{n-4+2m}-\sum^{m-1}_{i=1}F^{(m)}_{n-4+m+i}   \\
 \end{array}
\right)\\&=&-det\left(
\begin{array}{ccccccccccccccccc}
F^{(m)}_{n+m-3}&F^{(m)}_{n+m-2}  &F^{(m)}_{n+m-1}   &\cdots&F^{(m)}_{n+2m-4}    \\
        \vdots   &\vdots            &\vdots &      \vdots      &\vdots\\
 F^{(m)}_{n}&F^{(m)}_{n+1}    &F^{(m)}_{n+2}    &\cdots &F^{(m)}_{n+m-1}   \\
 F^{(m)}_{n-1}&F^{(m)}_{n}    &F^{(m)}_{n+1}     &\cdots&F^{(m)}_{n+m-2}   \\
 F^{(m)}_{n-2} &F^{(m)}_{n-1}  &F^{(m)}_{n}    &\cdots &F^{(m)}_{n+m-3}    \\
 \end{array}
\right)\\&=&-det\textbf{H}^{(m)}_{n}(-1,-1)=(-1)^{n(m-1)}(-F^{(m)}_{-1})F^{(m)}_{-1}=
(-1)^{n(m-1)}F^{(m)}_{-2}F^{(m)}_{-1}.\\
\end{eqnarray*}

\vspace{1.5mm}\noindent{\bf Lemma 4.9}~ {\it For all $n \in
\mathbb{Z}$, $d \in \mathbb{Z}^{+}$, and $l=1$, if $3 \leq d \leq
m$, then}
$$det\textbf{H}^{(m)}_{n}(-d,-1)=0=(-1)^{n(m-1)}F^{(m)}_{-d}F^{(m)}_{-1}.$$

\vspace{1.5mm}{\bf Proof}~ Expanding the bottom row of
$\textbf{H}^{(m)}_{n}(-d,-1)$ using
$F^{(m)}_{n}=2F^{(m)}_{n+m}-F^{(m)}_{n+m+1}$, and noting that for
$3\leq d\leq m$, the bottom row must be equal to twice some upper
row minus another row, we obtain
\begin{eqnarray*}
&&det\textbf{H}^{(m)}_{n}(-d,-1) =det\left(
\begin{array}{ccccccccccccccccc}
F^{(m)}_{n+m-3}&F^{(m)}_{n+m-2}  &F^{(m)}_{n+m-1}   &\cdots&F^{(m)}_{n+2m-4}    \\
        \vdots   &\vdots            &\vdots &      \vdots      &\vdots\\
 F^{(m)}_{n}&F^{(m)}_{n+1}    &F^{(m)}_{n+2}    &\cdots &F^{(m)}_{n+m-1}   \\
 F^{(m)}_{n-1}&F^{(m)}_{n}    &F^{(m)}_{n+1}     &\cdots&F^{(m)}_{n+m-2}   \\
 F^{(m)}_{n-1-d} &F^{(m)}_{n-d}  &F^{(m)}_{n+1-d}    &\cdots &F^{(m)}_{n+m-d}    \\
 \end{array}
\right)\\
\end{eqnarray*}
\begin{eqnarray*}
&=&det\left(
\begin{array}{ccccccccccccccccc}
F^{(m)}_{n+m-3}&F^{(m)}_{n+m-2}     &\cdots&F^{(m)}_{n+2m-4}    \\
        \vdots   &\vdots            &      \vdots      &\vdots\\
 F^{(m)}_{n}&F^{(m)}_{n+1}        &\cdots &F^{(m)}_{n+m-1}   \\
 F^{(m)}_{n-1}&F^{(m)}_{n}         &\cdots&F^{(m)}_{n+m-2}   \\
 2F^{(m)}_{n-1-d+m}-F^{(m)}_{n-d+m} &2F^{(m)}_{n-d+m}-F^{(m)}_{n-d+m+1}
 &\cdots &2F^{(m)}_{n+2m-d}-F^{(m)}_{n+2m-d+1}   \\
 \end{array}
\right)\\&=&-det\left(
\begin{array}{ccccccccccccccccc}
F^{(m)}_{n+m-3}&F^{(m)}_{n+m-2}  &F^{(m)}_{n+m-1}   &\cdots&F^{(m)}_{n+2m-4}    \\
        \vdots   &\vdots            &\vdots &      \vdots      &\vdots\\
 F^{(m)}_{n}&F^{(m)}_{n+1}    &F^{(m)}_{n+2}    &\cdots &F^{(m)}_{n+m-1}   \\
 F^{(m)}_{n-1}&F^{(m)}_{n}    &F^{(m)}_{n+1}     &\cdots&F^{(m)}_{n+m-2}   \\
 0 &0  &0    &\cdots &0    \\
 \end{array}
\right)\\&=&0=(-1)^{n(m-1)}F^{(m)}_{-d}F^{(m)}_{-1}.
\end{eqnarray*}

\vspace{1.5mm}\noindent{\bf Lemma 4.10}~ {\it For all $n \in
\mathbb{Z}$, $d \in \mathbb{Z}^{+}$, and $l=1$, we have}
$$det\textbf{H}^{(m)}_{n}(-d,-1)=(-1)^{n(m-1)}F^{(m)}_{-d}F^{(m)}_{-1}.$$

\vspace{1.5mm}{\bf Proof}~ For $1\leq d\leq m$, the result has been
established in Lemmas 4.8 and 4.9. Assume the proposition holds for
$d=k \geq m$. For $d=k+1$, we expand the bottom row using the
recurrence relation
$F^{(m)}_{n}=F^{(m)}_{n+m}-\sum^{m-1}_{i=1}F^{(m)}_{n+i}$, and by
applying properties of determinants, we obtain
\begin{eqnarray*}
&&det\textbf{H}^{(m)}_{n}\big(-(k+1, -1)\big)=det\left(
\begin{array}{ccccccccccccccccc}
F^{(m)}_{n+m-3}&F^{(m)}_{n+m-2}  &F^{(m)}_{n+m-1}   &\cdots&F^{(m)}_{n+2m-4}    \\
        \vdots &\vdots           & \vdots         & \vdots   &\vdots\\
 F^{(m)}_{n}  &F^{(m)}_{n+1}    &F^{(m)}_{n+2}    &\cdots &F^{(m)}_{n+m-1}   \\
 F^{(m)}_{n-1}&F^{(m)}_{n}    &F^{(m)}_{n+1}     &\cdots&F^{(m)}_{n+m-2}   \\
 F^{(m)}_{n-1-(k+1)} &F^{(m)}_{n-(k+1)}  &F^{(m)}_{n+1-(k+1)}    &\cdots &F^{(m)}_{n+m-(k+1)}    \\
 \end{array}
\right)\\
 &=& det\left(
\begin{array}{ccccccccccccccccc}
F^{(m)}_{n+m-3}&F^{(m)}_{n+m-2}   &\cdots   &F^{(m)}_{n+2m-4}    \\
        \vdots  &\vdots             & \vdots      &\vdots\\
 F^{(m)}_{n}&F^{(m)}_{n+1}       &\cdots &F^{(m)}_{n+m-1}   \\
 F^{(m)}_{n-1}&F^{(m)}_{n}        &\cdots&F^{(m)}_{n+m-2}   \\
F^{(m)}_{n-k-2+m} &F^{(m)}_{n-k-1+m}    &\cdots &F^{(m)}_{n-k-1+2m}    \\
 \end{array}
\right)\\
&&-det\left(
\begin{array}{ccccccccccccccccc}
F^{(m)}_{n+m-3}&F^{(m)}_{n+m-2}   &\cdots   &F^{(m)}_{n+2m-4}    \\
        \vdots  &\vdots             & \vdots      &\vdots\\
 F^{(m)}_{n}&F^{(m)}_{n+1}       &\cdots &F^{(m)}_{n+m-1}   \\
 F^{(m)}_{n-1}&F^{(m)}_{n}        &\cdots&F^{(m)}_{n+m-2}   \\
\sum^{m-1}_{i=1}F^{(m)}_{n-k-2+i} &\sum^{m-1}_{i=1}F^{(m)}_{n-k-1+i}    &\cdots &\sum^{m-1}_{i=1}F^{(m)}_{n-k-1+m+i}    \\
 \end{array}
\right)\\
\end{eqnarray*}
\begin{eqnarray*}
&=&det\textbf{H}^{(m)}_{n}\big(-(k+1-m), -1\big)-
\sum^{m-1}_{i=1}det\textbf{H}^{(m)}_{n}\big(-(k+1-i), -1\big)\\
&=&(-1)^{n(m-1)}F^{(m)}_{-(k+1-m)}F^{(m)}_{-1}-\sum^{m-1}_{i=1}\Big((-1)^{n(m-1)}F^{(m)}_{-(k+1-i)F^{(m)}_{-1}}\Big)\\
&=&(-1)^{n(m-1)}F^{(m)}_{-1}\Big(F^{(m)}_{-(k+1-m)}-\sum^{m-1}_{i=1}F^{(m)}_{-(k+1-i)}\Big)
=(-1)^{n(m-1)}F^{(m)}_{-(k+1)}F^{(m)}_{-1}.\\
\end{eqnarray*}
Thus, by the principle of mathematical induction, the statement
holds for all $d \in \mathbb{Z}^{+}$.

 For all $\forall~ n \in \mathbb{Z}$, $ d \in \mathbb{Z}^{+}$, following
analogous reasoning to Lemma 4.8(ii) and Lemma 4.9 regarding row
distance results and their proofs, we obtain for column distance
$l$:

When column distance $l=2$:
$$det\textbf{H}^{(m)}_{n}(-d,-2)=(-1)^{n(m-1)}F^{(m)}_{-d}F^{(m)}_{-2};$$

When column distance $3\leq l\leq m$:
$$det\textbf{H}^{(m)}_{n}(-d,-l)=(-1)^{n(m-1)}F^{(m)}_{-d}F^{(m)}_{-l}.$$
Consequently, for all $n \in \mathbb{Z}$, $d, l \in \mathbb{Z}^{+}$,
by applying mathematical induction on column distance $l$ similar to
Lemma 4.10, we establish the following Theorem 4.3.

\vspace{1.5mm}\noindent{\bf Theorem 4.3}~ {\it For all $n \in
\mathbb{Z}$ and $d, l \in \mathbb{Z}^{+}$, we have}
$$det\textbf{H}^{(m)}_{n}(-d,-l)=(-1)^{n(m-1)}F^{(m)}_{-d}F^{(m)}_{-l}.$$

\subsection*{4.2.4 \ Properties of
$det\textbf{H}^{(m)}_{n}(-d,r)$}

\noindent{\bf Theorem 4.4}~ {\it For all $n \in \mathbb{Z}$, $d \in
\mathbb{Z}^{+}$, and $r \in \mathbb{N}$, we have}
$$det\textbf{H}^{(m)}_{n}(-d,r)=(-1)^{(n-1)(m-1)}F^{(m)}_{-d}F^{(m)}_{r}.$$

\vspace{1.5mm}{\bf Proof} When $0\leq r\leq m-2$, the matrix
$\textbf{H}^{(m)}_{n}(-d,r)$ contains two identical columns (also
called a quasi full-dimensional submatrix), hence
$det\textbf{H}^{(m)}_{n}(-d,r)=0.$ Since
$F^{(m)}_{0}=F^{(m)}_{1}=F^{(m)}_{2}=\cdots=F^{(m)}_{m-2}=0,$ it
follows that for $0\leq r\leq m-2$,
$$det\textbf{H}^{(m)}_{n}(-d,r)=(-1)^{(n-1)(m-1)}F^{(m)}_{-d}F^{(m)}_{r}.$$

For the case $ r= m-1$, we expand the rightmost column using the
recurrence relation $F^{(m)}_{n}=2F^{(m)}_{n-1}-F^{(m)}_{n-m-1}$,
which yields
\begin{eqnarray*}
det\textbf{H}^{(m)}_{n}(-d,m-1) &=&det\left(
\begin{array}{ccccccccccccccccc}
F^{(m)}_{n+m-2}  &F^{(m)}_{n+m-1}    &\cdots &F^{(m)}_{n+2m-4}   & F^{(m)}_{n+2m-3} \\
F^{(m)}_{n+m-3}  &F^{(m)}_{n+m-2}   &\cdots&F^{(m)}_{n+2m-5}   & F^{(m)}_{n+2m-4} \\
        \vdots   &\vdots            &\vdots &      \vdots      &\vdots\\
F^{(m)}_{n+1}    &F^{(m)}_{n+2}    &\cdots &F^{(m)}_{n+m-1}     & F^{(m)}_{n+m} \\
F^{(m)}_{n}      &F^{(m)}_{n+1}    &\cdots &F^{(m)}_{n+m-2}    & F^{(m)}_{n+m-1}\\
F^{(m)}_{n-d}   &F^{(m)}_{n+1-d}   &\cdots&F^{(m)}_{n+m-2-d}   & F^{(m)}_{n+m-1-d}\\
 \end{array}
\right)\\&=&det\left(
\begin{array}{ccccccccccccccccc}
F^{(m)}_{n+m-2}  &F^{(m)}_{n+m-1}    &\cdots &F^{(m)}_{n+2m-4}   & 2F^{(m)}_{n+2m-4}-F^{(m)}_{n+m-4} \\
F^{(m)}_{n+m-3}  &F^{(m)}_{n+m-2}   &\cdots&F^{(m)}_{n+2m-5}   & 2F^{(m)}_{n+2m-5}-F^{(m)}_{n+m-5} \\
        \vdots   &\vdots            &\vdots &      \vdots      &\vdots\\
F^{(m)}_{n+1}    &F^{(m)}_{n+2}    &\cdots &F^{(m)}_{n+m-1}     & 2F^{(m)}_{n+m-1}-F^{(m)}_{n-1} \\
F^{(m)}_{n}      &F^{(m)}_{n+1}    &\cdots &F^{(m)}_{n+m-2}    & 2F^{(m)}_{n+m-2}-2F^{(m)}_{n-2}\\
F^{(m)}_{n-d}   &F^{(m)}_{n+1-d}   &\cdots&F^{(m)}_{n+m-2-d}   & 2F^{(m)}_{n+m-2-d}-F^{(m)}_{n-2-d}\\
 \end{array}
\right)\\
&=&-det\left(
\begin{array}{ccccccccccccccccc}
F^{(m)}_{n+m-2}  &F^{(m)}_{n+m-1}    &\cdots &F^{(m)}_{n+2m-4}   & F^{(m)}_{n+m-4} \\
F^{(m)}_{n+m-3}  &F^{(m)}_{n+m-2}   &\cdots&F^{(m)}_{n+2m-5}   & F^{(m)}_{n+m-5} \\
        \vdots   &\vdots            &\vdots &      \vdots      &\vdots\\
F^{(m)}_{n+1}    &F^{(m)}_{n+2}    &\cdots &F^{(m)}_{n+m-1}     & F^{(m)}_{n-1} \\
F^{(m)}_{n}      &F^{(m)}_{n+1}    &\cdots &F^{(m)}_{n+m-2}    & 2F^{(m)}_{n-2}\\
F^{(m)}_{n-d}   &F^{(m)}_{n+1-d}   &\cdots&F^{(m)}_{n+m-2-d}   & F^{(m)}_{n-2-d}\\
 \end{array}
\right).
\end{eqnarray*}
By successively exchanging the rightmost column with each of the
$m-1$ left columns until it reaches the far left position, we apply
Theorem 4.3 to obtain
\begin{eqnarray*}
det\textbf{H}^{(m)}_{n}(-d,m-1) &=&-(-1)^{m-1}det\left(
\begin{array}{ccccccccccccccccc}
 F^{(m)}_{n+m-4}&F^{(m)}_{n+m-2}  &F^{(m)}_{n+m-1}    &\cdots &F^{(m)}_{n+2m-4}    \\
 F^{(m)}_{n+m-5}&F^{(m)}_{n+m-3}  &F^{(m)}_{n+m-2}   &\cdots&F^{(m)}_{n+2m-5}    \\
        \vdots   &\vdots            &\vdots &      \vdots      &\vdots\\
 F^{(m)}_{n-1} &F^{(m)}_{n+1}    &F^{(m)}_{n+2}    &\cdots &F^{(m)}_{n+m-1}      \\
 2F^{(m)}_{n-2} &F^{(m)}_{n}      &F^{(m)}_{n+1}    &\cdots &F^{(m)}_{n+m-2}    \\
 F^{(m)}_{n-2-d}&F^{(m)}_{n-d}   &F^{(m)}_{n+1-d}   &\cdots&F^{(m)}_{n+m-2-d}   \\
 \end{array}
\right)\\&=&-(-1)^{m-1}det\textbf{H}^{(m)}_{n}(-d,-2)=-(-1)^{m-1}(-1)^{n(m-1)}F^{(m)}_{-d}F^{(m)}_{-2}\\
&=&-(-1)^{(n+1)(m-1)}F^{(m)}_{-d}(-F^{(m)}_{m-1})=(-1)^{(n-1)(m-1)}F^{(m)}_{-d}F^{(m)}_{m-1}.
\end{eqnarray*}

Assume the proposition holds for all non-negative integers $0 \leq r
\leq k$ (where $k \geq m-1$). For the case $r = k+1$, we expand the
rightmost column of $\det\textbf{H}^{(m)}_{n}(-d,r)$ using the
recurrence relation: $F^{(m)}_{n}=\sum^{m}_{i=1}F^{(m)}_{n-i}$, and
then apply the properties of determinants along with the induction
hypothesis to obtain
\begin{eqnarray*}
det\textbf{H}^{(m)}_{n}(-d,k+1) &=&det\left(
\begin{array}{ccccccccccccccccc}
F^{(m)}_{n+m-2}  &F^{(m)}_{n+m-1}   &\cdots&F^{(m)}_{n+2m-4}   & F^{(m)}_{n+m+k-1} \\
        \vdots   &\vdots            &\vdots &      \vdots      &\vdots\\
F^{(m)}_{n+2}    &F^{(m)}_{n+3}    &\cdots &F^{(m)}_{n+m}     & F^{(m)}_{n+k+3} \\
F^{(m)}_{n+1}    &F^{(m)}_{n+2}    &\cdots &F^{(m)}_{n+m-1}   & F^{(m)}_{n+k+2}\\
F^{(m)}_{n}      &F^{(m)}_{n+1}     &\cdots&F^{(m)}_{n+m-2}   & F^{(m)}_{n+k+1}\\
F^{(m)}_{n-d}   &F^{(m)}_{n+1-d}   &\cdots&F^{(m)}_{n+m-2-d}   & F^{(m)}_{n-d+k+1}\\
 \end{array}
\right)\\
\end{eqnarray*}
\begin{eqnarray*}
&=&det\left(
\begin{array}{ccccccccccccccccc}
F^{(m)}_{n+m-2}  &F^{(m)}_{n+m-1}   &\cdots&F^{(m)}_{n+2m-4}   & \sum^{m}_{i=1}F^{(m)}_{n+m+k-1-i} \\
        \vdots   &\vdots            &\vdots &      \vdots      &\vdots\\
F^{(m)}_{n+2}    &F^{(m)}_{n+3}    &\cdots &F^{(m)}_{n+m}     & \sum^{m}_{i=1}F^{(m)}_{n+k+3-i} \\
F^{(m)}_{n+1}    &F^{(m)}_{n+2}    &\cdots &F^{(m)}_{n+m-1}   & \sum^{m}_{i=1}F^{(m)}_{n+k+2-i}\\
F^{(m)}_{n}      &F^{(m)}_{n+1}     &\cdots&F^{(m)}_{n+m-2}   & \sum^{m}_{i=1}F^{(m)}_{n+k+1-i}\\
F^{(m)}_{n-d}   &F^{(m)}_{n+1-d}   &\cdots&F^{(m)}_{n+m-2-d}   & \sum^{m}_{i=1}F^{(m)}_{n-d+k+1-i}\\
 \end{array}
\right)\\&=&
\sum^{m}_{i=1}det\textbf{H}^{(m)}_{n}(-d,k+1-i)=\sum^{m}_{i=1}\Big((-1)^{(n-1)(m-1)}F^{(m)}_{-d}F^{(m)}_{k+1-i}\Big)\\
&=&(-1)^{(n-1)(m-1)}F^{(m)}_{-d}\sum^{m}_{i=1}F^{(m)}_{k+1-i}
=(-1)^{(n-1)(m-1)}F^{(m)}_{-d}F^{(m)}_{k+1}.\\
\end{eqnarray*}
The proof of the theorem is complete.

\subsection*{4.2.5~ A Unified Formula for the Determinant of
the $(m-1)$-Consecutive Full-dimensional Submatrix}

\vspace{1.5mm}From the above discussion on the four types of
$(m-1)$-consecutive full-dimensional submatrices, we can see that a
consecutive full-dimensional submatrix ($m$-consecutive) can be
viewed as either a $\urcorner$-type matrix, $\ulcorner$-type matrix,
$\llcorner$-type matrix, or $\lrcorner$-type matrix.

Moreover, the four $(m-1)$-consecutive full-dimensional submatrices
above can be considered as constructed by adding a movable $x$-th
row and a movable $y$-th column to a given principal submatrix
$\textbf{H}^{(m_{-1})}_{n}$. Thus, these four types of
$(m-1)$-consecutive full-dimensional submatrices are uniformly
denoted as $\textbf{H}^{(m)}_{n}(x,y)$, where $x, y \in \mathbb{Z}$
and $(x, y)$ corresponds to $(u, r), ~(u,-l), ~(-d,-l),$ or $
~(-d,r)$ according to the following mapping rule:
\begin{eqnarray*}
\left\{\begin{array}{cccccccc}
(x\geq0,y\geq0) &\Longleftrightarrow &(u,r)\\
(x\geq0,y<0) &\Longleftrightarrow &(u,-l)\\
(x<0,y<0)   &\Longleftrightarrow& (-d,-l) \\
(x<0,y\geq0) &\Longleftrightarrow&(-d,r) \\
\end{array}
 \right.
\end{eqnarray*}

If we further define $h(xy)=sgn(|xy|-xy)$, then all results from
Sections 4.2.1 to 4.2.4 can be unified into the following theorem:

\vspace{1.5mm}\noindent{\bf Theorem 4.5}~ {\it For all $n, x, y \in
\mathbb{Z}$ and $(m-1) \in \mathbb{Z}^{+}$, we have}
$$det\textbf{H}^{(m)}_{n}(x,y)=(-1)^{(n+h(xy))(m-1)}F^{(m)}_{x}F^{(m)}_{y}.$$

So far, we have discussed the determinant properties of the
$(m-1)$-consecutive matrix $\textbf{H}^{(m)}_{n}(x,y)$. As for the
properties of these types of matrices themselves, further
investigation is required.

 \subsection*{4.3 \ Other Types of Non-consecutive
Full-dimensional Submatrices of $\mathcal{H}^{(m)}$}

From our previous discussion, we observe that for $m=2$, any
non-consecutive full-dimensional submatrix of $\mathcal{H}^{(2)}$
can be viewed as the matrix $\textbf{H}^{(2)}_{n}(x,y)$. Therefore,
our results provide a satisfactory solution for computing
determinants of non-consecutive full-dimensional submatrices of
$\mathcal{H}^{(2)}$. However, for $m\geq 3$, the non-consecutive
full-dimensional submatrices of $\mathcal{H}^{(m)}$ exhibit much
greater diversity. Their general properties remain largely unknown.
To compute the determinant of an arbitrary non-consecutive
full-dimensional submatrix of $\mathcal{H}^{(m)}$, we typically need
to employ Definition 3.1, the recurrence formula from Theorem 3.1 or
even Corollary 3.1, along with determinant properties, to reduce it
to a linear combination of $det\textbf{H}^{(m)}_{n}(x,y)$. This
approach is illustrated in the following Example:

Using the recurrence formula from Corollary 3.1, the determinant
\begin{eqnarray*} &&\left|
\begin{array}{ccccccccccccccccc}
 F^{(3)}_{n-3}  &F^{(3)}_{n+2}    &F^{(3)}_{n+7}          \\
 F^{(3)}_{n-4}   &F^{(3)}_{n+1}    &F^{(3)}_{n+6}      \\
 F^{(3)}_{n-5}  &F^{(3)}_{n}   &F^{(3)}_{n+5}      \\
 \end{array}
\right|=\left|
\begin{array}{ccccccccccccccccc}
 8F^{(3)}_{n+6}-8F^{(3)}_{n+7}+2F^{(3)}_{n+8}-F^{(3)}_{n+1} &F^{(3)}_{n+2}    &F^{(3)}_{n+7}          \\
 8F^{(3)}_{n+5}-8F^{(3)}_{n+6}+2F^{(3)}_{n+7}-F^{(3)}_{n}   &F^{(3)}_{n+1}    &F^{(3)}_{n+6}      \\
 8F^{(3)}_{n+4}-8F^{(3)}_{n+5}+2F^{(3)}_{n+6}-F^{(3)}_{n-1}  &F^{(3)}_{n}   &F^{(3)}_{n+5}      \\
 \end{array}
\right|\\
\end{eqnarray*}
\begin{eqnarray*}
&=&-8\left|
\begin{array}{ccccccccccccccccc}
 F^{(3)}_{n+2}  &F^{(3)}_{n+6}   &F^{(3)}_{n+7}          \\
 F^{(3)}_{n+1} &F^{(3)}_{n+5}      &F^{(3)}_{n+6}      \\
 F^{(3)}_{n}   &F^{(3)}_{n+4}   &F^{(3)}_{n+5}      \\
 \end{array}
\right|+2\left|
\begin{array}{ccccccccccccccccc}
 F^{(3)}_{n+2}  &F^{(3)}_{n+7}   &F^{(3)}_{n+8}          \\
 F^{(3)}_{n+1} &F^{(3)}_{n+6}      &F^{(3)}_{n+7}      \\
 F^{(3)}_{n}   &F^{(3)}_{n+5}   &F^{(3)}_{n+6}      \\
 \end{array}
\right|-\left|
\begin{array}{ccccccccccccccccc}
 F^{(3)}_{n+1}  &F^{(3)}_{n+2}   &F^{(3)}_{n+7}          \\
 F^{(3)}_{n} &F^{(3)}_{n+1}      &F^{(3)}_{n+6}      \\
 F^{(3)}_{n-1}   &F^{(3)}_{n}   &F^{(3)}_{n+5}      \\
 \end{array}
\right|\\
&=&-8det\textbf{H}^{(3)}_{n+4}(2,-4)+2det\textbf{H}^{(3)}_{n+5}(2,-5)-det\textbf{H}^{(3)}_{n-1}(2,6)\\
&=&-8(-1)^{(n+4+h(2(-4)))(3-1)}F^{(3)}_{2}F^{(3)}_{-4}+2(-1)^{(n+5+h(2(-5)))(3-1)}F^{(3)}_{2}F^{(3)}_{-5}\\
&&-(-1)^{(n-1+h(2\cdot6))(3-1)}F^{(3)}_{2}F^{(3)}_{6}=-29.
\end{eqnarray*}

Next, we examine several other special types of non-consecutive
full-dimensional submatrices of $\bf\mathcal{H}^{(m)}$. First, by
applying the recurrence formula from Corollary 3.1 and utilizing
determinant properties, we can reduce them to determinants of either
$\urcorner$-type or $\llcorner$-type matrices (alternatively, they
may also be reduced to $\ulcorner$-type or $\lrcorner$-type matrix
determinants). These can then be computed using Theorem 4.5, as
demonstrated in the following examples:

\begin{eqnarray*} 1^{\circ}&&\left|
\begin{array}{ccccccccccccccccc}
 F^{(3)}_{n+2} &F^{(3)}_{n+4}    &F^{(3)}_{n+6}          \\
 F^{(3)}_{n} &F^{(3)}_{n+2}      &F^{(3)}_{n+4}      \\
 F^{(3)}_{n-2}&F^{(3)}_{n}   &F^{(3)}_{n+2}      \\
 \end{array}
\right|=\left|
\begin{array}{ccccccccccccccccc}
 2F^{(3)}_{n+5}-F^{(3)}_{n+6} &F^{(3)}_{n+4}    &F^{(3)}_{n+6}          \\
 2F^{(3)}_{n+3}-F^{(3)}_{n+4} &F^{(3)}_{n+2}      &F^{(3)}_{n+4}      \\
 2F^{(3)}_{n+1}-F^{(3)}_{n+2}&F^{(3)}_{n}   &F^{(3)}_{n+2}      \\
 \end{array}
\right|=-2\left|
\begin{array}{ccccccccccccccccc}
 F^{(3)}_{n+4}  &F^{(3)}_{n+5}   &F^{(3)}_{n+6}          \\
 F^{(3)}_{n+2} &F^{(3)}_{n+3}      &F^{(3)}_{n+4}      \\
 F^{(3)}_{n} &F^{(3)}_{n+1}   &F^{(3)}_{n+2}      \\
 \end{array}
\right|\hspace{2cm}\\&=&-2\left|
\begin{array}{ccccccccccccccccc}
 F^{(3)}_{n+4}  &F^{(3)}_{n+5}   &F^{(3)}_{n+6}          \\
 F^{(3)}_{n+2} &F^{(3)}_{n+3}      &F^{(3)}_{n+4}      \\
 2F^{(3)}_{n+3}-F^{(3)}_{n+4} &2F^{(3)}_{n+4}-F^{(3)}_{n+5}   &2F^{(3)}_{n+5}-F^{(3)}_{n+6}      \\
 \end{array}
\right|=4\left|
\begin{array}{ccccccccccccccccc}
 F^{(3)}_{n+4}  &F^{(3)}_{n+5}   &F^{(3)}_{n+6}          \\
  F^{(3)}_{n+3} &F^{(3)}_{n+4}   &F^{(3)}_{n+5}      \\
 F^{(3)}_{n+2} &F^{(3)}_{n+3}      &F^{(3)}_{n+4}      \\
 \end{array}
\right|\\&=&4det\textbf{H}^{(3)}_{n+2}(2,2)=4(-1)^{(n+2+h(2\cdot2))(3-1)}F^{(3)}_{2}F^{(3)}_{2}=4.
\end{eqnarray*}
Similarly, we obtain
\begin{eqnarray*}
2^{\circ}&&\left|
\begin{array}{ccccccccccccccccc}
 F^{(3)}_{n+3} &F^{(3)}_{n+6}    &F^{(3)}_{n+9}          \\
 F^{(3)}_{n} &F^{(3)}_{n+3}      &F^{(3)}_{n+6}      \\
 F^{(3)}_{n-3}&F^{(3)}_{n}   &F^{(3)}_{n+3}      \\
 \end{array}
\right|=\left|
\begin{array}{ccccccccccccccccc}
 F^{(3)}_{n+3} &2F^{(3)}_{n+9}-F^{(3)}_{n+10}     &F^{(3)}_{n+9}          \\
 F^{(3)}_{n}  &2F^{(3)}_{n+6}-F^{(3)}_{n+7}       &F^{(3)}_{n+6}      \\
 F^{(3)}_{n-3}&2F^{(3)}_{n+3}-F^{(3)}_{n+4}   &F^{(3)}_{n+3}      \\
 \end{array}
\right|=\left|
\begin{array}{ccccccccccccccccc}
 F^{(3)}_{n+3}    &F^{(3)}_{n+9}&F^{(3)}_{n+10}            \\
 F^{(3)}_{n}         &F^{(3)}_{n+6}&F^{(3)}_{n+7}      \\
 F^{(3)}_{n-3}  &F^{(3)}_{n+3} &F^{(3)}_{n+4}      \\
 \end{array}
\right|\hspace{2cm}\\&=&\left|
\begin{array}{ccccccccccccccccc}
 F^{(3)}_{n+4} &F^{(3)}_{n+10}   &F^{(3)}_{n+11}      \\
  F^{(3)}_{n+3}    &F^{(3)}_{n+9}&F^{(3)}_{n+10}            \\
 F^{(3)}_{n-3}  &F^{(3)}_{n+3} &F^{(3)}_{n+4}      \\
 \end{array}
\right|=det\textbf{H}^{(3)}_{n+9}(-6,-6)=(-1)^{(n+9+h(-6\cdot(-6)))(3-1)}F^{(3)}_{-6}F^{(3)}_{-6}=1.
\end{eqnarray*}
\begin{eqnarray*}
3^{\circ}&&\left|
\begin{array}{ccccccccccccccccc}
 F^{(3)}_{n+4} &F^{(3)}_{n+8}    &F^{(3)}_{n+12}          \\
 F^{(3)}_{n} &F^{(3)}_{n+4}      &F^{(3)}_{n+8}      \\
 F^{(3)}_{n-4}&F^{(3)}_{n}   &F^{(3)}_{n+4}      \\
 \end{array}
\right|=4det\textbf{H}^{(3)}_{n+10+h(-7\cdot(-7))}(-7,-7)=4(-1)^{(n+10)(3-1)}F^{(3)}_{-7}F^{(3)}_{-7}=64.\hspace{3cm}
\end{eqnarray*}
$4^{\circ}$
\begin{eqnarray*}
 \left|
\begin{array}{ccccccccccccccccc}
 F^{(3)}_{n+3} &F^{(3)}_{n+7}    &F^{(3)}_{n+10}          \\
 F^{(3)}_{n} &F^{(3)}_{n+4}      &F^{(3)}_{n+7}      \\
 F^{(3)}_{n-4}&F^{(3)}_{n}   &F^{(3)}_{n+3}      \\
 \end{array}
\right|&=&\left|
\begin{array}{ccccccccccccccccc}
 2^{2}F^{(3)}_{n+9}-2F^{(3)}_{n+10}-F^{(3)}_{n+7}&  F^{(3)}_{n+7}    &F^{(3)}_{n+10}          \\
 2^{2}F^{(3)}_{n+6}-2F^{(3)}_{n+7}-F^{(3)}_{n+4} & F^{(3)}_{n+4}      &F^{(3)}_{n+7}      \\
 2^{2}F^{(3)}_{n+2}-2F^{(3)}_{n+3}-F^{(3)}_{n}   & F^{(3)}_{n}      &F^{(3)}_{n+3}      \\
 \end{array}
\right|\\
&=&-2^{2}\left|
\begin{array}{ccccccccccccccccc}
F^{(3)}_{n+7}    &F^{(3)}_{n+9}&  F^{(3)}_{n+10}          \\
F^{(3)}_{n+4}      & F^{(3)}_{n+6} & F^{(3)}_{n+7}      \\
F^{(3)}_{n}      & F^{(3)}_{n+2}   & F^{(3)}_{n+3}      \\
 \end{array}
\right|\\
&=&2^{2}2^{2}\left|
\begin{array}{ccccccccccccccccc}
 F^{(3)}_{n+7} &F^{(3)}_{n+9}   &F^{(3)}_{n+10}      \\
  F^{(3)}_{n+6}    &F^{(3)}_{n+8}&F^{(3)}_{n+9}            \\
 F^{(3)}_{n+4}  &F^{(3)}_{n+6} &F^{(3)}_{n+7}      \\
 \end{array}
\right|\\
&=&2^{2}2^{2}det\textbf{H}^{(3)}_{n+8}(-2,-2)=(-1)^{(n+8+h(-2\cdot(-2)))(3-1)}F^{(3)}_{-2}F^{(3)}_{-2}=16.
\end{eqnarray*}

$5^{\circ}$\begin{eqnarray*}
&&\left|
\begin{array}{ccccccccccccccccc}
 F^{(4)}_{n+8} &F^{(4)}_{n+13}    &F^{(4)}_{n+17}  &F^{(4)}_{n+21}        \\
 F^{(4)}_{n+4} &F^{(4)}_{n+9}    &F^{(4)}_{n+13}  &F^{(4)}_{n+17}        \\
 F^{(4)}_{n}  &F^{(4)}_{n+5}      &F^{(4)}_{n+9}    &F^{(4)}_{n+13}   \\
 F^{(4)}_{n-5}&F^{(4)}_{n}       &F^{(4)}_{n+4}    &F^{(4)}_{n+8}  \\
 \end{array}
\right|\\&=&\left|
\begin{array}{ccccccccccccccccc}
2^{3}F^{(4)}_{n+20}-2^{2}F^{(4)}_{n+21}-2F^{(4)}_{n+17}-F^{(4)}_{n+13} &F^{(4)}_{n+13}    &F^{(4)}_{n+17}  &F^{(4)}_{n+21}        \\
2^{3}F^{(4)}_{n+16}-2^{2}F^{(4)}_{n+17}-2F^{(4)}_{n+13}-F^{(4)}_{n+9} &F^{(4)}_{n+9}    &F^{(4)}_{n+13}  &F^{(4)}_{n+17}        \\
2^{3}F^{(4)}_{n+12}-2^{2}F^{(4)}_{n+13}-2F^{(4)}_{n+9}-F^{(4)}_{n+5} &F^{(4)}_{n+5}      &F^{(4)}_{n+9}    &F^{(4)}_{n+13}   \\
2^{3}F^{(4)}_{n+7}-2^{2}F^{(4)}_{n+8}-2F^{(4)}_{n+4}-F^{(4)}_{n}&F^{(4)}_{n}       &F^{(4)}_{n+4}    &F^{(4)}_{n+8}  \\
 \end{array}
\right|\\
\end{eqnarray*}
\begin{eqnarray*}
&&=2^{3}\left|
\begin{array}{ccccccccccccccccc}
F^{(4)}_{n+13}    &F^{(4)}_{n+17}  &F^{(4)}_{n+20} &F^{(4)}_{n+21}        \\
F^{(4)}_{n+9}    &F^{(4)}_{n+13}  &F^{(4)}_{n+16} &F^{(4)}_{n+17}        \\
F^{(4)}_{n+5}      &F^{(4)}_{n+9} &F^{(4)}_{n+12}  &F^{(4)}_{n+13}   \\
F^{(4)}_{n}       &F^{(4)}_{n+4}  &F^{(4)}_{n+7}  &F^{(4)}_{n+8}  \\
 \end{array}
\right|=2^{3}\left|
\begin{array}{ccccccccccccccccc}
F^{(4)}_{n+13}   &2F^{(4)}_{n+21}-F^{(4)}_{n+22}  &F^{(4)}_{n+20} &F^{(4)}_{n+21}        \\
F^{(4)}_{n+9}    &2F^{(4)}_{n+17}-F^{(4)}_{n+18}  &F^{(4)}_{n+16} &F^{(4)}_{n+17}        \\
F^{(4)}_{n+5}    &2F^{(4)}_{n+13}-F^{(4)}_{n+14}   &F^{(4)}_{n+12}  &F^{(4)}_{n+13}   \\
F^{(4)}_{n}      &2F^{(4)}_{n+8}-F^{(4)}_{n+9}  &F^{(4)}_{n+7}  &F^{(4)}_{n+8}  \\
 \end{array}
\right|\\&&=-2^{3}\left|
\begin{array}{ccccccccccccccccc}
F^{(4)}_{n+13}     &F^{(4)}_{n+20} &F^{(4)}_{n+21}  &F^{(4)}_{n+22}      \\
F^{(4)}_{n+9}      &F^{(4)}_{n+16} &F^{(4)}_{n+17}  &F^{(4)}_{n+18}      \\
F^{(4)}_{n+5}      &F^{(4)}_{n+12}  &F^{(4)}_{n+13} &F^{(4)}_{n+14}  \\
F^{(4)}_{n}        &F^{(4)}_{n+7}  &F^{(4)}_{n+8}   &F^{(4)}_{n+9}  \\
 \end{array}
\right|=2^{3}2^{3}\left|
\begin{array}{ccccccccccccccccc}
F^{(4)}_{n+14}      &F^{(4)}_{n+21} &F^{(4)}_{n+22}  &F^{(4)}_{n+23}      \\
F^{(4)}_{n+13}      &F^{(4)}_{n+20} &F^{(4)}_{n+21}  &F^{(4)}_{n+22}      \\
F^{(4)}_{n+12}      &F^{(4)}_{n+19}  &F^{(4)}_{n+20} &F^{(4)}_{n+21}  \\
F^{(4)}_{n+5}       &F^{(4)}_{n+12}  &F^{(4)}_{n+13}   &F^{(4)}_{n+14}  \\
 \end{array}
\right|\\&&=2^{3}2^{3}det\textbf{H}^{(4)}_{n+19}(-7,-7)=2^{3}2^{3}(-1)^{(n+19+h(-7\cdot(-7)))
(4-1)}F^{(4)}_{-7}F^{(4)}_{-7}=(-1)^{3(n+1)}64.
\end{eqnarray*}

\subsection*{4.4 \ Dimensional-reducing Submatrices and
Dimensional-augmented Submatrices}

\subsection*{4.4.1 \ Dimensional-reducing Submatrices}

For the Fibonacci matrix $\mathcal{H}^{(m)}$ of the $m$-dimensional
harmonic-type golden section $T^{(m)}$, for all $n\in \mathbb{Z}$,
$(m-1)\in \mathbb{Z}^{+}$ and $k\in \mathbb{Z}^{+}$, when $1\leq
k\leq m-1$, its dimensional-reducing submatrix
$\textbf{H}^{(m_{-k})}_{n}$ is an $(m-k)$-order square matrix. When
$k=m-1$, $\textbf{H}^{(m_{-k})}_{n}$ becomes a first-order matrix
$(F^{(m)}_{n})$, which is a trivial case and requires no further
discussion. When $m=2$, the dimensional-reducing submatrix of
$\mathcal{H}^{(2)}$ is only the trivial case of the first-order
matrix $(F^{(2)}_{n})$.

Below, we consider a special case where $m\geq 3$ and $k=1$: the
principal submatrix (an $(m-1)$-order consecutive submatrix)
$\textbf{H}^{(m_{-1})}_{n}$ of the $(m-1)$-consecutive
full-dimensional submatrix $\textbf{H}^{(m)}_{n}(x, y)$ of
$\mathcal{H}^{(m)}$. For this, we have the following lemmas and
theorems.

\vspace{1.5mm}\noindent{\bf Lemma 4.11} ~ {\it Let $(m-2)\in
\mathbb{Z}^{+}$, then}
$$(i)\quad\quad\det\textbf{H}^{(m_{-1})}_{0}=(-1)^{(0-1)(m-1)}F^{(m)}_{-0},\quad\quad\quad\quad\quad\quad~~
(ii)\quad\quad\det\textbf{H}^{(m_{-1})}_{1}=(-1)^{(1-1)(m-1)}F^{(m)}_{-1}.$$

\vspace{1.5mm}{\bf Proof}
\begin{eqnarray*}(i)\hspace{-0.5cm}
&&det\textbf{H}^{(m_{-1})}_{0}\\
&&=\left|
\begin{array}{ccccccccccccccccc}
F^{(m)}_{0+m-2}  &F^{(m)}_{0+m-1}   &\cdots&F^{(m)}_{0+2m-4}   \\
        \vdots   &\vdots            &\ddots &      \vdots      \\
F^{(m)}_{0+2}    &F^{(m)}_{0+3}    &\cdots &F^{(m)}_{0+m}      \\
F^{(m)}_{0+1}    &F^{(m)}_{0+2}    &\cdots &F^{(m)}_{0+m-1}   \\
F^{(m)}_{0}      &F^{(m)}_{0+1}     &\cdots&F^{(m)}_{0+m-2}   \\
 \end{array}
\right|=\left|
\begin{array}{ccccccccccccccccc}
 0 &F^{(m)}_{m-1}   &\cdots&F^{(m)}_{2m-4}    \\
        \vdots   &\vdots            &\ddots &      \vdots      \\
 0   & 0   &\cdots &F^{(m)}_{m}      \\
 0   &  0  &\cdots &F^{(m)}_{m-1}   \\
0      & 0    &\cdots& 0  \\
 \end{array}
\right|=
0=(-1)^{(0-1)(m-1)}F^{(m)}_{-0}.\\
\end{eqnarray*}

\begin{eqnarray*}
(ii)~ det\textbf{H}^{(m_{-1})}_{1} =\left|
\begin{array}{ccccccccccccccccc}
F^{(m)}_{m-1}  &F^{(m)}_{m}   &\cdots&F^{(m)}_{2m-2}   \\
        \vdots   &\vdots            &\ddots &      \vdots      \\
F^{(m)}_{3}    &F^{(m)}_{4}    &\cdots &F^{(m)}_{m+1}      \\
F^{(m)}_{2}    &F^{(m)}_{3}    &\cdots &F^{(m)}_{m}   \\
F^{(m)}_{1}      &F^{(m)}_{2}     &\cdots&F^{(m)}_{m-1}   \\
 \end{array}
\right|=\left|
\begin{array}{ccccccccccccccccc}
 1 &F^{(m)}_{m-1}   &\cdots&F^{(m)}_{2m-3}    \\
        \vdots   &\vdots            &\ddots &      \vdots      \\
 0   & 0   &\cdots &F^{(m)}_{m+1}      \\
 0   &  0  &\cdots &F^{(m)}_{m}   \\
0      & 0    &\cdots& 1  \\
 \end{array}
\right|=
1=(-1)^{(1-1)(m-1)}F^{(m)}_{-1}.\\
\end{eqnarray*}

\vspace{1.5mm}\noindent{\bf Lemma 4.12} ~{\it Let $n\in\mathbb{Z},
(m-2)\in \mathbb{Z}^{+}$ and $-(m-2) \leq n \leq -1$, then}
$$det\textbf{H}^{(m_{-1})}_{n}=0=(-1)^{(n-1)(m-1)}F^{(m)}_{-n}.$$

\vspace{1.5mm}{\bf Proof} When $n = -1$, \begin{eqnarray*}
det\textbf{H}^{(m_{-1})}_{-1} &=&\left|
\begin{array}{ccccccccccccccccc}
F^{(m)}_{m-3}  &F^{(m)}_{m-2} &F^{(m)}_{m-1}  &\cdots&F^{(m)}_{2m-5}   \\
        \vdots   &\vdots         &\vdots    &\ddots &      \vdots      \\
F^{(m)}_{0+1}    &F^{(m)}_{2}   &F^{(m)}_{3} &\cdots &F^{(m)}_{m-1}      \\
F^{(m)}_{0}    &F^{(m)}_{1}    &F^{(m)}_{2}&\cdots &F^{(m)}_{m-2}   \\
F^{(m)}_{-1}      &F^{(m)}_{0}  &F^{(m)}_{1}   &\cdots&F^{(m)}_{m-3}   \\
 \end{array}
\right|=\left|
\begin{array}{ccccccccccccccccc}
 0 &0  &1 &\cdots&F^{(m)}_{2m-5}    \\
        \vdots   &\vdots            &\ddots &      \vdots      \\
 0   & 0 &0  &\cdots &1      \\
 0   &  0 &0 &\cdots &0   \\
1      & 0  &0  &\cdots& 0  \\
 \end{array}
\right|\\&=&
0=(-1)^{(-1-1)(m-1)}F^{(m)}_{-(-1)}.\\
\end{eqnarray*}
It is straightforward to verify that
$det\textbf{H}^{(m_{-1})}_{-1}=\cdots=det\textbf{H}^{(m_{-1})}_{-(m-2)}=0.$

\vspace{1.5mm}\noindent{\bf Lemma 4.13} ~{\it Let $(m-2)\in
\mathbb{Z}^{+}$. Then for $n=-(m-1)$ or $n=-m$, we have}
$$det\textbf{H}^{(m_{-1})}_{n}=(-1)^{(n-1)(m-1)}F^{(m)}_{-n}.$$

\vspace{1.5mm}{\bf Proof}~ For $n = -m+1$, \begin{eqnarray*}
det\textbf{H}^{(m_{-1})}_{-m+1} &=&\left|
\begin{array}{ccccccccccccccccc}
F^{(m)}_{-1}  &F^{(m)}_{0} &F^{(m)}_{1}  &\cdots&F^{(m)}_{m-3}   \\
F^{(m)}_{-2}    &F^{(m)}_{-1}   &F^{(m)}_{0} &\cdots &F^{(m)}_{m-4}      \\
        \vdots   &\vdots         &\vdots    &\ddots &      \vdots      \\
F^{(m)}_{-m+2}    &F^{(m)}_{-m+3}    &F^{(m)}_{-m+4}&\cdots &F^{(m)}_{0}   \\
F^{(m)}_{-m+1}      &F^{(m)}_{-m+2}  &F^{(m)}_{-m+3}   &\cdots&F^{(m)}_{-1}   \\
 \end{array}
\right|=\left|
\begin{array}{ccccccccccccccccc}
1 &0  &0 &\cdots&0   \\
-1   & 1 &0  &\cdots &0      \\
  \vdots   &\vdots            &\ddots &      \vdots      \\
  0   &  0 &0 &\cdots &0   \\
0   & 0  &0  &\cdots& 1  \\
 \end{array}
\right|\\&=&
1=(-1)^{(-m+1-1)(m-1)}F^{(m)}_{-(-m+1)}.\\
\end{eqnarray*}

For $n = -m$, \begin{eqnarray*} det\textbf{H}^{(m_{-1})}_{-m}
&=&\left|
\begin{array}{ccccccccccccccccc}
F^{(m)}_{-2}  &F^{(m)}_{-1} &F^{(m)}_{0}  &\cdots&F^{(m)}_{m-4}   \\
F^{(m)}_{-3}    &F^{(m)}_{-2}   &F^{(m)}_{-1} &\cdots &F^{(m)}_{m-5}      \\
        \vdots   &\vdots         &\vdots    &\ddots &      \vdots      \\
F^{(m)}_{-m+1}    &F^{(m)}_{-m+2}    &F^{(m)}_{-m+3}&\cdots &F^{(m)}_{-1}   \\
F^{(m)}_{-m}      &F^{(m)}_{-m+1}  &F^{(m)}_{-m+2}   &\cdots&F^{(m)}_{-2}   \\
 \end{array}
\right|=\left|
\begin{array}{ccccccccccccccccc}
-1 &1  &0 &\cdots&0   \\
0   & -1 &1  &\cdots &0      \\
  \vdots   &\vdots            &\ddots &      \vdots      \\
 0   &  0 &0 &\cdots &1   \\
0   & 0  &0  &\cdots& -1  \\
 \end{array}
\right|\\&=&
(-1)^{(m-1)}=(-1)^{(-m-1)(m-1)}F^{(m)}_{-(-m)}.\\
\end{eqnarray*}

\vspace{1.5mm}\noindent{\bf Theorem 4.6} ~ {\it For all integers $n
\in \mathbb{Z}$ and $(m-2)\in \mathbb{Z}^{+}$, then}
$$det\textbf{H}^{(m_{-1})}_{n}=(-1)^{(n-1)(m-1)}F^{(m)}_{-n}.$$

\vspace{1.5mm}{\bf Proof} ~ We proceed by mathematical induction.
From Lemmas 4.11, 4.12, and 4.13, the statement holds for all
integers $-m \leq n \leq 1$. Below, we separately prove the
statement for $n \geq 1$ and $n \leq -m$.

Let $k \in \mathbb{Z}^{+}$, and assume the statement holds for all
integers $-m \leq n \leq k$. For $n = k + 1$, we perform the
following operations in order from right to left: Expand the
rightmost column using the recurrence relation
$F^{(m)}_{n}=2F^{(m)}_{n-1}-F^{(m)}_{n-m-1}$, then add -2 times the
left column to this column. By properties of determinants, we obtain
\begin{eqnarray*}
det\textbf{H}^{(m_{-1})}_{k+1} &=&det\left(
\begin{array}{ccccccccccccccccc}
F^{(m)}_{k+m-1}  &F^{(m)}_{k+m}   &\cdots&F^{(m)}_{k+2m-5}&F^{(m)}_{k+2m-4} &F^{(m)}_{k+2m-3}  \\
        \vdots   &\vdots            &\vdots &      \vdots&      \vdots   &\vdots   \\
F^{(m)}_{k+3}    &F^{(m)}_{k+4}    &\cdots &F^{(m)}_{k+m-1}  &F^{(m)}_{k+m} &F^{(m)}_{k+m+1}   \\
F^{(m)}_{k+2}    &F^{(m)}_{k+3}    &\cdots&F^{(m)}_{k+m-2} &F^{(m)}_{k+m-1} &F^{(m)}_{k+m}  \\
F^{(m)}_{k+1}    &F^{(m)}_{k+2}     &\cdots&F^{(m)}_{k+m-3} &F^{(m)}_{k+m-2} &F^{(m)}_{k+m-1}   \\
 \end{array}
\right)\\
\end{eqnarray*}
\begin{eqnarray*}
&=&det\left(
\begin{array}{ccccccccccccccccc}
F^{(m)}_{k+m-1}  &F^{(m)}_{k+m}   &\cdots&F^{(m)}_{k+2m-5}&F^{(m)}_{k+2m-4}  & 2F^{(m)}_{k+2m-4}-F^{(m)}_{k+m-4} \\
        \vdots   &\vdots            &\vdots &      \vdots      &\vdots        &\vdots\\
F^{(m)}_{k+3}    &F^{(m)}_{k+4}    &\cdots &F^{(m)}_{k+m-1}  &F^{(m)}_{k+m}     & 2F^{(m)}_{k+m}-F^{(m)}_{k} \\
F^{(m)}_{k+2}    &F^{(m)}_{k+3}    &\cdots&F^{(m)}_{k+m-2} &F^{(m)}_{k+m-1}   & 2F^{(m)}_{k+m-1}-F^{(m)}_{k-1}\\
F^{(m)}_{k+1}    &F^{(m)}_{k+2}     &\cdots&F^{(m)}_{k+m-3} &F^{(m)}_{k+m-2}  & 2F^{(m)}_{k+m-2}-F^{(m)}_{k-2}\\
 \end{array}
\right)\\
\end{eqnarray*}
\begin{eqnarray*}
&=&-det\left(
\begin{array}{ccccccccccccccccc}
F^{(m)}_{k+m-1}  &F^{(m)}_{k+m}   &\cdots&F^{(m)}_{k+2m-5}&2F^{(m)}_{k+2m-5}-F^{(m)}_{k+m-5}  &F^{(m)}_{k+m-4} \\
        \vdots   &\vdots            &\vdots &      \vdots      &\vdots        &\vdots\\
F^{(m)}_{k+3}    &F^{(m)}_{k+4}    &\cdots &F^{(m)}_{k+m-1}  &2F^{(m)}_{k+m-1}-F^{(m)}_{k-1}     &F^{(m)}_{k} \\
F^{(m)}_{k+2}    &F^{(m)}_{k+3}    &\cdots&F^{(m)}_{k+m-2} &2F^{(m)}_{k+m-2}-F^{(m)}_{k-2}   &F^{(m)}_{k-1}\\
F^{(m)}_{k+1}    &F^{(m)}_{k+2}     &\cdots&F^{(m)}_{k+m-3} &2F^{(m)}_{k+m-3}-F^{(m)}_{k-3}  &F^{(m)}_{k-2}\\
 \end{array}
\right)\\
\end{eqnarray*}
\begin{eqnarray*}
&=&(-1)^{2}det\left(
\begin{array}{ccccccccccccccccc}
F^{(m)}_{k+m-1}  &F^{(m)}_{k+m}   &\cdots&2F^{(m)}_{k+2m-6}-F^{(m)}_{k+m-6}&F^{(m)}_{k+m-5}  &F^{(m)}_{k+m-4} \\
        \vdots   &\vdots            &\vdots &      \vdots      &\vdots        &\vdots\\
F^{(m)}_{k+3}    &F^{(m)}_{k+4}    &\cdots &2F^{(m)}_{k+m-2}-F^{(m)}_{k-2} &F^{(m)}_{k-1}    &F^{(m)}_{k} \\
F^{(m)}_{k+2}    &F^{(m)}_{k+3}    &\cdots&2F^{(m)}_{k+m-3}-F^{(m)}_{k-3}  &F^{(m)}_{k-2}   &F^{(m)}_{k-1}\\
F^{(m)}_{k+1}    &F^{(m)}_{k+2}     &\cdots&2F^{(m)}_{k+m-4}-F^{(m)}_{k-4} &F^{(m)}_{k-3}  &F^{(m)}_{k-2}\\
 \end{array}
\right)\\
&=&\cdots\\
\end{eqnarray*}
\begin{eqnarray*}
&=&(-1)^{m-2}det\left(
\begin{array}{ccccccccccccccccc}
F^{(m)}_{k+m-1}  &F^{(m)}_{k-1}   &\cdots&F^{(m)}_{k+m-6}&F^{(m)}_{k+m-5}  &F^{(m)}_{k+m-4} \\
        \vdots   &\vdots            &\vdots &      \vdots      &\vdots        &\vdots\\
F^{(m)}_{k+3}    &F^{(m)}_{k-m+3}    &\cdots &F^{(m)}_{k-2} &F^{(m)}_{k-1}    &F^{(m)}_{k} \\
F^{(m)}_{k+2}    &F^{(m)}_{k-m+2}    &\cdots&F^{(m)}_{k-3}  &F^{(m)}_{k-2}   &F^{(m)}_{k-1}\\
F^{(m)}_{k+1}    &F^{(m)}_{k-m+1}     &\cdots&F^{(m)}_{k-4} &F^{(m)}_{k-3}  &F^{(m)}_{k-2}\\
 \end{array}
\right)\\
\end{eqnarray*}
\begin{eqnarray*}
&=&(-1)^{m-2}det\left(
\begin{array}{ccccccccccccccccc}
F^{(m)}_{k+m-1}  &F^{(m)}_{k-1}   &\cdots&F^{(m)}_{k+m-6}&F^{(m)}_{k+m-5}  &F^{(m)}_{k+m-4} \\
        \vdots   &\vdots            &\vdots &      \vdots      &\vdots        &\vdots\\
F^{(m)}_{k+3}    &F^{(m)}_{k-m+3}    &\cdots &F^{(m)}_{k-2} &F^{(m)}_{k-1}    &F^{(m)}_{k} \\
F^{(m)}_{k+2}    &F^{(m)}_{k-m+2}    &\cdots&F^{(m)}_{k-3}  &F^{(m)}_{k-2}   &F^{(m)}_{k-1}\\
F^{(m)}_{k+1}    &F^{(m)}_{k-m+1}     &\cdots&F^{(m)}_{k-4} &F^{(m)}_{k-3}  &F^{(m)}_{k-2}\\
 \end{array}
\right).\\
\end{eqnarray*}
Expand the leftmost column using the same recurrence relation
$F^{(m)}_{n} = 2F^{(m)}_{n-1} - F^{(m)}_{n-m-1}$. By properties of
determinants, we obtain \begin{eqnarray*}
det\textbf{H}^{(m_{-1})}_{k+1}&=&(-1)^{m-2}det\left(
\begin{array}{ccccccccccccccccc}
2F^{(m)}_{k+m-2}-F^{(m)}_{k-2} &F^{(m)}_{k-1}   &\cdots&F^{(m)}_{k+m-6}&F^{(m)}_{k+m-5}  &F^{(m)}_{k+m-4} \\
        \vdots   &\vdots            &\vdots &      \vdots      &\vdots        &\vdots\\
2F^{(m)}_{k+2}-F^{(m)}_{k-m+2} &F^{(m)}_{k-m+3}    &\cdots &F^{(m)}_{k-2} &F^{(m)}_{k-1}    &F^{(m)}_{k} \\
2F^{(m)}_{k+1}-F^{(m)}_{k-m+1} &F^{(m)}_{k-m+2}    &\cdots&F^{(m)}_{k-3}  &F^{(m)}_{k-2}   &F^{(m)}_{k-1}\\
2F^{(m)}_{k}-F^{(m)}_{k-m} &F^{(m)}_{k-m+1}     &\cdots&F^{(m)}_{k-4} &F^{(m)}_{k-3}  &F^{(m)}_{k-2}\\
 \end{array}
\right)\\
\end{eqnarray*}
\begin{eqnarray*}
&=&(-1)^{m-2}2det\left(
\begin{array}{ccccccccccccccccc}
F^{(m)}_{k+m-2} &F^{(m)}_{k-1}   &\cdots&F^{(m)}_{k+m-6}&F^{(m)}_{k+m-5}  &F^{(m)}_{k+m-4} \\
        \vdots   &\vdots            &\vdots &      \vdots      &\vdots        &\vdots\\
F^{(m)}_{k+2} &F^{(m)}_{k-m+3}    &\cdots &F^{(m)}_{k-2} &F^{(m)}_{k-1}    &F^{(m)}_{k} \\
F^{(m)}_{k+1} &F^{(m)}_{k-m+2}    &\cdots&F^{(m)}_{k-3}  &F^{(m)}_{k-2}   &F^{(m)}_{k-1}\\
F^{(m)}_{k} &F^{(m)}_{k-m+1}     &\cdots&F^{(m)}_{k-4} &F^{(m)}_{k-3}  &F^{(m)}_{k-2}\\
 \end{array}
\right)\\
\end{eqnarray*}
\begin{eqnarray*}
&&+(-1)^{m-1}det\left(
\begin{array}{ccccccccccccccccc}
F^{(m)}_{k-2} &F^{(m)}_{k-1}   &\cdots&F^{(m)}_{k+m-6}&F^{(m)}_{k+m-5}  &F^{(m)}_{k+m-4} \\
        \vdots   &\vdots            &\vdots &      \vdots      &\vdots        &\vdots\\
F^{(m)}_{k-m+2} &F^{(m)}_{k-m+3}    &\cdots &F^{(m)}_{k-2} &F^{(m)}_{k-1}    &F^{(m)}_{k} \\
F^{(m)}_{k-m+1} &F^{(m)}_{k-m+2}    &\cdots&F^{(m)}_{k-3}  &F^{(m)}_{k-2}   &F^{(m)}_{k-1}\\
F^{(m)}_{k-m} &F^{(m)}_{k-m+1}     &\cdots&F^{(m)}_{k-4} &F^{(m)}_{k-3}  &F^{(m)}_{k-2}\\
 \end{array}
\right)\\
\end{eqnarray*}
Expand the leftmost column of the first determinant in the resulting
expression using the relation $F^{(m)}_{n} = \sum^{m}_{i=1}
F^{(m)}_{n-i}$, then multiply all right columns by -1 and add them
to the left column. By properties of determinants, we obtain
\begin{eqnarray*}
det\textbf{H}^{(m_{-1})}_{k+1}&=&(-1)^{m-2}2det\left(
\begin{array}{ccccccccccccccccc}
\sum^{m}_{i=1}F^{(m)}_{k+m-2-i} &F^{(m)}_{k-1}   &\cdots&F^{(m)}_{k+m-6}&F^{(m)}_{k+m-5}  &F^{(m)}_{k+m-4} \\
        \vdots   &\vdots            &\vdots &      \vdots      &\vdots        &\vdots\\
\sum^{m}_{i=1}F^{(m)}_{k+2-i} &F^{(m)}_{k-m+3}    &\cdots &F^{(m)}_{k-2} &F^{(m)}_{k-1}    &F^{(m)}_{k} \\
\sum^{m}_{i=1}F^{(m)}_{k+1-i} &F^{(m)}_{k-m+2}    &\cdots&F^{(m)}_{k-3}  &F^{(m)}_{k-2}   &F^{(m)}_{k-1}\\
\sum^{m}_{i=1}F^{(m)}_{k-i} &F^{(m)}_{k-m+1}     &\cdots&F^{(m)}_{k-4} &F^{(m)}_{k-3}  &F^{(m)}_{k-2}\\
 \end{array}
\right)\\
\end{eqnarray*}
\begin{eqnarray*}
&&+(-1)^{m-1}det\left(
\begin{array}{ccccccccccccccccc}
F^{(m)}_{k-2} &F^{(m)}_{k-1}   &\cdots&F^{(m)}_{k+m-6}&F^{(m)}_{k+m-5}  &F^{(m)}_{k+m-4} \\
        \vdots   &\vdots            &\vdots &      \vdots      &\vdots        &\vdots\\
F^{(m)}_{k-m+2} &F^{(m)}_{k-m+3}    &\cdots &F^{(m)}_{k-2} &F^{(m)}_{k-1}    &F^{(m)}_{k} \\
F^{(m)}_{k-m+1} &F^{(m)}_{k-m+2}    &\cdots&F^{(m)}_{k-3}  &F^{(m)}_{k-2}   &F^{(m)}_{k-1}\\
F^{(m)}_{k-m} &F^{(m)}_{k-m+1}     &\cdots&F^{(m)}_{k-4} &F^{(m)}_{k-3}  &F^{(m)}_{k-2}\\
 \end{array}
\right)\\
\end{eqnarray*}
\begin{eqnarray*}
&=&(-1)^{m-2}2det\left(
\begin{array}{ccccccccccccccccc}
F^{(m)}_{k+m-3}+F^{(m)}_{k-2} &F^{(m)}_{k-1}   &\cdots&F^{(m)}_{k+m-6}&F^{(m)}_{k+m-5}  &F^{(m)}_{k+m-4} \\
        \vdots   &\vdots            &\vdots &      \vdots      &\vdots        &\vdots\\
F^{(m)}_{k+1}+F^{(m)}_{k-m+2} &F^{(m)}_{k-m+3}    &\cdots &F^{(m)}_{k-2} &F^{(m)}_{k-1}    &F^{(m)}_{k} \\
F^{(m)}_{k}+F^{(m)}_{k-m+1} &F^{(m)}_{k-m+2}    &\cdots&F^{(m)}_{k-3}  &F^{(m)}_{k-2}   &F^{(m)}_{k-1}\\
F^{(m)}_{k-1}+F^{(m)}_{k-m} &F^{(m)}_{k-m+1}     &\cdots&F^{(m)}_{k-4} &F^{(m)}_{k-3}  &F^{(m)}_{k-2}\\
 \end{array}
\right)\\
\end{eqnarray*}
\begin{eqnarray*}
&&+(-1)^{m-1}det\left(
\begin{array}{ccccccccccccccccc}
F^{(m)}_{k-2} &F^{(m)}_{k-1}   &\cdots&F^{(m)}_{k+m-6}&F^{(m)}_{k+m-5}  &F^{(m)}_{k+m-4} \\
        \vdots   &\vdots            &\vdots &      \vdots      &\vdots        &\vdots\\
F^{(m)}_{k-m+2} &F^{(m)}_{k-m+3}    &\cdots &F^{(m)}_{k-2} &F^{(m)}_{k-1}    &F^{(m)}_{k} \\
F^{(m)}_{k-m+1} &F^{(m)}_{k-m+2}    &\cdots&F^{(m)}_{k-3}  &F^{(m)}_{k-2}   &F^{(m)}_{k-1}\\
F^{(m)}_{k-m} &F^{(m)}_{k-m+1}     &\cdots&F^{(m)}_{k-4} &F^{(m)}_{k-3}  &F^{(m)}_{k-2}\\
 \end{array}
\right)\\
\end{eqnarray*}
\begin{eqnarray*}
&=&(-1)^{m-2}2det\left(
\begin{array}{ccccccccccccccccc}
F^{(m)}_{k+m-3}&F^{(m)}_{k-1}   &\cdots&F^{(m)}_{k+m-6}&F^{(m)}_{k+m-5}  &F^{(m)}_{k+m-4} \\
        \vdots   &\vdots            &\vdots &      \vdots      &\vdots        &\vdots\\
F^{(m)}_{k+1} &F^{(m)}_{k-m+3}    &\cdots &F^{(m)}_{k-2} &F^{(m)}_{k-1}    &F^{(m)}_{k} \\
F^{(m)}_{k} &F^{(m)}_{k-m+2}    &\cdots&F^{(m)}_{k-3}  &F^{(m)}_{k-2}   &F^{(m)}_{k-1}\\
F^{(m)}_{k-1} &F^{(m)}_{k-m+1}     &\cdots&F^{(m)}_{k-4} &F^{(m)}_{k-3}  &F^{(m)}_{k-2}\\
 \end{array}
\right)\\
\end{eqnarray*}
\begin{eqnarray*}
&&+(-1)^{m-2}det\left(
\begin{array}{ccccccccccccccccc}
F^{(m)}_{k-2} &F^{(m)}_{k-1}   &\cdots&F^{(m)}_{k+m-6}&F^{(m)}_{k+m-5}  &F^{(m)}_{k+m-4} \\
        \vdots   &\vdots            &\vdots &      \vdots      &\vdots        &\vdots\\
F^{(m)}_{k-m+2} &F^{(m)}_{k-m+3}    &\cdots &F^{(m)}_{k-2} &F^{(m)}_{k-1}    &F^{(m)}_{k} \\
F^{(m)}_{k-m+1} &F^{(m)}_{k-m+2}    &\cdots&F^{(m)}_{k-3}  &F^{(m)}_{k-2}   &F^{(m)}_{k-1}\\
F^{(m)}_{k-m} &F^{(m)}_{k-m+1}     &\cdots&F^{(m)}_{k-4} &F^{(m)}_{k-3}  &F^{(m)}_{k-2}\\
 \end{array}
\right).\\
\end{eqnarray*}
By the induction hypothesis and the relation $F^{(m)}_{n} = 2F_{n+m}
- F_{n+m-1}$, we conclude \begin{eqnarray*}
det\textbf{H}^{(m_{-1})}_{k+1}&=& (-1)^{2(m-2)}2det\left(
\begin{array}{ccccccccccccccccc}
F^{(m)}_{k-1}   &\cdots&F^{(m)}_{k+m-6}&F^{(m)}_{k+m-5}  &F^{(m)}_{k+m-4}& F^{(m)}_{k+m-3}\\
        \vdots   &\vdots            &\vdots &      \vdots      &\vdots        &\vdots\\
F^{(m)}_{k-m+3}    &\cdots &F^{(m)}_{k-2} &F^{(m)}_{k-1}    &F^{(m)}_{k} &F^{(m)}_{k+1} \\
F^{(m)}_{k-m+2}    &\cdots&F^{(m)}_{k-3}  &F^{(m)}_{k-2}   &F^{(m)}_{k-1} &F^{(m)}_{k} \\
F^{(m)}_{k-m+1}     &\cdots&F^{(m)}_{k-4} &F^{(m)}_{k-3}  &F^{(m)}_{k-2} &F^{(m)}_{k-1} \\
 \end{array}
\right)\\
&&+(-1)^{m-2}(-1)^{(k-m-1)(m-1)}F_{-(k-m)}\\
&=&(-1)^{2(m-2)}(-1)^{(k-m+1-1)(m-1)}2F_{-(k-m+1)}
+(-1)^{(k-m)(m-1)-1}F_{-k+m}\\
&=&(-1)^{(k-m)(m-1)}\Big((-1)^{2(m-2)}2F_{-k+m-1}+(-1)^{-1}F_{-k+m}\Big)
\\
&=&(-1)^{(k+1-1)(m-1)-m(m-1)}\Big(2F_{-k+m-1}-F_{-k+m}\Big)\\
&=&(-1)^{(k+1-1)(m-1)}F_{-(k+1)}.
\end{eqnarray*}
By the principle of mathematical induction, the statement holds for
all positive integers $n \in \mathbb{Z}^{+}$.

Let $k \in \mathbb{Z}$ with $k \geq -1$, and assume the statement
holds for all integers $n$ satisfying $-k \leq n \leq 1$. For $n =
-(k+1)$, we proceed analogously to the proof for $n=k+1$:

Starting from the leftmost column of the determinant and
successively applying the recurrence relations:
$$F^{(m)}_{n}=2F^{(m)}_{n+m}-F^{(m)}_{n+m+1},$$
$$F^{(m)}_{n}=F^{(m)}_{n+m}-F^{(m)}_{n+m-1}-\cdots-F^{(m)}_{n+2}-F^{(m)}_{n+1},$$
$$F^{(m)}_{n}=2F^{(m)}_{n-1}-F^{(m)}_{n-m-1},$$
we can verify that the statement holds for $n = -(k+1)$.

By the principle of mathematical induction, the statement holds for
all integers $n \in \mathbb{Z}$.

Using the properties of dimensional-reducing submatrices from
Theorem 4.6, we can derive various relations for generalized
Fibonacci sequences. For instance:

For the 3-generalized Fibonacci sequence, we obtain a Cassini-type
formula analogous to the classical Fibonacci case:
\begin{eqnarray*}(F^{(3)}_{n})^{2}-F^{(3)}_{n-1}F^{(3)}_{n+1}=det\left(
\begin{array}{ccccccccccccccccc}
F^{(3)}_{n}   &F^{(3)}_{n+1}\\
F^{(3)}_{n-1}  &F^{(3)}_{n}\\
 \end{array}
\right)=detF^{(3_{-1})}_{n-1}=(-1)^{(n-2)(3-1)}F^{(3)}_{-n+1}=F^{(3)}_{-n+1}.\\
\end{eqnarray*}

As another example, consider the dimension-reduced 1-submatrix for
 the Fibonacci matrix $\mathcal{H}^{(3)}$ of the
3-dimensional harmonic-type golden section $T^{(3)}$:
\begin{eqnarray*}
\textbf{H}^{(3_{-1})}=\left(
\begin{array}{ccccccccccccccccc}
F^{(3)}_{n-1}      &F^{(3)}_{n+4}\\
F^{(3)}_{n-2}     &F^{(3)}_{n+3}\\
 \end{array}
\right),\\
\end{eqnarray*}
then its determinant is:
$$det\textbf{H}^{(3_{-1})}=F^{(3)}_{n-1}F^{(3)}_{n+3}-F^{(3)}_{n-2}F^{(4)}_{n+4}.$$
Furthermore, by Theorem 4.6 we have \begin{eqnarray*}
det\textbf{H}^{(3_{-1})}&=&det\left(
\begin{array}{ccccccccccccccccc}
2F^{(3)}_{n+2}-F^{(3)}_{n+3}      &F^{(3)}_{n+4}\\
2F^{(3)}_{n+1}-F^{(3)}_{n+2}     &F^{(3)}_{n+3}\\
 \end{array}
\right)=2det\left(
\begin{array}{ccccccccccccccccc}
F^{(3)}_{n+2}     &F^{(3)}_{n+4}\\
F^{(3)}_{n+1}     &F^{(3)}_{n+3}\\
 \end{array}
\right)-det\left(
\begin{array}{ccccccccccccccccc}
F^{(3)}_{n+3}      &F^{(4)}_{n+4}\\
F^{(3)}_{n+2}     &F^{(3)}_{n+3}\\
 \end{array}
\right)\\&=&2det\left(
\begin{array}{ccccccccccccccccc}
F^{(3)}_{n+2}     &F^{(3)}_{n+3}+F^{(3)}_{n+2}+F^{(3)}_{n+1}\\
F^{(3)}_{n+1}     &F^{(3)}_{n+2}+F^{(3)}_{n+1}+F^{(3)}_{n}\\
 \end{array}
\right)-F^{(3_{-1})}_{n+2}\\
&=&2det\left(
\begin{array}{ccccccccccccccccc}
F^{(3)}_{n+2}     &F^{(3)}_{n+3}\\
F^{(3)}_{n+1}     &F^{(3)}_{n+2}\\
 \end{array}
\right)-2det\left(
\begin{array}{ccccccccccccccccc}
F^{(3)}_{n+1}     &F^{(3)}_{n+2}\\
F^{(3)}_{n}     &F^{(3)}_{n+1}\\
 \end{array}
\right)-F^{(3)}_{-n-2}\\
&=&2F^{(3)}_{-n-1}-2F^{(3)}_{-n}-F^{(3)}_{-n-2}=4F^{(3)}_{-n-1}-F^{(3)}_{-n+2},
\end{eqnarray*} leading to the identity
$$F^{(3)}_{n-1}F^{(3)}_{n+3}-F^{(3)}_{n-2}F^{(4)}_{n+4}=4F^{(3)}_{-n-1}-F^{(3)}_{-n+2}.$$

\subsection*{4.4.2 \ Dimensional-augmented Submatrices}

\vspace{1.5mm}\noindent{\bf Theorem 4.7} ~ {\it For
$m\in\mathbb{Z}^{+}$ with $m \geq 2$,} $det\textbf{H}^{(m_{+})}=0.$

\vspace{1.5mm}{\bf Proof} ~ The dimensional-augmented matrix
$\textbf{H}^{(m_{+})}$ possesses a minimum of $m+1$ column vectors.
Based on the fundamental characteristics of $m$-generalized
Fibonacci sequences, these column vectors exhibit linear dependence.
It therefore follows that: $det\textbf{H}^{(m_{+})}=0.$

This theorem provides an effective approach for addressing
lower-dimensional Fibonacci problems through higher-order
determinant analysis. Specific applications include: The proofs of
equations $7^\circ$, $9^\circ$, and $14^\circ$ in Section 4.2.1, as
well as equation (18) in Section 4.1, can be straightforwardly
derived using Theorem 4.7.
\begin{eqnarray*}
\begin{array}{ll}
7^{\circ}~F_{n+4}F_{n-2}-F_{n+3}F_{n-1}=(-1)^{n-1}5 &\hspace{2cm} 9^{\circ}~F_{n}F_{n+4}-F_{n+1}F_{n+3}=2(-1)^{n-1}\\
14^{\circ}~F^{2}_{n+1}-F_{n+4}F_{n-2}=(-1)^{n}4 &\hspace{2cm} F^{(2)}_{n}F^{(2)}_{n+5}+2F^{(2)}_{n+1}
F^{(2)}_{n+4}=3F^{(2)}_{n+2}F^{(2)}_{n+3}~~~  (18)\\
\end{array}
\end{eqnarray*}

For equation $7^\circ$, consider the dimensional-augmented
1-submatrix of the Fibonacci matrix $\mathcal{H}^{(2)}$ associated
with the 2-dimensional harmonic-type golden section $T^{(2)}$
\begin{eqnarray*}
\textbf{H}^{(2_{+1})}=\left(
\begin{array}{ccccccccccccccccc}
F^{(2)}_{2} &F^{(2)}_{n}        &F^{(2)}_{n+5} \\
F^{(2)}_{1} &F^{(2)}_{n-1}      &F^{(2)}_{n+4}\\
F^{(2)}_{0} &F^{(2)}_{n-2}     &F^{(2)}_{n+3}\\
 \end{array}
\right).\\
\end{eqnarray*}
By Theorem 4.7, we have: $\det\textbf{H}^{(2_{+1})} = 0.$ Expanding
$\det\textbf{H}^{(2_{+1})}$
  along its leftmost column yields
  \begin{eqnarray*}
det\textbf{H}^{(2_{+1})}=F^{(2)}_{n-1}F^{(2)}_{n+3}-F^{(2)}_{n-2}F^{(2)}_{n+4}-\left|
\begin{array}{ccccccccccccccccc}
F^{(2)}_{n}        &F^{(2)}_{n+5} \\
F^{(2)}_{n-2}     &F^{(2)}_{n+3}\\
 \end{array}
\right|=F^{(2)}_{n-1}F^{(2)}_{n+3}-F^{(2)}_{n-2}F^{(2)}_{n+4}-(-1)^{n-2}F^{(2)}_{5}F^{(2)}_{2}.
\end{eqnarray*}
Thus $$F_{n+4}F_{n-2}-F_{n+3}F_{n-1}=(-1)^{n-1}5.$$

For equation (18), consider the dimensionally-augmented 1-submatrix
of the Fibonacci matrix $\mathcal{H}^{(2)}$
\begin{eqnarray*}
\textbf{H}^{(2_{+1})}=\left(
\begin{array}{ccccccccccccccccc}
F^{(2)}_{0} &F^{(2)}_{n+2}        &F^{(2)}_{n+5} \\
F^{(2)}_{-1} &F^{(2)}_{n+1}      &F^{(2)}_{n+4}\\
F^{(2)}_{-2} &F^{(2)}_{n}       &F^{(2)}_{n+3}\\
 \end{array}
\right).\\
\end{eqnarray*}
By Theorem 4.7, we immediately obtain
\begin{eqnarray*}
det\textbf{H}^{(2_{+1})}&=&F^{(2)}_{n}F^{(2)}_{n+5}-F^{(2)}_{n+2}F^{(2)}_{n+3}+F^{(2)}_{n+1}F^{(2)}_{n+5}-
F^{(2)}_{n+2}F^{(2)}_{n+4}\\
&=&F^{(2)}_{n}F^{(2)}_{n+5}-F^{(2)}_{n+2}F^{(2)}_{n+3}+F^{(2)}_{n+1}\big(2F^{(2)}_{n+4}-F^{(2)}_{n+2}\big)-
F^{(2)}_{n+2}\big(2F^{(2)}_{n+3}-F^{(2)}_{n+1}\big)\\
&=&F^{(2)}_{n}F^{(2)}_{n+5}-3F^{(2)}_{n+2}F^{(2)}_{n+3}+2F^{(2)}_{n+1}F^{(2)}_{n+4}=0.
\end{eqnarray*}
Thus
$$F^{(2)}_{n}F^{(2)}_{n+5}+2F^{(2)}_{n+1}F^{(2)}_{n+4}=3F^{(2)}_{n+2}F^{(2)}_{n+3}.$$

\section{ Multiple Online Leasing Problem}

\subsection*{5.1 Problem Statement and Competitive Analysis Model}

A skiing interest group consisting of multiple people needs $m$
pairs of skis to practice skiing for an unknown continuous duration
$T$. Whenever the ski group needs to ski, they can choose to: $(i)$
rent all $m$ pairs of skis, $(ii)$ buy all $m$ pairs, $(iii)$
purchase some and rent the remaining pairs, or other possible
combinations... The rental cost per unit time for each pair of skis
is 1 unit of currency (in the rent option, when the rental duration
is not an integer multiple of the unit time, the cost is
proportional to the usage time). The purchase price for each pair of
skis is $s$ units of currency, where $s$ is greater than one. Once
the group buys a pair of skis at a certain time, they no longer need
to pay rental fees for that pair when $m$ pairs of skis are required
for future skiing sessions. The question is: facing an uncertain
skiing duration $T$, what strategy should the skiing group adopt to
minimize their expenses?

We refer to this problem as the $m$-fold online leasing problem, or
$m$-fold online leasing. Online leasing problems with two or more
folds are collectively referred to as {\bf multiple online leasing
problems}, or {\bf multiple online leasing}, abbreviated as $MOLP$.

For the $m$-fold online leasing problem, it is easy to see that the
optimal offline strategy is: when $T < s$, rent $m$ pairs of skis;
and when $T \geq s$, purchase $m$ pairs of skis at the beginning of
the activity. Thus, the optimal offline cost is
$$C_{OPT}(T)=\left\{
               \begin{array}{ll}
                 mT, & T<s; \\
                 ms, & T\geq s.
               \end{array}
             \right.$$

For the $m$-fold online leasing problem, the cost of any online
strategy $ALG$ is greater than the cost of the optimal offline
algorithm. Therefore, the performance of an online algorithm is
measured by comparing its cost to the cost of the optimal offline
strategy for the same skiing time input $T$. The worst-case ratio
between these two costs is called the competitive ratio of the
online algorithm $ALG$, i.e., the competitive ratio of the online
algorithm is
$$\sup\limits_{T>0}\frac{C_{ALG}(T)}{C_{OPT}(T)},$$
where $C_{ALG}(T)$ is the cost incurred by the online strategy $ALG$
for skiing. If there exists an online algorithm $A$ such that
$$\sup\limits_{T>0}\frac{C_{A}(T)}{C_{OPT}(T)}=\inf_{ALG\in\Omega}\sup
_{T>0}\frac{C_{ALG}(T)}{C_{OPT}(T)}=\alpha,
$$
then the online algorithm $A$ is called the optimal competitive
algorithm for this problem, and $\alpha$ is the optimal competitive
ratio for this problem. Here, $\Omega$ represents the set of all
online strategies. Therefore, under the competitive framework, the
goal is to design the optimal online strategy as much as possible,
or otherwise to develop an online algorithm that minimizes its
competitive ratio for online skiing.

\subsection*{5.2 \ Optimal Competitive Algorithm for Multiple
Online Leasing}

For the $m$-fold online leasing problem ($m-MOLP$), where each pair
of skis is rented for 1 unit of currency per unit time and each pair
of skis is purchased for $s$ units of currency, we address the
online characteristic of an uncertain skiing duration $T$ by
designing the following control (risk) and balancing (risk)
strategy.

\vspace{1.5mm}\noindent{\bf Control and Balancing Strategy ($CB$)}

{\bf Control }: At the initial phase, lease all $m$ pairs of skis to
prevent potential losses from premature purchases should the
activity end early. Should the skiing activity continue until time
$s$, procure all $m$ pairs within the interval $(0,s]$ to control
the risk of incurring high rental fees if the activity continues
further.

{\bf Balancing}:~Within the time interval $(0,s]$, select $m$ time
points $0<t_{1}<t_{2}<\ldots<t_{m-1}<t_{m}=s$. As skiing continues
within $(0,s]$, at each time point $t_{i}$, purchase only one pair
of skis and rent $m-i$ pairs, gradually transitioning from renting
to buying (as illustrated in Fig.6). This balances and disperses
the risk of purchasing multiple pairs of skis only for the activity
to end shortly afterward. The selection of each decision time point
$t_{i}$ ensures that the ratio of the online cost $C_{CB}(t_{i})$ to
the offline cost $C_{OPT}(t_{i})$,
$\frac{C_{CB}(t_{i})}{C_{OPT}(t_{i})}$, remains equal each time a
new pair of skis is purchased. This balances the peak ratio of
online to offline costs caused by the offline adversary terminating
the activity immediately after each $t_{i}$ (as shown in Fig.7).

%%%%%%%%%%%%%%%%%%%%%%%%%%%%%%%%%%%%%%%%%%%%%%%%%%%%%%%%图象插入%%%%%%%%%%%%%%%%%%
\begin{figure}[H]
  \centering
  \includegraphics[width=\textwidth, trim=5 210 5 240, clip]{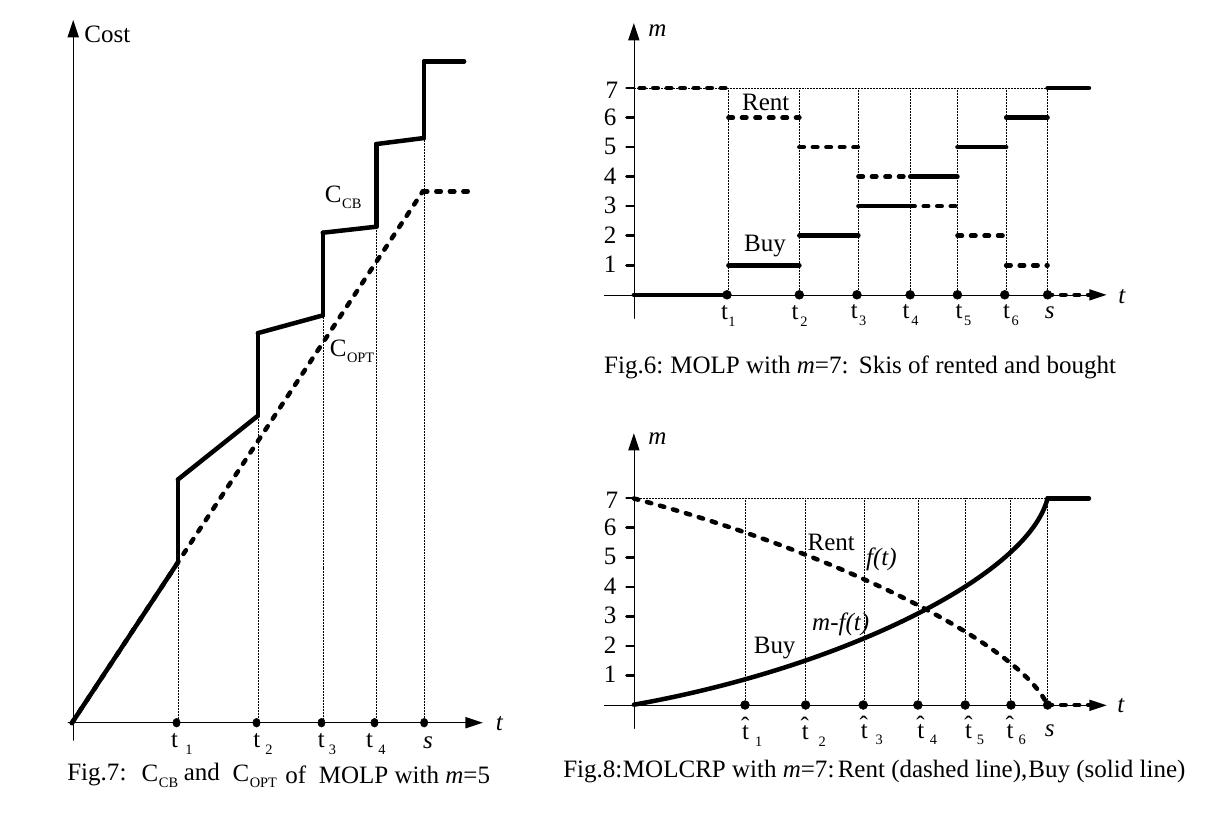}
  \end{figure}
%%%%%%%%%%%%%%%%%%%%%%%%%%%%%%%%%%%%%%%%%%%%%%%%%%%%%%%%%%%%%%%%%%%%%%%%%%%%%%%%%%%

\vspace{1.5mm}\noindent{\bf Theorem 5.1}~{\it For the m-fold online
leasing problem, the $CB$ strategy is a competitive algorithm, and
its competitive ratio $\alpha$ is the unique solution of the
following equation on $( e/(e-1) ,2]$
\begin{eqnarray}
\frac{1}{m(\alpha-1)}+\frac{1}{m(\alpha-1)+1}+\cdots+\frac{1}{m(\alpha-1)+m-1}=1.
\end{eqnarray}}

\vspace{1.5mm}{\bf Proof}~Let $t_{1}, t_{2}, \ldots, t_{m-1}$ be the
moments at which the $CB$ algorithm sequentially purchases a pair of
skis within the time interval $(0,s]$, and let $\alpha$ be the ratio
of online to offline costs when the competitive adversary
sequentially sets the skiing duration $T=t_{1}, t_{2}, \ldots,
t_{m-1}, t_{m}$. According to the $CB$ algorithm, we have
\begin{eqnarray}
\alpha&=&\frac{mt_{1}+s}{mt_{1}}\nonumber\\
&=&\frac{mt_{1}+s+(m-1)(t_{2}-t_{1})+s}{mt_{2}}\nonumber\\
&=&\ldots\nonumber\\
&=&\frac{mt_{1}+s+(m-1)(t_{2}-t_{1})+s+\cdots+2(t_{m-1}-t_{m-2})+s}{mt_{m-1}}\nonumber\\
&=&\frac{mt_{1}+s+(m-1)(t_{2}-t_{1})+s+\cdots+(t_{m}-t_{m-1})+s}{mt_{m}},
\end{eqnarray}
from $\alpha=\frac{mt_{1}+s}{mt_{1}}$, we have $\alpha>1$. Now,
assuming that the skiing group's skiing termination time $T>0,$
according to the $CB$ strategy, we have the following cases:

$(i)$ If $T = t_{i}$ for all $i = 1, 2, \ldots, m$, then the total
cost of renting and purchasing skis for the skiing group is:
$$mt_{1}+s+(m-1)(t_{2}-t_{1})+s+\cdots+(m-i+1)(t_{i}-t_{i-1})+s,$$
while the optimal offline cost is $mt_{i}$. The cost of the $CB$
strategy is $\alpha$ times the cost of the optimal offline
algorithm.

$(ii)$ If $T > s$, then from time $t_{m}$ onward, the skiing group
uses the $m$ pairs of skis purchased in the first $m$ transactions.
The total cost of renting and purchasing skis from time $t_{0} = 0$
until the end of skiing at $T$ is:
$$mt_{1}+s+(m-1)(t_{2}-t_{1})+s+\cdots+(t_{m}-t_{m-1})+s,$$
while the optimal offline cost is $ms$. Thus, the ratio of the cost
of the online algorithm $CB$ to that of the optimal offline
algorithm is $\alpha$.

$(iii)$ If $t_{i-1} < T < t_{i}$ for any $i \in \{1, 2, \ldots,
m\}$, then the total cost of renting and purchasing skis
is:$$mt_{1}+s+(m-1)(t_{2}-t_{1})+s+\cdots+(m-i+2)(t_{i-1}-t_{i-2})+s+(m-i+1)(T-t_{i-1}),$$
while the optimal offline cost is $mT=mt_{i-1}+m(T-t_{i-1})$.
 Since:
\begin{eqnarray*}
\frac{mt_{1}+s+(m-1)(t_{2}-t_{1})+s+\cdots+(m-i+2)(t_{i-1}-t_{i-2})+s}{mt_{i-1}}=\alpha>1
\end{eqnarray*}
and
\begin{eqnarray*}
\frac{(m-i+1)(T-t_{i-1})}{m(T-t_{i-1})}\leq 1,
\end{eqnarray*}
the ratio of the cost of the $CB$ strategy to that of the optimal
offline algorithm satisfies:
\begin{eqnarray*}
\frac{mt_{1}+s+(m-1)(t_{2}-t_{1})+s+\cdots+(m-i+2)(t_{i-1}-t_{i-2})+s(m-i+1)(T-t_{i-1})}{mT}<\alpha.
\end{eqnarray*}

Based on the above analysis, it can be concluded that the worst-case
scenario for the $CB$ strategy is when $T=t_{i} \ \hbox{for all} \
i=1, 2, \ldots, m,$ or $T>s$. Therefore, the competitive ratio of
the $CB$ strategy is
$$\sup\limits_{T}\{\frac{C_{BC}(T)}{C_{OPT}(T)}\}=\alpha.$$
In (24), let $i=1$. From
\begin{eqnarray*}
\alpha=\frac{mt_{1}+s}{mt_{1}},
\end{eqnarray*}
we have
\begin{eqnarray*} t_{1}=\frac{1}{m(\alpha-1)}s.
\end{eqnarray*}
For general $1 \leq i \leq m$, using the iterative method in (24),
we obtain
\begin{eqnarray*}
t_{i}=\Big[\frac{1}{m(\alpha-1)}+\frac{1}{m(\alpha-1)+1}+\cdots+\frac{1}{m(\alpha-1)+i-1}\Big]s.
\end{eqnarray*}
Thus, for $i=m$, we have
\begin{eqnarray*}
t_{m}=\Big[\frac{1}{m(\alpha-1)}+\frac{1}{m(\alpha-1)+1}+\cdots+\frac{1}{m(\alpha-1)+m-1}\Big]s.
\end{eqnarray*}
According to the $CB$ strategy, it holds that $t_{m} = s$.
Therefore, we derive the following equation for $\alpha$:
\begin{eqnarray*}
\frac{1}{m(\alpha - 1)} + \frac{1}{m(\alpha - 1) + 1} + \cdots +
\frac{1}{m(\alpha - 1) + m - 1} = 1.
\end{eqnarray*}

If we let $\alpha-1= \varphi$ in equation (23) and then set $x=
m\varphi$, we obtain the $m$-dimensional harmonic-type golden ratio
equation (6).

By Theorem 2.3, it is known that equation (6) has a unique solution
$\omega$ in the interval $\left(\frac{m}{e - 1}, m\right)$. That is,
$\alpha$ is the unique solution to equation (6) for all positive
integers $m \geq 2$ in the interval $\left(\frac{e}{e - 1},
2\right)$. When $m = 1$, there exists a unique solution $\alpha = 2$
in the interval $\left(\frac{e}{e - 1}, 2\right]$. The proof of the
theorem is complete.

\vspace{1.5mm}\noindent{\bf Corollary 5.1} {\it For the $m$-fold
online leasing problem, the competitive ratio of the $CB$ strategy
is
\begin{eqnarray}
\alpha = 1 + \frac{t_{1} + t_{2} + \cdots + t_{i} + \cdots + t_{m-1}
+ s}{ms}.
\end{eqnarray}}

\vspace{1.5mm}\noindent{\bf Theorem 5.2} \ {\it For the $m$-fold
online leasing problem, the Control and Balancing Strategy is the
unique optimal competitive algorithm for this problem.}

\vspace{1.5mm}{\bf Proof.} Let $t_{1}, t_{2}, \ldots, t_{m}$ denote
the $m$ purchase times of the Control and Balancing Strategy for an
$m$-fold online leasing problem. Let $ALG$ be an optimal online
algorithm different from the $CB$ algorithm for the $m$-fold online
leasing problem, and let $0 \leq t'_{1} \leq t'_{2} \leq \ldots \leq
t'_{m}$ represent the $m$ purchase times of algorithm $ALG$. Below,
we shall prove that algorithm $ALG$ can only achieve a competitive
ratio $\alpha'$ that is exactly equal to $\alpha$.

First, we use mathematical induction to prove that $t'_{i} \geq
t_{i}$ for any $i \in \{1, 2, \ldots, m\}$. Otherwise, algorithm
$ALG$ would obtain a competitive ratio strictly greater than
$\alpha$.

When $i=1$, if $t'_{1} < t_{1}$, then
$$
\alpha^{\prime}\geq\frac{mt^{\prime}_{1}+s}{mt'_{1}}=1+\frac{s}{mt'_{1}}>1+\frac{s}{mt_{1}}=
\frac{mt_{1}+s}{mt_{1}}=\alpha.
$$
Thus, for algorithm $ALG$ to achieve a competitive ratio no worse
than that of the $CB$ algorithm, it must hold that $t'_{1} \geq
t_{1}$.

Now, assume that for any $i < k \in \{2, \ldots, m\}$, it holds that
$t'_{i} \geq t_{i}$.

Below, we use proof by contradiction to show that when $i=k$,
$t'_{k}\geq t_{k}$ must also hold. Suppose $t'_{k} < t_{k}$. Then,
we have
\begin{eqnarray*}
\alpha^{\prime}&\geq&
\frac{mt'_{1}+s+(m-1)(t'_{2}-t'_{1})+s+\cdots+(m-k+1)(t'_{k}-t'_{k-1})+s}{mt'_{k}}\\
&=&\frac{1}{mt'_{k}}\Big[mt'_{1}+m(t'_{2}-t'_{1})+\cdots+m(t'_{k-1}-t'_{k-2})-(t'_{2}-t'_{1})-2(t'_{3}-t'_{2})
-\cdots\\
&~&-(k-2)(t'_{k-1}-t'_{k-2})-mt'_{k-1}+(k-1)t'_{k-1}+(m-k+1)t'_{k}+ks\Big]\\
&=&\frac{m-i+1}{m}+\frac{t'_{1}+t'_{2}\cdots+t'_{k-1}+ks}{mt'_{k}}\\
&>&\frac{m-i+1}{m}+\frac{t_{1}+t_{2}\cdots+t_{k-1}+ks}{mt_{k}}\\
&=&\frac{mt_{1}+s+(m-1)(t_{2}-t_{1})+s+\cdots+(m-i+1)(t_{k}-t_{k-1})+s}{mt_{k}}\\
&=&\alpha.
\end{eqnarray*}
This contradicts the assumption that $ALG$ is an optimal competitive
algorithm. Therefore, $t'_{k} \geq t_{k}$. This proves that for any
$i \in \{1, 2, \ldots, m\}$, it must hold that $t'_{i} \geq t_{i}$.

In particular, $t'_{m} \geq t_{m} = s$. However, if $t'_{m} > t_{m}
= s$, then
\begin{eqnarray*}
\alpha^{\prime}&\geq&
\frac{mt'_{1}+s+(m-1)(t'_{2}-t'_{1})+s+\cdots+(t'_{m}-t'_{m-1})+s}{ms}\\
&=&\frac{t'_{1}+t'_{2}\cdots+t'_{m-1}+t'_{m}+ms}{ms}\\
&>&\frac{t_{1}+t_{2}\cdots+t_{m-1}+ms+t_{m}}{ms}\\
&=&\frac{mt_{1}+s+(m-1)(t_{2}-t_{1})+s+\cdots+(t_{m}-t_{m-1})+s}{ms}\\
&=&\alpha.
\end{eqnarray*}
Thus, for the optimal algorithm $ALG$ to achieve a competitive ratio
no worse than that of the $CB$ algorithm, it must hold that $t'_{m}
= t_{m}$. In fact, for any $i \in \{1, 2, \ldots, m\}$, it must be
true that $t'_{i} = t_{i}$. If not, suppose there exists some $k_{0}
\in \{1, 2, \ldots, m-1\}$ such that $t'_{k_{0}} > t_{k_{0}}$. Then,
we have
\begin{eqnarray*}
\alpha^{\prime}&\geq&
\frac{mt'_{1}+s+(m-1)(t'_{2}-t'_{1})+s+\cdots+(m-k_{0}+1)(t'_{k_{0}}-t'_{k_{0}-1})+s
+\ldots+(t'_{m}-t'_{m-1})+s}{ms}\\
&=&\frac{t'_{1}+t'_{2}\cdots+t'_{k_{0}}+\ldots+t'_{m}+ms}{ms}\\
&>&\frac{t_{1}+t_{2}\cdots+t_{k_{0}}+\ldots+t_{m}+ms}{ms}\\
&=&\frac{mt_{1}+s+(m-1)(t_{2}-t_{1})+s+\cdots+(t_{m}-t_{m-1})+s}{ms}\\
&=&\alpha.
\end{eqnarray*}
This shows that the optimal algorithm $ALG$ is identical to the $CB$
Strategy. This completes the proof of Theorem 5.2.

From Theorem 2.4 in Section 2, we directly obtain the following
Theorem 5.3.

\vspace{1.5mm}\noindent{\bf Theorem 5.3} {\it For all positive
integers $m$, the competitive ratio $\alpha \in \left(\frac{e}{e-1},
2\right]$ of the Control and Balancing Strategy for the $m$-fold
online leasing problem strictly decreases monotonically as $m$
increases and converges to the limit $\frac{e}{e-1}$.}

\vspace{1.5mm}\noindent{\bf Remark 5.1} In the above, we have
extended the classic online ski leasing problem for a single pair of
skis to the multiple online leasing problem. For the multiple online
leasing problem, we proposed the Control and Balancing Strategy, and
proved that this strategy is the uniquely optimal competitive
strategy for this problem. To determine the optimal competitive
ratio for a multiple online leasing problem and the optimal decision
points for purchasing each additional pair of skis, one needs to
solve Equation (23).

For Equation (23), when $m=1$, there exists a unique solution
$\alpha = 2$ in the interval $\left( \frac{e}{e-1}, 2 \right]$
(which corresponds to the classic online leasing problem for a
single pair of skis and its optimal competitive ratio). When $m=2$,
Equation (23) has a unique solution $\alpha=\frac{5+\sqrt{5}}{4}$ in
$\left( \frac{e}{e-1}, 2 \right)$. The optimal time for the first
purchase of skis is $t_{1} = \frac{\sqrt{5}-1}{2}s$ (this aligns
with the double online leasing problem and its golden section
algorithm in Section 1.3).

In fact, if we let $m(\alpha-1)=x$, then Equation (23) is equivalent
to Equation (6). Therefore, solving for the optimal competitive
ratio $\alpha$ of the multiple online leasing problem reduces to
solving the harmonic-type higher-dimensional golden section equation
(6). Determining the optimal decision points for purchasing each
additional pair of skis corresponds to finding the $m$-dimensional
harmonic-type golden section points (sometimes referred to as golden
times (points)) in the interval $[0, s]$.

Thus, the Control and Balancing Strategy for the multiple online
leasing problem is essentially the harmonic-type higher-dimensional
golden section method. For this reason, in the following discussion,
we will not distinguish between the multiple online leasing problem
and the harmonic-type higher-dimensional golden section method, and
will refer to the Control and Balancing Strategy as the
Harmonic-type Higher-dimensional Golden Section Strategy (or
algorithm), abbreviated as $HHGS$.

\section{ The Continuous Relaxation Problem of Multiple Online
Leasing}

The conclusion of Section 5 shows that to solve the competitive
ratio $\alpha$ of an $m$-fold online leasing problem and the first
$m-1$ purchase times $t_{i}$ for all $i=1,2,\ldots,m-1$, we must
first solve the harmonic-type higher-dimensional golden section
equation (6). However, equation (6) may have no explicit solution
when $m\geq5$, making the $HHGS$ Strategy difficult to implement. To
address this, this section first considers the continuous relaxation
of this problem.

An investor has a total usage of $m$ units of a continuously
divisible asset. The rental cost per unit of asset per unit time is
1 unit of currency, and the purchase price per unit of asset is $s$
units of currency. When needed, the investor can rent or purchase
any amount of the asset at any time, but the total amount of rented
and purchased assets at any given time must be $m$ units. Given the
investor's uncertainty about the usage time $T$, the optimal
investment strategy is considered. This problem is referred to as
the continuous relaxation problem of multiple online leasing,
abbreviated as $MOLCRP$.

Clearly, the optimal offline strategy for $MOLCRP$ is: when $T < s$,
the investor rents $m$ units of the asset in total; when $T \geq s$,
the investor purchases $m$ units of the asset at the outset. Thus,
the total cost of the optimal offline strategy is
$$C_{OPT}=\left\{
                 \begin{array}{ll}
                 mT, & T<s; \\
                 ms, & T\geq s.
                 \end{array}
                 \right.
$$

\vspace{1.5mm}\noindent{\textbf{Control and Smoothing Strategy
(CS)}}

\textbf{Control}: Take any variable $t \in [0, s]$. Suppose that at
time $t$, the amount of rented assets is $f(t)$, then the amount of
purchased assets is $m-f(t)$. At $t=0$, the investor rents the full
$m$ units of the asset to avoid immediate investment failure. By
$t=s$, the total amount of purchased assets should equal the total
required usage of $m$ units, to control the risk that after $t=s$,
the usage time $T$ may continue, and the investor still has
unpurchased assets for which rental fees must be paid.

\textbf{Smoothing}: As the usage time $T$ continues, the amount of
purchased assets $m-f(t)$ is a monotonically increasing function,
while the rented amount $f(t)$ is a monotonically decreasing
function. At any time $t \in [0, s]$, the ratio of the online
investor's cost $C_{CS}(t)$ to the optimal offline cost
$C_{OPT}(t)$, $\frac{C_{CS}(t)}{C_{OPT}(t)}$, remains constant. From
$t=0$ to $t=s$, the transition from renting to purchasing is
achieved uniformly and smoothly, as shown in Fig.8. Here, $f(t)$
is a continuously differentiable function that decreases uniformly
and smoothly with $t$.

\vspace{1.5mm}\noindent{\bf Theorem 6.1} \ {\it For the continuous
relaxation problem of multiple online leasing, the competitive ratio
of the Control and Smoothing Strategy is
\begin{eqnarray*}
\alpha=\frac{e}{e-1}.
\end{eqnarray*}}

\vspace{1.5mm}{\bf Proof} \ When $T\leq s$, for any $0\leq t\leq s$,
according to the $CS$ Strategy, the investor's leased value of the
asset is $f(t)$, and the total amount of the asset purchased by the
investor from the start until time $t$ is $m-f(t)$. When the
investor stops using the asset at time $t$, there exists a constant
$\alpha$ such that
\begin{eqnarray*}
\frac{\int^{t}_{0}f(\tau)d\tau+s\big(m-f(t)\big)}{mt}=\alpha.
\end{eqnarray*}
That is,
\begin{eqnarray*}
\int^{t}_{0}f(\tau)d\tau-sf(t)-m\alpha t+sm=0.
\end{eqnarray*}
Differentiating both sides with respect to $t$ yields
\begin{eqnarray*}
f^{\prime}(t)-\frac{1}{s}f(t)+\frac{1}{s}m\alpha =0.
\end{eqnarray*}
The general solution to this differential equation is
\begin{eqnarray*}
f(t)=m\alpha+Ce^{\frac{t}{s}}.
\end{eqnarray*}
The solution to the initial value problem satisfying the initial
condition $f(t)\mid_{t=0}=m$ is
\begin{eqnarray}
f(t)=m\alpha+m(1-\alpha)e^{\frac{t}{s}}.
\end{eqnarray}
By the $CS$ Strategy, we have
\begin{eqnarray}
f(t)\mid_{t=s}=0.
\end{eqnarray}
Substituting (27) into (26) yields
\begin{eqnarray}
\alpha=\frac{e}{e-1}.
\end{eqnarray}
Substituting (28) into (26) gives
\begin{eqnarray*}
f(t)=\frac{m}{e-1}(e-e^{\frac{t}{s}}).
\end{eqnarray*}

When $T>s$, since the $CS$ Strategy has completed the purchase of
all $m$ units of the required asset by time $t=s$, at any time $t
(s<t\leq T)$, the investor no longer needs to lease or purchase the
asset. Therefore, for $t>s$, the total cost of the $CS$ Strategy
equals the total cost at time $t=s$, i.e.,
\begin{eqnarray*}
C_{t>s}&=&\int^{s}_{0}f(\tau)d\tau+s\big(m-f(s)\big)\\
&=&\int^{s}_{0}\frac{m}{e-1}(e-e^{\frac{\tau}{s}})d\tau+s\big(m-0\big)\\
&=&ms\frac{e}{e-1}.
\end{eqnarray*}
Since for $T>s$, the optimal offline cost is $ms$, the ratio of the
$CS$ Strategy's cost $C_{t>s}$ to the optimal offline cost $ms$
remains $\frac{e}{e-1}$.

Combining the cases $T\leq s$ and $T>s$, the competitive ratio of
the Control and Smoothing Strategy for this problem is
$\alpha=\frac{e}{e-1}$. The proof of Theorem 6.1 is complete.

\vspace{1.5mm}\noindent{\bf Corollary 6.1} \ {\it For the continuous
relaxation problem of multiple online leasing, we have:

$(i)$~ The leased value of the Control and Smoothing Strategy at any
time $t$ is
\begin{eqnarray} f(t)=\left\{
       \begin{array}{ll}
         \frac{m}{e-1}(e-e^{\frac{t}{s}}), & t\leq s; \\
         0, & t>s.
       \end{array}
     \right.
\end{eqnarray}

$(ii)$ Within the time interval $[0,t]$, the cumulative purchased
value of the asset under the Control and Smoothing Strategy is
\begin{eqnarray*}
\left\{
  \begin{array}{ll}
    \frac{e^{\frac{t}{s}}-1}{e-1}m, & t\leq s; \\
    m, & t>s.
  \end{array}
\right.
\end{eqnarray*}}

Using a proof method similar to that of Theorem 5.2, we obtain the
following Theorem 6.2.

\vspace{1.5mm}\noindent{\bf Theorem 6.2} {\it For the continuous
relaxation problem of multiple online leasing, the competitive ratio
$e/(e-1)$ of the Control and Smoothing Strategy is the minimal
achievable competitive ratio for this problem. Moreover, the $CS$
Strategy is the optimal strategy for this problem.}

\section{Relaxed Solution and Approximation Analysis for
Multiple Online Leasing}

In this section, we utilize the optimal rental formula (29) of the
continuous relaxation problem for multiple online leasing ($MOLCRP$)
to first derive the approximate purchase times $t'_1, t'_2, t'_3,
\cdots, t'_{m-1}, t'_m$ for multiple online leasing. These are
referred to as the relaxed times. Subsequently, we obtain an
approximate formula for the competitive ratio of multiple online
leasing, termed the relaxed competitive ratio. Additionally, by
establishing the range and trends of the relaxed competitive ratio,
we conduct numerical analysis for cases where $m$ takes different
values. This demonstrates the feasibility of approximating the
optimal competitive ratio with the relaxed competitive ratio and
discusses the resulting errors.

For multiple online leasing, the $CB$ strategy involves purchasing
each pair of skis at the critical times $t_i (i = 1, 2, \cdots, m),$
which are the golden times (as shown by the left endpoints of the
solid lines in Figure 6). When the skiing time ends at any $t_i$,
the ratio of the online cost to the offline cost equals $\alpha$ (as
illustrated in Figure 7). In contrast, the $CS$ strategy of $MOLCRP$
involves the continuous and gradual accumulation of divisible assets
according to the continuous function $m - f(t)$ (as shown by the
solid line in Figure 8). Here, we delay the critical times $t_i$ of
the $CB$ strategy until the times $t'_i$ when the $CS$ strategy of
$MOLCRP$ accumulates $i$ units of the asset (as indicated by the
solid points on the $t$-axis in Figure 8). That is, $m - f(t) = i$,
or equivalently, $f(t) = m - i ~(i = 1, 2, \cdots, m).$ Therefore,
in (29), we sequentially set
$$f(t) = m - 1, \ f(t) = m - 2, \ \cdots, \ f(t) = m - i, \
\cdots, \ f(t) = m - m,$$  and solve for $$t = s \ln \frac{m + e -
1}{m}, \, t = s \ln \frac{m + 2(e - 1)}{m}, \, \cdots, \, t = s \ln
\frac{m + i(e - 1)}{m}, \, \cdots, \, t = s \ln \frac{m + m(e -
1)}{m} = s.$$ These are then denoted sequentially as
$$t'_{1}=s\ln\frac{m+e-1}{m},~t'_{2}=s\ln\frac{m+2(e-1)}{m},\cdots,
~t'_{i}=s\ln\frac{m+i(e-1)}{m},\cdots,~t'_{m}=s\ln\frac{m+m(e-1)}{m}=s.$$
Using $t'_1, t'_2, \cdots, t'_m$ as the approximate times for
purchasing each pair of skis in $m$ stages---referred to as the
relaxed times for multiple online leasing---we define the relaxed
competitive ratio as the ratio of online to offline costs when the
offline adversary chooses to stop skiing at $t'_m = 1$:
\begin{eqnarray*}
\alpha' &=& \frac{m t'_1 + s + (m - 1)(t'_2 - t'_1) + s + \cdots +
(m - i + 1)(t'i - t'{i-1}) + s + \cdots + (t'm - t'{m-1}) + s}{m
t'_m} \\ &=& \frac{t'_1 + t'2 + \cdots + t'{m-1} + t'_m + m s}{m
t'_m}.
\end{eqnarray*}
Substituting each $t'_i$
  into the above expression yields the relaxed competitive ratio:
\begin{eqnarray*}
\alpha' = 1 + \frac{\sum^{m}_{i=1} \ln \left[ 1 + \frac{i(e - 1)}{m}
\right]}{m}.
\end{eqnarray*}

Here, we refer to the relaxed times $t'_{i} (i=1,2,\cdots,m)$ and
the relaxed competitive ratio $\alpha'$ collectively as the relaxed
solution for the $m$-fold online leasing problem, and we call
$x'=m(\alpha'-1)$ the relaxed solution of the $m$-dimensional
harmonic-type golden section equation (6).

For simplicity, in the following analysis and derivations, we set
$s=1$, effectively treating the original interval $[0, s]$ as $[0,
1]$. Under this simplification, the relaxed times and the relaxed
competitive ratio become:
\begin{eqnarray}
t'_i = \ln \frac{m + i(e - 1)}{m}, \quad i = 1, 2, \cdots, m.
\end{eqnarray}
\begin{eqnarray}
\hat{\alpha}=1+\frac{t'_{1}+t'_{2}+\cdots +t'_{m-1}+1}{m}~~~
or~~~1+\frac{\sum^{m}_{i=1}ln[1+\frac{i(e-1)}{m}]}{m}.
\end{eqnarray}

We observe that the relaxed competitive ratio $\alpha'$ corresponds
to the ratio of online to offline costs only when the offline
adversary chooses to stop skiing at $t'_m = 1$. To refine this
competitive ratio, we need to adjust each $t'_i$. Let
\begin{eqnarray}
t''_{i}=\sum^{i}_{j=1}\frac{1}{m(\alpha'-1)+j-1}, ~i=1,2,\cdots,m.
\end{eqnarray}
This yields the corrected purchase times $t''_i$ for each pair of
skis, corresponding to the relaxed competitive ratio $\alpha'$.

\vspace{1.5mm}\noindent{\bf Theorem 7.1}~ {\it Relaxed competitive
ratio $\frac{e}{e-1} < \alpha' \leq 2$ strictly decreases as $m$
increases, converging to $\frac{e}{e-1}$ as the limit.}

\vspace{1.5mm}{\bf Proof}\ To prove Theorem 7.1, we need to show
that the sequence
$\{\frac{\sum^{m}_{i=1}ln[1+\frac{i(e-1)}{m}]}{m}\}$ lies within the
interval $\left( \frac{1}{e-1}, 1 \right]$, strictly decreases with
$m$, and converges to $\frac{1}{e-1}$ as its limit.

Let $\{\frac{\sum^{m}_{i=1}ln[1+\frac{i(e-1)}{m}]}{m}\}$ be denoted
as $\{x_{m}\}$, and rewrite its general term as
\begin{eqnarray*}
 x_{m}&=&\frac{\sum^{m}_{i=1}ln[1+\frac{i(e-1)}{m}]}{m}\\
 &=& \Big\{\frac{m}{m(m+1)}ln[1+\frac{1}{m}(e-1)]+\frac{1}{m(m+1)}ln[1+\frac{1}{m}(e-1)]\Big\}\\
 &&+\Big\{\frac{m-1}{m(m+1)}ln[1+\frac{2}{m}(e-1)]+\frac{2}{m(m+1)}ln[1+\frac{2}{m}(e-1)]\Big\}\\
 &&+\cdots \\
 &&+\Big\{\frac{m-i+1}{m(m+1)}ln[1+\frac{i}{m}(e-1)]+\frac{i}{m(m+1)}ln[1+\frac{i}{m}(e-1)]\Big\}\\
 &&+\cdots \\
 &&+\Big\{\frac{2}{m(m+1)}ln[1+\frac{m-1}{m}(e-1)]+\frac{m-1}{m(m+1)}ln[1+\frac{m-1}{m}(e-1)]\Big\}\\
 &&+\Big\{\frac{1}{m(m+1)}ln[1+\frac{m}{m}(e-1)]+\frac{m}{m(m+1)}ln[1+\frac{m}{m}(e-1)]\Big\}\\
\end{eqnarray*}
For any $1 \leq i \leq m$, we denote successively
\begin{eqnarray*}
a_{2i-1}=\frac{m-i+1}{m(m+1)}ln[1+\frac{i}{m}(e-1)]
\end{eqnarray*}
\begin{eqnarray*}
\hbox{and} \hspace{5cm}
a_{2i}=\frac{i}{m(m+1)}ln[1+\frac{i}{m}(e-1)],\hspace{5cm}
\end{eqnarray*}
\begin{eqnarray*}
\hbox{then}
\hspace{5cm}x_{m}=\frac{\sum^{m}_{i=1}ln[1+\frac{i(e-1)}{m}]}{m}=\sum^{m}_{i=1}
\{a_{2i-1}+ a_{2i}\}\hspace{5cm}
\end{eqnarray*}
Rewriting $x_{m+1}$ as
\begin{eqnarray*}
 x_{m+1}&=&\frac{\sum^{m+1}_{i=1}ln[1+\frac{i(e-1)}{m+1}]}{m+1}\\
 &=& \frac{m}{m(m+1)}ln[1+\frac{1}{m+1}(e-1)]\\
 &&+\Big\{\frac{1}{m(m+1)}ln[1+\frac{2}{m+1}(e-1)]+\frac{m-1}{m(m+1)}ln[1+\frac{2}{m+1}(e-1)]\Big\}\\
 &&+\cdots \\
 &&+\Big\{\frac{i-1}{m(m+1)}ln[1+\frac{i}{m+1}(e-1)]+\frac{m-i+1}{m(m+1)}ln[1+\frac{i}{m+1}(e-1)]\Big\}\\
 &&+\Big\{\frac{i}{m(m+1)}ln[1+\frac{i+1}{m+1}(e-1)]+\frac{m-i}{m(m+1)}ln[1+\frac{i+1}{m+1}(e-1)]\Big\}\\
 &&+\cdots \\
 &&+\Big\{\frac{m-1}{m(m+1)}ln[1+\frac{m}{m+1}(e-1)]+\frac{1}{m(m+1)}ln[1+\frac{m}{m+1}(e-1)]\Big\}\\
 &&+\frac{m}{m(m+1)}ln[1+\frac{m+1}{m+1}(e-1)].\\
\end{eqnarray*}
Let
\begin{eqnarray*}
b_{1} = \frac{m}{m(m+1)} \ln \left[ 1 + \frac{1}{m+1}(e-1) \right]
\end{eqnarray*}
and
\begin{eqnarray*}
b_{2m} = \frac{m}{m(m+1)} \ln \left[ 1 + \frac{m+1}{m+1}(e-1)
\right].
\end{eqnarray*}
For any $2 \leq i \leq m$, define successively:
\begin{eqnarray*}
b_{2i-2} = \frac{i-1}{m(m+1)} \ln \left[ 1 + \frac{i}{m+1}(e-1)
\right],
\end{eqnarray*}
and
\begin{eqnarray*}
 b_{2i-1} = \frac{m-i+1}{m(m+1)} \ln \left[ 1 +
\frac{i}{m+1}(e-1) \right].
\end{eqnarray*}
Then,
\begin{eqnarray*}
x_{m+1}=\frac{\sum^{m+1}_{i=1}ln[1+\frac{i(e-1)}{m+1}]}{m+1}=b_{1}+\sum^{m}_{i=2}
\{b_{2i-2}+ b_{2i-1}\}+b_{2m}.
\end{eqnarray*}
Since \begin{eqnarray*}
(a_{1}-b_{1})+(a_{2}-b_{2})&=&\Big(\frac{m}{m(m+1)}ln[1+\frac{1}{m}(e-1)]-\frac{m}{m(m+1)}ln[1+\frac{1}{m+1}(e-1)]\Big)\\
&&+\Big(\frac{1}{m(m+1)}ln[1+\frac{1}{m}(e-1)]-\frac{1}{m(m+1)}ln[1+\frac{2}{m+1}(e-1)]\Big)\\
&=&\frac{m}{m(m+1)}ln\Big(\frac{m+1}{m}\frac{m+e-1}{m+e}\Big)+
\frac{1}{m(m+1)}ln\Big(\frac{m+1}{m}\frac{m+e-1}{m+2e-1}\Big)\\
&=&\frac{1}{m(m+1)}ln\Big((1+\frac{1}{m})^{m+1}\frac{1}{(1+\frac{1}{m+e-1})^{m+e}}
(1+\frac{1}{m+e-1})^{e}\frac{m+e-1}{m+2e-1}\Big)\\
&>&\frac{1}{m(m+1)}ln\Big((1+\frac{1}{m})^{m+1}\frac{1}{(1+\frac{1}{m})^{m+1}}
(1+\frac{1}{m+e-1})^{e}\frac{m+e-1}{m+2e-1}\Big)\\
&>&\frac{1}{m(m+1)}ln\Big((1+\frac{e}{m+e-1})\frac{m+e-1}{m+2e-1}\Big)\\
&=&0
\end{eqnarray*}
and for any $2 \leq i \leq m-1$, we have
\begin{eqnarray*}
\noindent&&(a_{2i-1}-b_{2i-1})+(a_{2i}-b_{2i})\\
&=&\Big(\frac{m-i+1}{m(m+1)}ln[1+\frac{i}{m}(e-1)]-\frac{m-i+1}{m(m+1)}ln[1+\frac{i}{m+1}(e-1)]\Big)\\
&&+\Big(\frac{i}{m(m+1)}ln[1+\frac{i}{m}(e-1)]-\frac{i}{m(m+1)}ln[1+\frac{i+1}{m+1}(e-1)]\Big)\\
&=&\frac{m-i+1}{m(m+1)}ln\Big(\frac{m+1}{m}\frac{m+i(e-1)}{m+1+i(e-1)}\Big)+
\frac{i}{m(m+1)}ln\Big(\frac{m+1}{m}\frac{m+i(e-1)}{m+1+(i+1)(e-1)}\Big)\\
&=&\frac{1}{m(m+1)}ln\Big((1+\frac{1}{m})^{m+1}\frac{1}{(1+\frac{1}{m+i(e-1)})^{m+i(e-1)+1}}
(1+\frac{1}{m+i(e-1)})^{ie}(\frac{m+i(e-1)}{m+1+(i+1)(e-1)})^{i}\Big)\\
&>&\frac{1}{m(m+1)}ln\Big((1+\frac{1}{m})^{m+1}\frac{1}{(1+\frac{1}{m})^{m+1}}
(1+\frac{1}{m+i(e-1)})^{ie}(\frac{m+i(e-1)}{m+1+(i+1)(e-1)})^{i}\Big)\\
&>&\frac{1}{m(m+1)}ln\Big((1+\frac{e}{m+i(e-1)})^{i}(\frac{m+i(e-1)}{m+1+(i+1)(e-1)})^{i}\Big)\\
&=&0,
\end{eqnarray*}
Moreover,
\begin{eqnarray*}
&&(a_{2m-1}-b_{2m-1})+(a_{2m}-b_{2m})\\
&=&\Big(\frac{1}{m(m+1)}ln[1+\frac{m}{m}(e-1)]-\frac{1}{m(m+1)}ln[1+\frac{m}{m+1}(e-1)]\Big)\\
&&+\Big(\frac{m}{m(m+1)}ln[1+\frac{m}{m}(e-1)]-\frac{m}{m(m+1)}ln[1+\frac{m+1}{m+1}(e-1)]\Big)\\
&>&0.
\end{eqnarray*}
Therefore, \begin{eqnarray*}
x_{m}-x_{m+1}=\sum_{i=1}^{m}[(a_{2i-1}-b_{2i-1})+(a_{2i}-b_{2i})]>0.
\end{eqnarray*}
Moreover, since $x_{1} = 1$ and
\begin{eqnarray*}
\lim_{m\rightarrow\infty}x_{m}=\lim_{m\rightarrow\infty}\frac{\sum^{m}_{i=1}ln[1+\frac{i(e-1)}{m}]}{m}
=\int^{1}_{0}ln[1+(e-1)x]dx=\frac{1}{e-1}.
\end{eqnarray*}
Therefore, the sequence
$\{\frac{\sum^{m}_{i=1}ln[1+\frac{i(e-1)}{m}]}{m}\}$ is strictly
decreasing within the interval $\left(\frac{1}{e-1}, 1\right]$ and
converges to $\frac{1}{e-1}$ as its limit.

By comparing Theorem 7.1 and Theorem 5.3, the relaxed competitive
ratio for the multiple online leasing problem remains within the
same range as the optimal competitive ratio and follows the same
trend as $m$ varies. This indicates that approximating the optimal
competitive ratio with the relaxed competitive ratio derived from
Equations (30), (31), and (32) is effective and feasible, especially
when $m$ is large.

$$\footnotesize
\begin{tabular}{c|cccccccccc}
\multicolumn{11}{c}{\bf Table 1: ~$m$-fold Online Leasing Software
Solution: ${t}_{i}$, $\alpha$;
Relaxed Solution: $t'_{i}$, $\alpha'$; Corrected Time $t''_{i}$}\\
\hline
  $m$              &1      &2       &3       &4       &5       &6       &7       &8       &9       &10 \\\hline
  $t_{1}$          &1      &0.61803 &0.45160 &0.35665 &0.29494 &0.25153 &0.21930 &0.19442 &0.17461 &0.15848\\
  $t'_{1}$         &1      &0.62011 &0.45283 &0.35737 &0.29539 &0.25183 &0.21951 &0.19457 &0.17473 &0.15857\\
$t''_{1}$      &1      &0.61724 &0.45122 &0.35645 &0.29482 &0.25146 &0.21925 &0.19438 &0.17459 &0.15846\\
  \hline
  $t_{2}$          &       & 1      &0.76271 &0.61955 &0.52270 &0.45251 &0.39916 &0.35718 &0.32327 &0.29528\\
  $t'_{2}$         &       & 1      &0.76338 &0.62011 &0.52314 &0.45283 &0.39941 &0.35737 &0.32342 &0.29539\\
$t''_{2}$      &       & 0.99890&0.76214 &0.61923 &0.52251 &0.45239 &0.39908 &0.35713 &0.32323 &0.29525\\
  \hline
  $t_{3}$          &       &        & 1      &0.82772 &0.70821 &0.61985 &0.55160 &0.49717 &0.45268 &0.41561\\
  $t'_{3}$         &       &        & 1      &0.82799 &0.70851 &0.62011 &0.55182 &0.49735 &0.45283 &0.41574\\
$t''_{3}$      &       &        & 0.99932&0.82732 &0.70798 &0.61970 &0.55150 &0.49709 &0.45263 &0.41558\\
  \hline
  $t_{4}$          &       &        &        & 1      &0.86469 &0.76321 &0.68387 &0.61997 &0.56727 &0.52303\\
  $t'_{4}$         &       &        &        & 1      &0.86484 &0.76338 &0.68404 &0.62011 &0.56740 &0.52314\\
$t''_{4}$      &       &        &        & 0.99957&0.86443 &0.76304 &0.68376 &0.61987 &0.56721 &0.5298\\
  \hline
  $t_{5}$          &       &        &        &        & 1      &0.88859 &0.80070 &0.72933 &0.67008 &0.62002\\
  $t'_{5}$         &       &        &        &        & 1      &0.88867 &0.80081 &0.72944 &0.67019 &0.62011\\
$t''_{5}$      &       &        &        &        & 0.99971&0.88840 &0.80057 &0.72923 &0.67000 &0.61997\\
  \hline
  $t_{6}$          &       &        &        &        &        & 1      &0.90530 &0.82791 &0.76330 &0.70844\\
  $t'_{6}$         &       &        &        &        &        &1       &0.90536 &0.82799 &0.76338 &0.70851\\
$t''_{6}$      &       &        &        &        &        &0.99980 &0.90516 &0.82780 &0.76322 &0.70838\\
  \hline
  $t_{7}$          &       &        &        &        &        &        & 1      &0.91765 &0.84857 &0.78967\\
  $t'_{7}$         &       &        &        &        &        &        & 1      &0.91769 &0.84863 &0.78973\\
$t''_{7}$      &       &        &        &        &        &        & 0.99985&0.91753 &0.84848 &0.78961\\
  \hline
  $t_{8}$          &       &        &        &        &        &        &        & 1      &0.92715 &0.86481\\
  $t'_{8}$         &       &        &        &        &        &        &        & 1      &0.92718 &0.86484\\
$t''_{8}$      &       &        &        &        &        &        &        & 0.99987&0.92705 &0.86474\\
  \hline
 $t_{9}$           &       &        &        &        &        &        &        &        &1       &0.93469\\
 $t'_{9}$          &       &        &        &        &        &        &        &        &1       &0.93470\\
$t''_{9}$      &       &        &        &        &        &        &        &        &0.99990 &0.93462\\
  \hline
  $t_{10}$         &       &        &        &        &        &        &        &        &        &  1    \\
  $t'_{10}$        &       &        &        &        &        &        &        &        &        &  1    \\
$t''_{10}$     &       &        &        &        &        &        &        &        &        &0.99993\\
  \hline
  $\alpha  $       &2.00000&1.80902 &1.73811 &1.70097 &1.67811 &1.66261 &1.65142 &1.64295 &1.63633 &1.63100\\
  $\alpha'$   &2.00000&1.81006 &1.73874 &1.70137 &1.67838 &1.66280 &1.65156 &1.64307 &1.63642 &1.63107\\
   \hline
 \end{tabular}$$

In Table~1, the values of $t_i$ and $\alpha$ for $m = 1, 2, \dots,
10$ were obtained using a binary search algorithm and computational
software, representing the golden times and optimal competitive
ratios for $m$-fold online leasing, accurate to $10^{-5}$. The
values $t'_i$ and $\alpha'$ are the relaxed times and relaxed
competitive ratios derived from Equations (30) and (31), while
$t''_i$ (the corrected times) are the adjusted purchase times
obtained by applying Equation (32) to the relaxed times $t'_i$,
corresponding to the relaxed competitive ratio $\alpha'$.

We conducted a detailed comparison between the golden times $t_i$
from the $CB$ strategy and the relaxed times $t'_i$ in Table 1.
Except for the final purchase time, where $t_m = t'_m = 1$, the
relaxed times $t'_i$ for the first $m-1$ purchases are all delayed
compared to the golden times $t_i  (i = 1, 2, \dots, m-1).$ In
contrast, the corrected times $t''_i$ are all earlier than the
golden times $t_i$. Whether delayed or advanced, both deviations
disrupt the equilibrium of the $CB$ strategy and increase the
competitive ratio (balance ratio).

From Table 1, it is evident that the errors between $t'_i$, $t''_i$
and $t_i$, as well as between $\alpha'$ and $\alpha,$ are
approximately on the order of $10^{-3}$. Admittedly, this level of
accuracy could be improved. However, this method of using natural
logarithms of $m$ and $i$ to derive relaxed times and relaxed
competitive ratios (or balance ratios) may have potential
applications in solving higher-dimensional (five or more)
harmonic-type golden section equations.

\section{\hspace{-4mm}Mantissa Weighted Method for Solving Harmonic-type
Higher-dimensional Golden Section Equations}

In Section 7, We use the optimal solution of the hiring function
$f(t)$ for continuously divisible assets in the continuous
relaxation problem of multiple online leasing to obtain the relaxed
solution for multiple online leasing or the harmonic-type
higher-dimensional golden section equation. However, approximation
analysis shows that directly substituting the relaxed solution for
the optimal solution does not yield very satisfactory accuracy.
Therefore, in this section, we propose a method that better fits the
roots of harmonic-type higher-dimensional golden section
equations--the {\bf mantissa weighted method} of golden section.

From Section 2, we know that the section intervals $\Delta_{1},
~\Delta_{2}, ~\cdots, ~\Delta_{m}$ of the $m$-dimensional
harmonic-type golden section correspond to the terms on the
left-hand side of the $m$-dimensional harmonic-type golden section
equation:
$$\frac{1}{x}+\frac{1}{x+1}+\frac{1}{x+2}
+\cdots+\frac{1}{x+m-1}=1$$ Specifically,
$$t_{i}-t_{i-1}=\frac{1}{x+i-1},~i=1, 2, \cdots, m,$$
where $t_{0}=0.$

Based on the relaxed time formula (30) from the previous section, we
have the following approximation for $i = 1, 2, \dots, m$:
$$\ln\big[1+\frac{e-1}{m+(i-1)(e-1)}\big]=\frac{1}{x+i-1}.$$
This further leads to the approximate formula:
\begin{eqnarray*}
x=\ln^{-1}\big[1+\frac{e-1}{m+(i-1)(e-1)}\big]-i+1, ~i=1, 2, \cdots,
m
\end{eqnarray*}
which we refer to as the relaxed solution set of the $m$-dimensional
harmonic-type golden section equation (6), or simply the relaxed
solution.

\vspace{1.5mm}\noindent{\bf Theorem 8.1}~ {\it Let $i = 1, 2, \dots,
m$. Then, the $m$ relaxed solutions of the harmonic-type
high-dimensional golden section equation (6),
\begin{eqnarray}
x_{i}=\ln^{-1}\big[1+\frac{e-1}{m+(i-1)(e-1)}\big]-i+1,
\end{eqnarray}
satisfy
$$\frac{m}{e-1}+\frac{1}{2}-\frac{e-1}{4m+2(e-1)}<x_{1}<x_{2}<\cdots<x_{m}<\frac{m}{e-1}+\frac{1}{2}.$$}

{\bf Proof}~ We first prove $x_{1}<x_{2}<\cdots<x_{m}$. Let
$$f(i)=\ln^{-1}\big[1+\frac{e-1}{m+(i-1)(e-1)}\big]-i+1.$$
The derivative of $f(i)$ is
$$f'(i)=\frac{\frac{(e-1)^{2}}{[m+(i-1)(e-1)][m+i(e-1)]}-\ln^{2}[1+\frac{e-1}{m+(i-1)(e-1)}]}
{\big(\ln[m+i(e-1)]-\ln[m+(i-1)(e-1)]\big)^{2}}.$$ Using the
inequality
\begin{eqnarray}
\frac{2}{2x+1}<\ln(1+\frac{1}{x})<\frac{1}{\sqrt{x^{2}+x}},
\end{eqnarray}
we obtain
$$\Big[\frac{2(e-1)}{2m+(2i-1)(e-1)}\Big]^{2}<\ln^{2}[1+\frac{e-1}{m+(i-1)(e-1)}]<
\frac{(e-1)^{2}}{[m+(i-1)(e-1)][m+i(e-1)]}.$$ Thus, $f'(i) > 0$,
meaning $f(i)$ is strictly increasing, and hence
$x_{1}<x_{2}<\cdots<x_{m}$.

Next, we prove that for all $i \in \{1, 2, \dots, m\}$,
$$\frac{m}{e-1}+\frac{1}{2}-\frac{e-1}{4m+2(e-1)}<x_{i}<\frac{m}{e-1}+\frac{1}{2}.$$
From inequality (34), we have
$$L=\frac{\sqrt{1+\frac{e-1}{m+(i-1)(e-1)}}}{\frac{e-1}{m+(i-1)(e-1)}}
<\ln^{-1}\big[1+\frac{e-1}{m+(i-1)(e-1)}\big]<
\frac{\frac{e-1}{m+(i-1)(e-1)}+2}{\frac{2(e-1)}{m+(i-1)(e-1)}}=R.$$
For the upper bound:
$$R=\frac{2m+(2i-1)(e-1)}{2(e-1)}=\frac{m}{(e-1)}+i-\frac{1}{2}$$
and thus$$x_{i}<R-i+1=\frac{m}{e-1}+\frac{1}{2}.$$ For the lower
bound:
\begin{eqnarray*}L&=&\frac{\sqrt{[m+(i-1)(e-1)][m+i(e-1)]}}{e-1}>
\frac{2}{e-1}\frac{1}{\frac{1}{m+(i-1)(e-1)}+\frac{1}{m+i(e-1)}}\\
&=&\frac{1}{e-1}\Big[m+(i-1)(e-1)+
\frac{0.5(e-1)[m+(i-1)(e-1)]}{m+(i-0.5)(e-1)}\Big]\\
&=&\frac{m}{e-1}+i-1+\frac{1}{2}-\frac{0.5(e-1)}{2m+(2i-1)(e-1)}\\
&>&\frac{m}{e-1}+i-1+\frac{1}{2}-\frac{e-1}{4m+2(e-1)}.
\end{eqnarray*}
Simplifying further,
$$x_{i}>L-i+1=\frac{m}{e-1}+\frac{1}{2}-\frac{e-1}{4m+2(e-1)}.$$
This completes the proof.

\vspace{1.5mm}\noindent{\bf Corollary 8.1}~ {\it The unique positive
real root $\omega$ of the $m$-dimensional harmonic-type golden
section equation
$$\frac{1}{x}+\frac{1}{x+1}+\frac{1}{x+2}
+\cdots+\frac{1}{x+m-1}=1$$
satisfies$$\omega\in\big(\frac{m}{e-1}+\frac{1}{2}-\frac{e-1}{4m+2(e-1)},~
\frac{m}{e-1}+\frac{1}{2}\big).$$ The balance ratio $\alpha$
satisfies$$\alpha\in\big(\frac{e}{e-1}+\frac{1}{2m}-\frac{e-1}{2m(2m+e-1)},~
\frac{e}{e-1}+\frac{1}{2m}\big).$$}

\vspace{1.5mm}{\bf Proof}~ Since each relaxed time $t'_{i}$ is
delayed compared to the golden time $t_{i}$~($~i=1, 2, \cdots,
m-1$), the section intervals satisfy $\Delta_{1}<\Delta'_{1}$ and
$\Delta_{m}>\Delta'_{m}$. Thus, $$\frac{1}{x}<\frac{1}{x_{1}},~
\frac{1}{x+m-1}>\frac{1}{x_{m}+m-1},$$ which implies
$x_{1}<x<x_{m}$. The conclusion follows directly from Theorem 8.1.

\vspace{1.5mm} \noindent{\bf Relaxed Solutions Golden-Section
Mantissa Weighted Method}

For a positive integer $m$, we perform the golden section to obtain
the total weight: $w=\frac{\sqrt{5}-1}{2}m$, where the integer part
$\lfloor\frac{\sqrt{5}-1}{2}m\rfloor$ is denoted as $m_{G}$, and the
fractional part (mantissa)
$\frac{\sqrt{5}-1}{2}m-\lfloor\frac{\sqrt{5}-1}{2}m\rfloor$ is
denoted as $g_{m}$. We select the first $m_{G}+1$ terms from the $m$
relaxed solutions $\{x_{1}, x_{2}, \cdots, x_{m}\}$ namely $x_{1},
x_{2}, \cdots, x_{m_{G}}, x_{m_{G}+1}$, and assign weights as
follows: $w_{1}=w_{2}=\cdots=w_{m_{G}}=1, w_{m_{G}+1}=g_{m}.$ Taking
the weighted average, we obtain an approximate solution (section
basis) for the $m$-dimensional harmonic-type golden section
equation:
\begin{eqnarray}
\hat{\omega}=\frac{\sqrt{5}+1}{2m}\Big[\sum_{i=1}^{m_{G}}x_{i}
+g_{m}x_{m_{G}+1}\Big].
\end{eqnarray}
Substituting the expression for $x_i$ from (33),
\begin{eqnarray*}
x_{i}=\ln^{-1}\big[1+\frac{e-1}{m+(i-1)(e-1)}\big]-i+1,
\end{eqnarray*}
into (35), we derive the approximate formula for the section basis:
\begin{eqnarray}
\hat{\omega}=\frac{\sqrt{5}+1}{2m}\bigg [\sum_{i=1}^{m_{G}}
\big(\ln^{-1}[1+\frac{e-1}{m+(i-1)(e-1)}]-i+1\big)
+g_{m}\big(\ln^{-1}[1+\frac{e-1}{m+m_{G}(e-1)}]
-m_{G}\big)\bigg].
\end{eqnarray}
This method is termed the {\bf Relaxed Solutions Golden-Section
Mantissa Weighted Method}, or simply the {\bf Mantissa Weighting
Method}.

After obtaining $\hat{\omega}$ using formula (36), we substitute
$\hat{\omega}$ for $x$ in the following equation to derive the
golden section points for the $m$-dimensional harmonic-type golden
section:$$\hat{t}_{i}=\frac{1}{x}+\frac{1}{x+1}
+\cdots+\frac{1}{x+i-1}~~(i=1, 2, \cdots, m-1),$$ From this, we
further obtain the section intervals and the balance ratio:
$$\hat{\Delta}_{i}~(i=1, 2, \cdots, 20)~~\hbox{and}
~~\hat{\alpha}=\frac{\hat{\omega}}{m}+1.$$

Example: For $m = 20$, applying the golden section yields the total
weight:$w=10(\sqrt{5}-1)$, with integer part $m_{G}=12$, and
mantissa $g_{m}=10\sqrt{5}-22\approx0.361$. Here, $m_{G}+1=13$.
Using formula (33), we compute\begin{eqnarray*}&&x_{1}=12.13266641,
x_{2}=12.1331895,
 x_{3}=12.13363853, x_{4}=12.13402818, x_{5}=12.13436951,\\
&&x_{6}=12.13467089, x_{7}=12.13493919, x_{8}=12.135117936,
x_{9}=12.13539566, x_{10}=12.13559149,\\
&&x_{11}=12.13576962, x_{12}=12.13593235, x_{13}=12.13608195.
\end{eqnarray*}
The weights are assigned as follows:
$$w_{1}=w_{2}=\cdots=w_{12}=1,~
w_{13}=g_{m}=10\sqrt{5}-22.$$ Using formula (36), we obtain the
approximate solution (section basis) for the 20-dimensional
harmonic-type golden section equation:
$$\sum_{i=1}^{20}\frac{1}{x+i-1}=1.$$
The calculated partition basis is:
$$
\hat{\omega}=\frac{\sqrt{5}+1}{40}\Big[\sum_{i=1}^{12}x_{i}
+(10\sqrt{5}-22)x_{13}\Big]=12.13465705.
$$
From this, we obtain the golden section points $\hat{t}_{i}$ and
section intervals $\hat{\Delta}_{i}$ (for $i=1, 2, \cdots, 20$),
along with the balance ratio $\hat{\alpha}=1.6067328525$. When
comparing these results ($\hat{\omega}, \hat{\alpha}$) with the
high-precision solutions obtained via the binary search algorithm
($\omega=12.134641511$, $\alpha=1.6067320756$), we observe the
following errors: Absolute errors: $E_{\omega}=1.554\times10^{-5}$
for $\omega$, $E_{\alpha}=7.764\times10^{-7}$ for $\alpha$. The
relative errors are identical for both quantities.

Interestingly, if we take the first 12 or the first 13 terms of the
sequence $\{x_{1}, x_{2}, \cdots, x_{20}\}$ to compute the
arithmetic mean, we obtain $$
\hat{\omega}^{(12)}=\frac{1}{12}\sum_{i=1}^{12}x_{i}=12.13461423,~
\hbox{or}~~
\hat{\omega}^{(13)}=\frac{1}{13}\sum_{i=1}^{12}x_{i}=12.13472713.
$$
The corresponding errors are $E_{\omega}=-2.728\times10^{-5}$ or
$E_{\omega}=8.562\times10^{-5},$ respectively. This perfectly
illustrates the saying ``a hair's breadth less would be deficient, a
hair's breadth more would be excessive."  Moreover, this phenomenon
is not isolated. The pattern holds unless the fractional part of
$\frac{\sqrt{5}-1}{2}m$ is extremely small while the root of the
equation is very large. For example, in Table 2, when $m=500,$ the
fractional part of $250\sqrt{5}-250$ is merely 0.017, which makes
the mantissa weighting negligible for the large root of 291.488.
This demonstrates that the Relaxed Solutions Golden-Section Mantissa
Weighted Method fits well with the harmonic-type higher-dimensional
golden section equations.
$$\footnotesize
\begin{tabular}{c|ccccccc}
\multicolumn{7}{c}{\bf Table 2:  Software solutions ($\omega$,
$\alpha$), ~mantissa weighted method($\hat{\omega}$, $\hat{\alpha}$)
; errors $E_{\omega}$, $E_{\alpha}$}\\
\hline
$m$            & 4         &5          &6          &7           &8          &9          &10   \\\hline
$\omega$       &2.803886359&3.390541233&3.975681966&4.559938135 &5.143633852&5.726952553&6.310005691\\
$\hat{\omega}$ &2.803733108&3.390416647&3.975655994&4.559934223 &5.143626836&5.726968103&6.310018790\\
$E_{\omega}$   &$-1.54\times10^{-4}$&$-1.25\times10^{-4}$&$-2.59\times10^{-5}$&$-3.91\times10^{-6}$
&$-7.01\times10^{-6}$&$1.56\times10^{-5}$&$1.31\times10^{-5}$\\
$\alpha$       &1.700971590&1.678108247&1.662613661&1.651419734 &1.642954232&1.636328061&1.631000569\\
$\hat{\alpha}$ &1.700933277&1.678083329&1.662609137&1.651419175 &1.642953355&1.636329789&1.631001879\\
$E_{\alpha}$   &$-3.83\times10^{-5}$&$-2.49\times10^{-5}$&$-4.52\times10^{-6}$&$-5.59\times10^{-7}$
&$-8.77\times10^{-7}$&$1.73\times10^{-6}$&$1.31\times10^{-6}$\\
\hline
$m$            & 15        &20          & 30         &   40       &    50      &   60       & 70         \\\hline
$\omega$       &9.223131264&12.134641511&17.956037962&23.776620443&29.596876934&35.416970342&41.236970532\\
$\hat{\omega}$ &9.223149045&12.134657045&17.956049228&23.776629091&29.596883893&35.416976190&41.236975587\\
$E_{\omega}$&$1.78\times10^{-5}$&$1.55\times10^{-5}$&$1.13\times10^{-5}$&$8.65\times10^{-6}$&$6.97\times10^{-6}$
&$5.85\times10^{-6}$&$5.06\times10^{-6}$\\
$\alpha$       &1.614875418&1.606732076 &1.598534599 &1.594415511 &1.591937539 &1.590282839&1.589099579\\
$\hat{\alpha}$ &1.614876603&1.606732852 &1.598534974 &1.594415727 &1.591937678 &1.590282937&1.589099651\\
$E_{\alpha}$   &$1.19\times10^{-6}$&$7.76\times10^{-7}$&$3.75\times10^{-7}$&$2.16\times10^{-7}$&$1.39\times10^{-7}$
&$9.80\times10^{-8}$&$7.20\times10^{-8}$\\
\hline
$m$            & 80         &90          &100         &200          &300          &400          &500 \\\hline
$\omega$       &47.056912450&52.876815518&58.696691384&116.894851693&175.092685600&233.290437945&291.488157569\\
$\hat{\omega}$ &47.056916891&52.876819467&58.696694952&116.894853501&175.092686802&233.290438796&291.488158312\\
$E_{\omega}$&$4.44\times10^{-6}$&$3.95\times10^{-6}$&$3.57\times10^{-6}$&$1.81\times10^{-6}$&$1.20\times10^{-6}$
&$8.55\times10^{-7}$&$1.31\times10^{-7}$\\
$\alpha$       &1.588211406 &1.587520172 &1.586966914 &1.584474260  &1.583642285  &1.583226095  &1.582976315\\
$\hat{\alpha}$ &1.588211461 &1.587520216 &1.586966949 &1.584474268  &1.583642289  &1.583226097  &1.582976317\\
$E_{\alpha}$   &$5.50\times10^{-8}$&$4.40\times10^{-8}$&$3.50\times10^{-8}$&$8.00\times10^{-9}$&$4.00\times10^{-9}$
&$2.00\times10^{-9}$&$2.00\times10^{-9}$\\
\hline
\end{tabular}$$

Table 2 shows the values of $m$ and the corresponding solutions for
the $m$-dimensional harmonic-type golden section equation (6). It
includes the bisection search algorithm software solutions for
$\omega$ and $\alpha$, the approximate solutions $\hat{\omega}$ and
$\hat{\alpha}$, obtained by the Mantissa Weighted Method, as well as
the errors $E_{\omega}=\hat{\omega}-\omega$ and
$E_{\alpha}=\hat{\alpha}-\alpha$.

\vspace{1.5mm} \noindent{\bf Remark 8.1}~ In this section, we employ
natural logarithms of $m$ and $i$ to define relaxed time points,
narrowing the range of the largest positive real root (the partition
basis $\omega$) of the $m$-dimensional harmonic-type golden section
equation to: $\big(\frac{m}{e-1}+\frac{1}{2}-\frac{e-1}{4m+2(e-1)},~
\frac{m}{e-1}+\frac{1}{2}\big).$ Moreover, as $m$ increases,
$\omega$ gets closer to the right endpoint
$\frac{m}{e-1}+\frac{1}{2}$ of the interval, while the balance ratio
$\alpha$ gets closer to the left endpoint
$\frac{e}{e-1}+\frac{1}{2m}-\frac{e-1}{2m(2m+e-1)}$ of the interval
$\big(\frac{e}{e-1}+\frac{1}{2m}-\frac{e-1}{2m(2m+e-1)},~
\frac{e}{e-1}+\frac{1}{2m}\big).$ Examples and the data in Table 2
show that the approximate formulas for the roots of the
harmonic-type higher-dimensional golden section equations obtained
by the Relaxed Solutions Golden-Section Mantissa Weighted Method fit
the true solutions very closely. The method yields small errors and
high approximation accuracy when solving harmonic-type
higher-dimensional golden section equations. It also indicates that
our generalized harmonic-type higher-dimensional golden section is
closely related to the classical golden section.

\section{ Conclusions and Future Work}

In summary,  the core content and main contributions of this paper
are as follows:

1)~ {\it Definition and Modeling of Higher-Dimensional Golden
Sections}~ The classical golden section divides a unit line segment
into two parts. This paper defines and discusses the geometric-type
higher-dimensional golden section $\tau^{(m)}$ and the harmonic-type
higher-dimensional golden section $T^{(m)}$, which divide a unit
line segment into $m$ parts. From the definitions, we derive the
equation for the geometric-type higher-dimensional golden section
ratio $\lambda$:
$$\frac{1}{\lambda}+\frac{1}{\lambda^{2}}+\cdots+\frac{1}{\lambda^{m}}=1$$
and the equation for the harmonic-type higher-dimensional golden
section basis $\omega$:
$$\frac{1}{\omega}+\frac{1}{\omega+1}+\cdots+
\frac{1}{\omega+m-2}+\frac{1}{\omega+m-1}=1.$$ When $m=2$, both
cases reduce to the classical (two-dimensional) golden section. Our
generalization aims to expand the application space of this
classical problem.

2)~{\it Connection Between Generalized Fibonacci Sequences and
Higher-dimensional Golden Section}~ To better characterize the
properties of higher-dimensional golden sections, especially the
harmonic-type, we extend the $k$-generalized Fibonacci sequence
${F^{(k)}_{n}}_{n=0}^{+\infty}$ from existing literature to a
bidirectionally infinite $m$-generalized Fibonacci sequence
${F^{(m)}_{n}}_{n=0}^{\pm\infty}$ and define the bidirectionally
infinite generating sequence ${T^{(m)}_{n}}_{n=0}^{\pm\infty}$ for
harmonic-type higher-dimensional golden section points. This
establishes a crucial connection between the harmonic-type
higher-dimensional golden section $T^{(m)}$ and the $m$-generalized
Fibonacci sequence ${F^{(m)}_{n}}_{n=0}^{\pm\infty}$:
$$T^{(m)}_{n}=F^{(m)}_{n+m-1}t_{m}+F^{(m)}_{n+m-2}t_{m-1}+F^{(m)}_{n+m-3}t_{m-2}+
 \cdots+F^{(m)}_{n+1}t_{2}+F^{(m)}_{n}t_{1}.$$
Furthermore, we introduce the concept of the Fibonacci matrix for
harmonic-type higher-dimensional golden section points. Considering
the compatibility between determinant properties and the recurrence
relations of (generalized) Fibonacci sequences, we discuss the
determinants of submatrices of the Fibonacci matrix and obtain
significant results (e.g.,
$det\textbf{H}^{(m)}_{n}(x,y)=(-1)^{(n+h(xy))(m-1)}F^{(m)}_{x}F^{(m)}_{y}$).
When applied to the two-dimensional case, these results provide
concise and elegant proofs for certain key properties of the
classical Fibonacci sequence. Although the Fibonacci matrix is a
newly defined concept, our preliminary findings offer fresh
perspectives and approaches for studying (generalized) Fibonacci
sequences or higher-dimensional golden sections, with potential
academic value.

3)~{\it Mathematical Model and Optimal Strategy for Multiple Online
Leasing Problems}~ Another contribution of this paper is extending
the classical online leasing for a single pair of skis to the
multiple online leasing for $m$ pairs of skis. We design control and
balancing strategies for multiple online leasing problem and prove
their optimality using competitive analysis for online problems. The
optimal competitive ratio $\alpha$ is determined by the equation:
$$\frac{1}{m(\alpha-1)}+\frac{1}{m(\alpha-1)+1}+\cdots+\frac{1}{m(\alpha-1)+m-2}+\frac{1}{m(\alpha-1)+m-1}=1$$
Clearly, this is identical to the equation for the harmonic-type
higher-dimensional golden section basis $\omega$, i.e., the
$m$-dimensional harmonic-type golden section equation (6), by noting
that $\omega = m(\alpha-1)$. Moreover, this optimal strategy
coincides with the $m$-dimensional harmonic-type golden section
method: the ski group purchases one pair of skis at each golden
section point $t_{1}, t_{2}, \cdots, t_{m-1}, t_{m} = s$ in the
interval $[0, s]$. This provides insights and foundations for the
application of higher-dimensional golden section in optimization and
broader fields.

4)~{\it Continuous Relaxation Problem for Multiple Online Leasing
and the Relaxed Solutions Golden-Section Mantissa Weighted Method}~
Since the $m$-dimensional harmonic-type golden section equation (6)
for dimensions 5 and higher may not have an explicit solution, we
consider the continuous relaxation problem of the multiple online
leasing and its optimal strategy. We explore the feasibility of
approximating the optimal solution for multiple online leasing or
the largest real root of the harmonic-type higher-dimensional golden
section equation using transformed relaxation solutions. Numerical
and approximate analysis show that directly using the relaxed
solution as an approximation to the optimal solution does not yield
fully satisfactory outcomes. Fortunately, we discovered a
method---the relaxed solutions golden-section mantissa weighted
method: take the first $\frac{\sqrt{5}-1}{2}m$ relaxation solutions
and compute their weighted average, i.e.,
\begin{eqnarray*}
\hat{\omega}=\frac{\sqrt{5}+1}{2m}\Big[\sum_{i=1}^{m_{G}}x_{i} +
g_{m}x_{m_{G}+1}\Big].
\end{eqnarray*}
Interestingly, averaging the first
$\lfloor\frac{\sqrt{5}-1}{2}m\rfloor$ or
$\lceil\frac{\sqrt{5}-1}{2}m\rceil$ relaxation solutions typically
increases the error. This approximation formula may provide a
reference and new pathway for solving algebraic equations similar to
the harmonic-type higher-dimensional golden section equation (6).
Error analysis shows that the mantissa weighting method achieves
high precision in solving harmonic-type higher-dimensional golden
section equations. For $m>5$, the absolute error between the
competitive ratio for multiple online leasing or the harmonic-type
higher-dimensional golden section basis and the true solution is on
the order of $10^{-5}$, with relative errors below 6.515 parts per
million. Moreover, the error decreases further as $m$ increases.

The golden section and Fibonacci sequences are symbols and stepping
stones of human knowledge and civilization, while online problems
and competitive algorithms are key research topics in contemporary
theoretical computer science, artificial intelligence, and robotics.
Online leasing is one of the most classical problems in online
problems and competitive algorithms. The higher-dimensional golden
section proposed in this paper can advance research in robot search
algorithms and optimal decision-making for online leasing.
Conversely, the study of multiple online leasing problem offers new
pathways for (approximate) explicit solutions to algebraic equations
like the harmonic-type golden section equation (6). Extending the
golden section, Fibonacci sequences, and online leasing to
higher-dimensional and multiple scenarios holds both academic value
and practical potential.

This paper leaves several open questions for future research:

$\bullet$ Empirical studies on the biological properties of
high-dimensional golden sections.

$\bullet$ Investigation of architectural and artistic
characteristics of high-dimensional golden sections and their
applications.

$\bullet$ Applications of higher-dimensional golden sections in
optimization and competitive algorithms.

$\bullet$ General properties of Fibonacci matrices and their
submatrices for harmonic-type higher-dimensional golden section
points.

$\bullet$ Explicit solutions or higher-precision approximate
solutions for harmonic-type higher-dimensional golden section
equations.

$\bullet$ Optimal competitive algorithms for multi-device online
leasing problems and their relationship with harmonic-type
higher-dimensional golden sections.

%%%%%%%%%%%%%%%%%%%%%References%%%%%%%%%%%%%%%%%%%%%%%%%%%%%%%%%%%%%%%%

\end{document}